\documentclass[11pt,a4paper]{article}
\usepackage[T1]{fontenc}
\usepackage[utf8]{inputenc}
\usepackage{lmodern,microtype}
\usepackage[margin=27mm]{geometry}
\usepackage{amsmath,amssymb,amsthm,mathtools,mathrsfs}
\usepackage{enumitem,array,longtable,booktabs,xcolor}
\usepackage{xurl}
\usepackage[colorlinks=true,linkcolor=blue!50!black,citecolor=blue!50!black,urlcolor=blue!50!black]{hyperref}
\hypersetup{pdftitle={An Ito--Michler theorem and the prime graph of conjugacy-class dimensions in braided fusion categories},pdfdisplaydoctitle=true}
\allowdisplaybreaks[2]
\numberwithin{equation}{section}
\newtheorem{theorem}{Theorem}[section]
\newtheorem{proposition}[theorem]{Proposition}
\newtheorem{lemma}[theorem]{Lemma}
\newtheorem{corollary}[theorem]{Corollary}
\theoremstyle{definition}

\newtheorem{problem}[theorem]{Problem}
\theoremstyle{remark}
\newtheorem{remark}[theorem]{Remark}
\newtheorem{example}[theorem]{Example}
\newcommand{\C}{\mathcal C}
\newcommand{\D}{\mathcal D}
\newcommand{\E}{\mathcal E}
\newcommand{\M}{\mathcal M}
\newcommand{\Z}{\mathcal Z}
\newcommand{\one}{\mathbf1}
\providecommand{\Vec}{}
\renewcommand{\Vec}{\mathrm{Vec}}
\newcommand{\sVec}{\mathrm{sVec}}
\newcommand{\Rep}{\operatorname{Rep}}
\newcommand{\FPdim}{\operatorname{FPdim}}
\newcommand{\Irr}{\operatorname{Irr}}
\newcommand{\Hom}{\operatorname{Hom}}
\newcommand{\End}{\operatorname{End}}

\newcommand{\ad}{\mathrm{ad}}
\newcommand{\pt}{\mathrm{pt}}
\newcommand{\rev}{\mathrm{rev}}

\newcommand{\Dcl}{\Delta_{\mathrm{cl}}}

\newcommand{\kk}{\mathbb C}

\makeatletter
\newcommand{\eqlabelalias}[2]{\begingroup
  \edef\@currentlabel{\theequation}\edef\@currentHref{equation.\theHequation}%
  \label{#1}\endgroup}
\makeatother
\title{An It\^o--Michler  type theorem and the prime  class graph of weakly integral braided fusion categories}
\makeatletter
\renewcommand{\@maketitle}{%
  \begin{center}\LARGE\@title\par
  \ifx\@author\@empty\else\vskip .9em{\large\@author\par}\fi
  \end{center}\vskip .6em}
\makeatother
\author{Sebastian Burciu \thanks{This work was supported by a grant of the Ministry of Research, Innovation and Digitization, CNCS -UEFISCDI, project number PN-IV-P1-PCE-2023-0264, within PNCDI IV 
}} 
\date{\today}
\begin{document}
\maketitle
\begin{abstract}
Let $\C$ be a weakly integral braided fusion category. A prime $p$ is
absent from all conjugacy-class dimensions if and only if the Sylow
$p$-subgroup of the invertible objects generates a braided Deligne factor
of dimension equal to the full $p$-part of $\FPdim(\C)$. This condition
is also equivalent to $p\nmid\FPdim(\C_{\ad})$. We prove that the prime graph of the conjugacy-class dimensions of $\mathcal C$ has at most two connected components and that its disconnection implies solvability. The bound on the number of connected components is sharp for the Drinfeld double of $S_3$. 
\end{abstract}

\section{Introduction}
\label{sec:graph-introduction}\label{sec:im-introduction-20261024}

Let $G$ be a finite group and $p$ a prime divisor of $|G|$. The
It\^o--Michler theorem characterizes the absence of $p$ from the degrees
of irreducible characters of $G$ by the existence of a normal abelian
Sylow $p$-subgroup \cite{Ito1951,Michler1986}. For conjugacy-class sizes
the criterion has a different group-theoretic conclusion
\cite[Lemma~2.1]{CasoloTombari2012}.
\begin{equation}\label{explicit-alpha:a45450d4ea:eq0002}
 p\nmid |g^G|\text{ for every }g\in G
 \quad\Longleftrightarrow\quad
 \text{a Sylow $p$-subgroup of $G$ is central}.
\end{equation}
Such a Sylow subgroup is a direct factor. Our Theorem~\ref{thm:main} gives a categorical version of 
\eqref{explicit-alpha:a45450d4ea:eq0002}. Theorem~\ref{thm:unified-components-solvability-20261102}
bounds the number of connected components formed by the primes that can occur in the prime class graph.

Let $\C$ be a weakly integral braided fusion category over $\mathbb C$.
Its conjugacy classes are by definition the simple summands $Z_i$ of the canonical
algebra $I(\one)$ in its Drinfeld center, where $I$ is right adjoint to
the forgetful functor $F\colon Z(\C)\to\C$. Their Frobenius--Perron
dimensions are positive integers. We denote by $\Dcl(\C)$ the graph
whose vertices are the primes dividing some $\FPdim(Z_i)$, with an edge
between $p$ and $q$ when $pq$ divides a single $\FPdim(Z_i)$.
For $\C=\Rep(G)$ the class dimensions are the ordinary class sizes $|g^G|$. As usually, we write $\C_{\ad}$ for the subcategory generated by the objects
$X\otimes X^*$ with $X$ simple, and $\C_{\pt}$ for the subcategory of
invertible objects and their direct sums.

\begin{theorem}[Class It\^o--Michler theorem]\label{thm:main}
Let $\C$ be a weakly integral braided fusion category over $\kk$. Let $Z_i$ range over
its conjugacy-class objects in $Z(\C)$, and let $\C_{\ad}$ denote its
adjoint subcategory. Write $\FPdim(\C)=p^a N$ with $p$ prime,
$a\ge1$ and $N$ a positive integer coprime to $p$.
The following conditions are equivalent.
\begin{enumerate}[label=\textup{(\roman*)},leftmargin=*]
\item $p\nmid\FPdim(Z_i)$ for every conjugacy class of $\C$.
\item There is a braided equivalence
\begin{equation}\label{explicit:80658146e8:eq0011}
 \C\simeq\mathcal P\boxtimes\mathcal B,\qquad
 \mathcal P\text{ pointed},\quad
 \FPdim(\mathcal P)=p^a,\quad\FPdim(\mathcal B)=N.
\end{equation}\eqlabelalias{explicit-alpha:a45450d4ea:eq0007}{explicit:80658146e8:eq0011}
\item $p\nmid\FPdim(\C_{\ad})$.
\end{enumerate}
Here $\mathcal P$ is generated by the Sylow $p$-subgroup of
$\operatorname{Inv}(\C)$, and $\mathcal B$ is the sum of the universal
grading components whose degrees have order prime to $p$. In particular,
\begin{equation}\label{eq:unified-vertices-20261102}
 \operatorname{Vert}\Dcl(\C)=\pi\bigl(\FPdim(\C_{\ad})\bigr),
\end{equation}
where $\pi(n)$ denotes the prime divisors of $n$.
\end{theorem}

Under hypothesis \textup{(i)}, one already has
$p^a\mid\FPdim(\C_{\pt})$ by \cite[Theorem~1.2]{JIMJ25}.
The modular pointed splitting and its iteration were given in
\cite[Proposition~1.3 and Theorem~1.4]{JIMJ25} in terms of simple
dimensions. In the weakly integral modular case, the prime divisors
must be taken from the integer squares $\FPdim(X)^2$, as explained in
Section~\ref{sec:im-checks-20261024}. Theorem~\ref{thm:main} gives the
braided splitting for possibly degenerate $\C$, specifies its
complement, and characterizes it by the adjoint subcategory.
The pointed factor can be degenerate. If $\C$ is nondegenerate, both
factors are nondegenerate.

The \emph{conjugacy-class prime graph} $\Dcl(\C)$ has vertex set
$\bigcup_i \pi(\FPdim(Z_i))$. Two distinct vertices $p$ and $q$ are adjacent
if and only if $pq\mid\FPdim(Z_i)$ for some conjugacy class $Z_i$.

\begin{theorem}[Components and solvability]\label{thm:unified-components-solvability-20261102}
Let $\C$ be a weakly integral braided fusion category over $\kk$. Then
$\Dcl(\C)$ has at most two connected components. If it has two, then $\C$ is solvable.
\end{theorem}

We note that this bound is sharp. In $\Rep(D(S_3))$, the nontrivial class dimensions
are $4$ and $9$, so the graph consists of the two isolated vertices $2$
and $3$ (Example~\ref{ex:class-graph-double-s3-20261102}).

The history of the finite-group bound involves two graphs. In 1981
Kazarin classified groups admitting a pair of isolated conjugacy classes
\cite[main theorem]{Kazarin}. In 1990 Bertram, Herzog and Mann proved
that the graph on noncentral conjugacy classes, with adjacency defined
by a common prime divisor of their sizes, has at most two components
\cite{BHM}. This also bounds the number of components of the prime
graph. Indeed, both are obtained from the incidence graph between primes and
noncentral classes by joining vertices on one side with a common
neighbor. Alfandary gave the result directly for the graph on primes
in 1994 and characterized the disconnected case in terms of a Frobenius
quotient \cite{Alfandary}. Theorem~\ref{thm:unified-components-solvability-20261102}
extends the component bound from $\Rep(G)$ to arbitrary weakly
integral braided fusion categories.

For comparison, in \cite{NPR}, Natale and Rodr\'iguez studied the graph built
from the dimensions of \emph{simple objects}. For integral nondegenerate
braided categories their Theorem~1.1 gives at most three components
in general and at most two under solvability. The bound of two also
holds when all Tannakian subcategories are pointed
\cite[Remark~2.2 and Proposition~10.4]{NPR} and for every
group-theoretical modular category \cite[Theorem~1.2(ii)]{NPR}.
For a modular category the class dimensions are precisely
$\FPdim(X)^2$, $X\in\Irr(\C)$ (Lemma~\ref{lem:class-square-minimal-center-20261026}\textup{(1)}). In the integral modular case, squaring
preserves every prime support, so these two graphs coincide. Theorem~\ref{thm:unified-components-solvability-20261102} therefore gives
the bound of two for the simple-dimension graph of every integral
modular category. For strictly weakly integral modular categories
the graph here uses the integer squares $\FPdim(X)^2$.

The solvability assertion also has a precursor in the group-theoretical
modular case. The proof of \cite[Theorem~1.2(ii)]{NPR} realizes
$\C$ as a braided factor of a twisted double $\Rep(D^\omega(G))$ with
pointed complement, so the two categories have the same prime graph.
When that graph is disconnected, $G/Z(G)$ is Frobenius with abelian
kernel and complement \cite[Theorem~8.5(iii)]{NPR}. Thus $G$ is
solvable. Solvability passes from $\Vec_G^\omega$ to its center and
then to $\C$ by \cite[Proposition~4.5(i),(ii)]{ENO}.
Theorem~\ref{thm:unified-components-solvability-20261102} obtains
solvability from disconnectedness for every weakly integral braided
fusion category.

Two properties of the class objects enter the graph arguments.
First, Proposition~\ref{prop:diagonal} identifies the underlying class
objects $Y_i=F(Z_i)$ as internal endomorphism objects, and
Proposition~\ref{prop:positivity} proves that their character values are
nonnegative. Using orthogonality, Theorem~\ref{thm:productintro} gives the
regular-product formula
\begin{equation}\label{eq:unified-intro-regular-product-20261102}
 \prod_i[Y_i]=
 \frac{\prod_i\FPdim(Z_i)}{\FPdim(\C_{\ad})}
 \sum_{X\in\Irr(\C_{\ad})}\FPdim(X)[X].
\end{equation}
The coefficient of the tensor unit gives the divisibility used in
Theorem~\ref{thm:main}. In the modular case this recovers
\cite[Theorem~1.16]{BP}. Second, Theorem~\ref{thm:braided-divisibility} gives the individual
divisibility~\eqref{eq:braided-divisibility} under restriction to a braided
fusion subcategory. Under equivariantization of a
non-pointed base, every component of the new graph contains a prime from
the base graph (Theorem~\ref{thm:class-graph-connected-equivariantization-20261027}).
The pointed case follows from the character-fiber formula~\eqref{explicit:80658146e8:eq0117}
of Proposition~\ref{prop:equivariant-class-fiber-formula-20261027}.

The proof of Theorem~\ref{thm:unified-components-solvability-20261102}
has two parts. Theorem~\ref{thm:class-graph-modular-two-20261027} establishes the bound
for modular categories. Proposition~\ref{prop:class-graph-modular-realization-20261027}
constructs a modular extension with the same class graph for each weakly
integral slightly degenerate category. De-equivariantization by the even part of the
M\"uger center reduces the general bound to these two cases
(Theorem~\ref{thm:general-two-components-20261102}). For solvability,
Proposition~\ref{prop:class-graph-minimal-counterexample-wgt-20261027} uses
centralizer descent to reduce a hypothetical minimal counterexample to
a weakly group-theoretical category, where the structure of its core
gives the contradiction
(Theorem~\ref{thm:disconnected-class-graph-braided-solvable-20261027}).
The stronger assertion that class characters vanish across components
holds for categories tensor equivalent to $\Rep(G)$ and for untwisted
Drinfeld doubles (Propositions~\ref{prop:symmetric-cross-component-vanishing-20261102}
and~\ref{prop:untwisted-double-cross-component-vanishing-20261102}).
In the modular case it also holds when every Tannakian subcategory
is pointed, by \eqref{eq:modular-cross-component-class-value-20261102}
in the proof of Theorem~\ref{thm:class-graph-modular-two-20261027}.
The general vanishing question remains open.

Section~\ref{sec:prelim} proves that character values on class images
are nonnegative and that class dimensions divide under restriction.
Section~\ref{sec:product} derives the Sylow splitting of Theorem~\ref{thm:main}
from the regular-product identity~\eqref{eq:productintro} of Theorem~\ref{thm:productintro}. Section~\ref{sec:unified-prime-supports}
treats products and gradings and proves solvability for prime-power
class dimensions (Theorem~\ref{thm:class-prime-powers-solvable-20261026}).
Section~\ref{sec:unified-tannakian} uses coprime class dimensions to prove
the two-component bound when every Tannakian subcategory is pointed
(Theorem~\ref{thm:class-graph-pointed-tannakians-20261027}).
Section~\ref{sec:unified-modular} treats modular categories, slightly degenerate
categories, and categories with cyclic Tannakian M\"uger center. In the
slightly degenerate case, the two fermion quotients give all class dimensions (Proposition~\ref{prop:fermion-sign-class-dimensions}).
Section~\ref{sec:unified-equivar} proves the general two-component bound
by equivariantization (Theorem~\ref{thm:general-two-components-20261102}),
and Section~\ref{sec:unified-solvability} deduces solvability from
disconnectedness by centralizer descent
(Theorem~\ref{thm:disconnected-class-graph-braided-solvable-20261027}).
Section~\ref{sec:component-bound-vanishing-20261102} studies vanishing across
components for finite-group representation categories and untwisted
Drinfeld doubles. Appendix~\ref{sec:general-relative} uses the relative
vanishing theorem (Theorem~\ref{thm:general-relative}) to bound the total
weight of characters vanishing on a class image (Corollary~\ref{cor:image-zero}).

\section{Class objects, formal codegrees and restriction}\label{sec:prelim}\label{app:class-graph-inputs-20261102}
For a fusion category $\C$, write $\Irr(\C)$ for its simple isomorphism
classes and $[V:X]=\dim\Hom(X,V)$ for the multiplicity of a simple
$X$ in $V$. The pointed subcategory $\C_{\pt}$ is generated by the invertible
simple objects, equivalently those with dimension one. Their
isomorphism classes form a group $\operatorname{Inv}(\C)$, which is abelian when
$\C$ is braided. The adjoint subcategory $\C_{\ad}$ is generated by
$X\otimes X^*$ for $X\in\Irr(\C)$.
The universal grading is the faithful grading
\begin{equation}\label{explicit-alpha:a45450d4ea:eq0011}
 \C=\bigoplus_{g\in U(\C)}\C_g,\qquad
 \C_e=\C_{\ad},\qquad
 \FPdim(\C)=|U(\C)|\FPdim(\C_{\ad}).
\end{equation}\eqlabelalias{eq:gradingdim}{explicit-alpha:a45450d4ea:eq0011}
For weakly integral $\C$, the subcategory $\C_{\ad}$ is integral
\cite{GN}. The Deligne tensor product is denoted by $\boxtimes$;
its simple objects are external tensor products and its dimension is
the product of the dimensions of its factors.

Recall that a fusion category is \emph{nilpotent}
if repeated passage to the adjoint subcategory eventually gives $\Vec$.
It is \emph{weakly group-theoretical} if it is categorically Morita
equivalent to a nilpotent fusion category, and \emph{solvable} if it is
Morita equivalent to one obtained from $\Vec$ by successive extensions
by cyclic groups of prime order. Here categorical Morita equivalence
means equivalence to the dual category of an indecomposable semisimple
module category; see~\cite{ENO}.

All categories are over $\kk=\mathbb C$. A fusion category $\mathcal B$
is weakly integral if $\FPdim(\mathcal B)$ is an integer, and integral
if every simple object has integer Frobenius--Perron dimension.
A braided category is nondegenerate if its M\"uger center is $\Vec$,
slightly degenerate if its M\"uger center is braided equivalent to $\sVec$,
and symmetric if its M\"uger center is the entire category. For
weakly integral categories we use the canonical positive spherical
structure, under which a nondegenerate braided fusion category is
modular. Unless otherwise specified, $\C$ is weakly integral and
braided. Write
\begin{equation}\label{explicit:80658146e8:eq0001}
 F:Z(\C)\longrightarrow\C,\qquad
 I:\C\longrightarrow Z(\C),
\end{equation}\eqlabelalias{explicit-alpha:a45450d4ea:eq0004}{explicit:80658146e8:eq0001}
where $I$ is the right adjoint of the forgetful functor $F$.
Since $K_0(\C)$ is commutative, the canonical algebra is
multiplicity-free \cite[Theorem~2.13]{OstrikRankThree}. It decomposes as
\begin{equation}\label{explicit:80658146e8:eq0002}
 \mathscr A=I(\one)=\bigoplus_{i\in J}Z_i.
\end{equation}\eqlabelalias{explicit-alpha:a45450d4ea:eq0005}{explicit:80658146e8:eq0002}
The $Z_i$ are the conjugacy-class objects. Put
\begin{equation}\label{explicit:80658146e8:eq0003}
 Y_i=F(Z_i),\qquad \FPdim(Z_i)=\FPdim(Y_i),\qquad
 f_i=\frac {\FPdim(\C)}{\FPdim(Z_{i})}.
\end{equation}
Under weak integrality both $\FPdim(Z_i)$ and $f_i$ are positive integers. Here $f_i$ is the
formal codegree of the associated algebra character $\mu_i$ of
$K(\C)=K_0(\C)\otimes_{\mathbb Z}\kk$.
Here a character is a unital algebra homomorphism to $\mathbb C$.
For a braided fusion category $\mathcal B$ and a character $\nu$ of
$K(\mathcal B)$, write $Z_\nu^{\mathcal B}$ for the corresponding
conjugacy-class object in $\mathcal Z(\mathcal B)$. 
We use unnormalized $S$-matrices throughout. We denote the set of
connected components of a graph by $\pi_0$ and the M\"uger center by
$Z_2(\C)$ (also written $\C'$).

For $n\in\mathbb Z_{>0}$ let $\pi(n)$ be the set of prime divisors of $n$;
thus $\pi(1)=\varnothing$. The graph $\Dcl(\C)$ has vertex set
$\bigcup_i\pi(\FPdim(Z_{i}))$, with an edge between distinct primes $p,q$ if
$pq\mid \FPdim(Z_{i})$ for some $i$. A class dimension $\FPdim(Z_{i})>1$ is said to belong
to a component when its prime support does; its support is a clique and
hence lies in a unique component. We count zero components for the empty graph and call a graph disconnected when it has at least two nonempty components.
The common-divisor graph on the nontrivial class dimensions has the same
number of connected components, but not necessarily the same diameter.

Let $X^+=\iota_+(X)$ and $X^-=\iota_-(X)$ denote the two canonical lifts
to the center. Their images satisfy
\begin{equation}\label{explicit:80658146e8:eq0004}
 (\C^+)'=\C^-,\qquad (\C^-)'=\C^+.
\end{equation}\eqlabelalias{eq:mutual}{explicit:80658146e8:eq0004}
We write $N_{X,Y}^W=[X\otimes Y:W]$ for the fusion multiplicities.
For a full fusion subcategory $\mathcal B$ and an object $V$, the notation
$V\cap\mathcal B$ denotes the sum of the simple constituents of $V$
belonging to $\mathcal B$, with multiplicities.

For a commutative fusion ring, let $\tau$ denote the coefficient of
the unit. The primitive idempotents, the trace decomposition, and
character orthogonality are given by
\begin{equation}\label{eq:app-class-idempotents-20261102}
 \begin{gathered}
 e_\mu=\frac1{f_\mu}\sum_{X\in\Irr(\C)}\mu(X^*)[X],
 \qquad \tau(v)=\sum_\mu\frac{\mu(v)}{f_\mu},\\
 \sum_{X\in\Irr(\C)}\mu(X)\overline{\nu(X)}=\delta_{\mu,\nu} f_\mu,
 \qquad f_\mu=\sum_{X\in\Irr(\C)}|\mu(X)|^2 .
 \end{gathered}
\end{equation}\eqlabelalias{eq:idempotents}{eq:app-class-idempotents-20261102}\eqlabelalias{eq:traceexpansion}{eq:app-class-idempotents-20261102}
Indeed, the simple basis is orthonormal for $\tau(uv^*)$ and
multiplication by $[X^*]$ is the adjoint of multiplication by $[X]$.
Simultaneous diagonalization gives~\eqref{eq:app-class-idempotents-20261102}.

For $W\in\Irr(\Z(\C))$, write
$\chi_W(V)=S_{V,W}/\FPdim(W)$. The labelling theorem
\cite[Theorem~2.13]{OstrikRankThree} identifies the character
$\mu_i\circ F$ with $\chi_{Z_i}$ and gives
\begin{equation}\label{eq:app-master-20261102}
 S_{W,Z_i}=\FPdim(Z_i)\mu_i(F(W)),\qquad
 f_i=\frac{\FPdim(\C)}{\FPdim(Z_i)},\qquad
 [Y_i:\one]=[I(\one):Z_i]=1 .
\end{equation}\eqlabelalias{eq:master}{eq:app-master-20261102}

By the usual description of induction,
$\bigoplus_iY_i\simeq F(I(\one))\simeq
\bigoplus_{X\in\Irr(\C)}XX^*$. Thus the $Y_i$ lie in and generate
$\C_{\ad}$. If $\C$ is weakly integral, this subcategory is integral by the dimension grading
\cite[Theorem~3.10]{GN}. In particular, the class dimensions and
their formal codegrees are integers, the latter being rational
algebraic integers.

\subsection{Restricted modular orthogonality}\label{app:class-structure-20261102}
\begin{lemma}[Orthogonality for a fusion subcategory]\label{lem:orthog}

Let $\mathcal L$ be a pseudo-unitary nondegenerate braided fusion
category over $\kk$ with its canonical positive spherical structure.
Let $S$ be its unnormalized modular $S$-matrix,
$\mathcal T\subseteq\mathcal L$ a full fusion subcategory with the
inherited braiding, and $\mathcal T'_{\mathcal L}$ its M\"uger
centralizer in $\mathcal L$. For any $U,V\in\Irr(\mathcal L)$, one has
\begin{equation}\label{eq:app-centralizer-orthogonality-20261102}
 \sum_{T\in\Irr(\mathcal T)}S_{T,U}\overline{S_{T,V}}
 =\FPdim(\mathcal T)
   \sum_{W\in\Irr(\mathcal T'_{\mathcal L})}
          \FPdim(W) N_{W,V}^{U}.
\end{equation}\eqlabelalias{eq:orthog}{eq:app-centralizer-orthogonality-20261102}
\end{lemma}
\begin{proof}
Evaluate
$\sum_{T\in\Irr(\mathcal T)}\FPdim(T)[T]$ under $\chi_W$.
This regular element has value $\FPdim(\mathcal T)$ if
$\chi_W|_{\mathcal T}=\FPdim$ and zero otherwise. The former
condition is equivalent to $W\in\mathcal T'_{\mathcal L}$ by the
pseudo-unitary centralizer criterion~\cite[Section~3]{DGNO}.
The character property of $\chi_T$ gives
$S_{T,U}S_{T,V^*}=\FPdim(T)\sum_W N_{U,V^*}^{W}S_{T,W}$.
Sum over $T$, use $S_{T,V^*}=\overline{S_{T,V}}$ and apply
$N_{U,V^*}^{W}=N_{W,V}^{U}$.
\end{proof}

\subsection{Internal endomorphisms and positivity}\label{sec:positive}
\begin{proposition}[Diagonal multiplicities]\label{prop:diagonal}
Let $\C$ be a pseudo-unitary braided fusion category over $\kk$ with its canonical positive spherical structure. Let
$F\colon Z(\C)\to\C$ be the forgetful functor with right adjoint $I$,
let $\mathscr A=I(\one)=\bigoplus_{i\in J}Z_i$ be the decomposition
into simple class objects, and put $Y_i=F(Z_i)$. For a simple
$X\in\C$ let $X^+=\iota_+(X)$ be its central lift using the given
braiding of $\C$. Then for $i,j\in J$,
\begin{equation}\label{explicit:80658146e8:eq0005}
 N_{X^+,Z_j}^{Z_i}=0\quad(i\ne j),\qquad
 [Y_i:X]=N_{X^+,Z_i}^{Z_i}.
\end{equation}\eqlabelalias{eq:offdiagonal}{explicit:80658146e8:eq0005}
Thus $Y_i$ is the internal endomorphism object of $Z_i$ for the
$\C$-action on $Z(\C)$ through $\iota_+$, and carries a connected
separable algebra structure.
\end{proposition}
\begin{proof}
Take $\mathcal L=\Z(\C)$, $\mathcal T=\C^-$, $U=Z_i$ and $V=Z_j$
in the orthogonality relation~\eqref{eq:app-centralizer-orthogonality-20261102}
of Lemma~\ref{lem:orthog}.
Since $(\C^-)'=\C^+$ and $F(X^-)=X$,
\eqref{eq:app-master-20261102} and character orthogonality~\eqref{eq:traceexpansion} yield
\begin{equation}\label{eq:app-diagonal-sum-20261102}
 \sum_{X\in\Irr(\C)}\FPdim(X) N_{X^+,Z_j}^{Z_i}
 =\frac{\FPdim(Z_i)\FPdim(Z_j)}{\FPdim(\C)}
   \sum_{X\in\Irr(\C)}\mu_i(X)\overline{\mu_j(X)}
 =\delta_{ij}\FPdim(Z_i).
\end{equation}
Nonnegativity gives the off-diagonal vanishing. Adjunction and the
projection formula~\eqref{eq:projection-right-adjoint} give
$[Y_i:X]=[I(X):Z_i]=[X^+ I(\one):Z_i]=N_{X^+,Z_i}^{Z_i}$.
These multiplicities characterize $Y_i$ as the internal end of
$Z_i$ for the action through $\iota_+$.

We recall why these internal ends are connected separable algebras.
For a fusion subcategory $\mathcal E\subseteq\mathcal L$ acting on
$\mathcal L$ by left tensoring and a simple $W\in\mathcal L$, the
internal end of $W$ is $(WW^*)\cap\mathcal E$. Its multiplication
is induced by evaluation on $WW^*$, and its unit has multiplicity one
by rigidity. For $Q=\underline{\End}_{\mathcal E}(W)$, the category
of right $Q$-modules in $\mathcal E$ is the indecomposable semisimple
module component generated by $W$ \cite[Theorem~7.10.1]{EGNO}.
Its dual, the $Q$-bimodule category up to tensor reversal, is
semisimple \cite[Corollary~9.3.3]{EGNO}. Multiplication
$Q\otimes Q\to Q$ is an epimorphism of bimodules, so it splits.
This proves separability, both in $\mathcal E$ and in a containing
fusion category.
\end{proof}

\begin{proposition}[Positivity and support]\label{prop:positivity}
Let $\C$ be a pseudo-unitary braided fusion category over $\kk$ with its
canonical positive spherical structure. For the forgetful functor
$F\colon Z(\C)\to\C$ and its right adjoint $I$, write
$I(\one)=\bigoplus_i Z_i$ and $Y_i=F(Z_i)$. Let $\mu_i$ denote the
algebra character of $K(\C)$ corresponding to $Z_i$. For every
algebra character $\mu$ of $K(\C)$ and every $i$,
\begin{equation}\label{explicit:80658146e8:eq0006}
 \mu(Y_i)\ge0,\qquad [Y_i:\one]=1,\qquad Y_i\in\C_{\ad}.
\end{equation}\eqlabelalias{eq:positive}{explicit:80658146e8:eq0006}
The $Y_i$ are self-dual and generate $\C_{\ad}$. With $S$ the
unnormalized $S$-matrix of $Z(\C)$ and $X^+$ the lift of
$X\in\Irr(\C)$ through the given braiding, the center-character
normalizations are
\begin{equation}\label{explicit:80658146e8:eq0007}
 \frac{S_{Z_i,X^+}}{\FPdim(Z_{i})}=\mu_i(X),\qquad
 \frac{S_{Z_i,Z_j}}{\FPdim(Z_{i})}=\mu_i(Y_j).
\end{equation}
\end{proposition}
\begin{proof}
Let $\pi$ be the action by $X^+$ on the vector space with
orthonormal basis $v_W$, $W\in\Irr(\Z(\C))$. Rigidity gives
$\pi([X])^\dagger=\pi([X^*])$, so $\pi(e_\mu)$ is an orthogonal
projection. Its diagonal coefficient at $Z_i$ is, by the multiplicity
formula~\eqref{eq:offdiagonal} of Proposition~\ref{prop:diagonal},
\begin{equation}\label{eq:app-positive-projection-20261102}
 \frac{\mu(Y_i^*)}{f_\mu}
 =\frac1{f_\mu}\sum_X\mu(X^*)[Y_i:X]
 =\|\pi(e_\mu)v_{Z_i}\|^2 .
\end{equation}\eqlabelalias{eq:positive-norm}{eq:app-positive-projection-20261102}
Thus $\mu(Y_i^*)=\overline{\mu(Y_i)}$ is real and nonnegative.
The characters separate elements, giving $Y_i\simeq Y_i^*$ as claimed.

The character normalizations in \eqref{explicit:80658146e8:eq0007} follow from \eqref{eq:app-master-20261102} with $W=X^+$ and $W=Z_j$. This projection argument uses only
pseudo-unitarity.
\end{proof}

\subsection{Restriction of class characters}\label{sec:restriction}
\begin{proposition}[Extent of a character fiber]\label{prop:extent}
Let $\C$ be a pseudo-unitary braided fusion category over $\kk$ with its canonical positive spherical structure. Let
$\mathcal B\subseteq\C$ be a full fusion subcategory with the inherited
braiding. For each algebra character $\nu$ of $K(\mathcal B)$, summing
over the algebra characters $\mu$ of $K(\C)$ that restrict to $\nu$
gives
\begin{equation}\label{explicit:80658146e8:eq0008}
 \sum_{\mu|_{\mathcal B}=\nu}\FPdim(Z_{\mu}^{\C})
 =\frac{\FPdim(\C)}{\FPdim(\mathcal B)}\FPdim(Z_{\nu}^{\mathcal B}).
\end{equation}\eqlabelalias{eq:extent-general}{explicit:80658146e8:eq0008}
Every character of $K(\mathcal B)$ has an extension to $K(\C)$.
\end{proposition}
\begin{proof}
For a braided inclusion $\mathcal B\subseteq\C$, the image of a
primitive idempotent $e_\nu\in K(\mathcal B)$ is
$\sum_{\mu|_{\mathcal B}=\nu}e_\mu$. This sum is nonzero by
injectivity, so each character extends. The unit-coefficient traces
agree under inclusion and give
$1/f_\nu=\sum_{\mu|_{\mathcal B}=\nu}1/f_\mu$.
Multiplying by $\FPdim(\C)$ and using~\eqref{eq:app-master-20261102}
in both categories gives the formula.
\end{proof}

Let $C_\mu$ and $C_\nu^{\mathcal B}$ denote the class sums indexed
by these characters, and let $\pi_e$ be the map on central elements
dual to inclusion of the character algebras. Normalizing the class sums
to be dual to the primitive idempotents gives
\begin{equation}\label{eq:normalized-projection}
 \pi_e\!\left(\frac{C_\mu}{\FPdim(Z_\mu^{\C})}\right)
 =\frac{C_{\mu|_{\mathcal B}}^{\mathcal B}}
        {\FPdim(Z_{\mu|_{\mathcal B}}^{\mathcal B})}.
\end{equation}
Indeed, pairing each side with $e_\nu$ gives one if
$\mu|_{\mathcal B}=\nu$ and zero otherwise
\cite[Proposition~3.1 and Theorem~3.3]{JIMJ25}.
Write $\widehat{\Lambda}_{\C}=\sum_\mu C_\mu$ and
$\widehat{\Lambda}_{\mathcal B}=\sum_\nu C_\nu^{\mathcal B}$.
With the canonical positive spherical structure, the identity
$\pi_e(\widehat{\Lambda}_{\C})
=(\FPdim(\C)/\FPdim(\mathcal B))\,\widehat{\Lambda}_{\mathcal B}$
of \cite[(3.4)]{JLMS26}, together with \eqref{eq:normalized-projection},
gives, on comparing the coefficients of $C_\nu^{\mathcal B}$,
\begin{equation}\label{eq:im-fibre-normalized-mass-20261102}
 \sum_{\mu|_{\mathcal B}=\nu}
 \frac{\FPdim(Z_\mu^{\C})}{\FPdim(Z_\nu^{\mathcal B})}
 =\frac{\FPdim(\C)}{\FPdim(\mathcal B)} .
\end{equation}
Multiplying by $\FPdim(Z_\nu^{\mathcal B})$ gives
\eqref{eq:extent-general}. The missing factor $\FPdim(Z_\nu^{\mathcal B})$
on the right-hand side of \cite[Proposition~3.1, (3.3)]{JLMS26} is a
typographical omission in the published formula. The same omission
appears in its formal-order version \cite[Remark~3.2, (3.6)]{JLMS26}.
For the dimension character $\nu=\FPdim|_{K(\mathcal B)}$, the
downstairs class is the unit, so the unnormalized mass equals the index.
For the identity inclusion of $\Rep(S_3)$, the fiber over evaluation at a
transposition is a singleton of mass three and the index is one. For
$\Rep(S_3)\subset\Rep(S_3)\boxtimes\Vec_{C_2}$, the same downstairs
character has two extensions of class dimension three each. Dividing
their total mass by three recovers the index two.

In the modular case, Bantay's extent of a restriction class is the
reciprocal of the sum of the squared normalized vacuum-row entries in
that class \cite[Definition~3, (3.4)]{BantayFC}. In our notation these
entries are $\FPdim(X)/\sqrt{\FPdim(\C)}$, so the extent of the fiber
over $\nu$ is its downstairs formal codegree $f_\nu$ by
\eqref{eq:extent-general} and \eqref{eq:app-master-20261102}.

\begin{theorem}[Exact modular quotient]\label{thm:modular-quotient}
Let $\mathcal M$ be a pseudo-unitary nondegenerate braided fusion
category over $\kk$ with its canonical positive spherical structure,
and let $\mathcal B\subseteq\mathcal M$ be a full fusion subcategory
with the inherited braiding. For $W\in\Irr(\mathcal M)$ let $\nu$
be the restriction of its Verlinde character, given on
$X\in\Irr(\mathcal B)$ by $\nu(X)=S_{X,W}/\FPdim(W)$ for the
unnormalized $S$-matrix of $\mathcal M$. If
$\mathcal B'_{\mathcal M}$ denotes the M\"uger centralizer in
$\mathcal M$, then
\begin{equation}\label{explicit:80658146e8:eq0010}
 \frac{{\FPdim(W)}^2}{\FPdim(Z_{\nu}^{\mathcal B})}
 =\FPdim\bigl((WW^*)\cap\mathcal B'_{\mathcal M}\bigr)
 .
\end{equation}\eqlabelalias{eq:modular-quotient}{explicit:80658146e8:eq0010}
The quotient is a positive algebraic integer. If $\mathcal M$ is weakly
integral, it is a positive integer.
\end{theorem}
\begin{proof}
Let $\nu=\chi_W|_{\mathcal B}$ for $W\in\Irr(\M)$. Applying
Lemma~\ref{lem:orthog} with $U=V=W$ in~\eqref{eq:app-centralizer-orthogonality-20261102}
gives
\begin{equation}\label{eq:app-modular-quotient-20261102}
 \frac{\FPdim(W)^2}{\FPdim(Z_\nu^{\mathcal B})}
 =\frac{\FPdim(W)^2 f_\nu}{\FPdim(\mathcal B)}
 =\sum_{T\in\Irr(\mathcal B'_{\M})}\FPdim(T)N_{W,W^*}^{T}.
\end{equation}
The right-hand side is the Frobenius--Perron dimension of an object
containing the unit. Under weak integrality, each contributing $T$ belongs
to the integral subcategory $\M_{\ad}$.
\end{proof}

\begin{theorem}[Individual braided divisibility]\label{thm:braided-divisibility}
Let $\C$ be a pseudo-unitary braided fusion category over $\kk$ with
its canonical positive spherical structure. Let
$\mathcal B\subseteq\C$ be a full fusion subcategory with the inherited
braiding, and let $\mu$ be an algebra character of $K(\C)$.
The quotient $\FPdim(Z_\mu^{\C})/
\FPdim(Z_{\mu|_{\mathcal B}}^{\mathcal B})$ is a positive algebraic
integer. If $\C$ is weakly integral, this gives ordinary integer
divisibility
\begin{equation}\label{explicit:80658146e8:eq0009}
 \FPdim(Z_{\mu|_{\mathcal B}}^{\mathcal B})\mid \FPdim(Z_{\mu}^{\C}).
\end{equation}\eqlabelalias{eq:braided-divisibility}{explicit:80658146e8:eq0009}
In the weakly integral case $\Dcl(\mathcal B)$ is a subgraph of
$\Dcl(\C)$.
\end{theorem}
\begin{proof}
For $\mathcal B\subseteq\C$ take $W=Z_i\in\Z(\C)$ and put
$\nu=\mu_i|_{\mathcal B}$. In $\Z(\C)$ let
$\mathcal E=(\C^+)'\subseteq\mathcal L=(\mathcal B^+)'$.
Write $Q_{\mathcal E}=(WW^*)\cap\mathcal E$ and
$Q_{\mathcal L}=(WW^*)\cap\mathcal L$.
Each is the internal end of $W$ for the indicated subcategory.
Because both multiplications are restrictions of that on $WW^*$,
$Q_{\mathcal E}\subseteq Q_{\mathcal L}$ is a unital subalgebra.
By~\eqref{eq:app-master-20261102}, $\chi_W$ restricts to $\mu_i$
on $\C^+$, and hence to $\nu$ on $\mathcal B^+$. Applying
the modular quotient formula~\eqref{eq:modular-quotient}
of Theorem~\ref{thm:modular-quotient} twice gives
\begin{equation}\label{eq:app-internal-ratio-20261102}
 \FPdim(Q_{\mathcal E})=\FPdim(Z_i),\qquad
 \FPdim(Q_{\mathcal L})=\frac{\FPdim(Z_i)^2}{\FPdim(Z_\nu^{\mathcal B})}.
\end{equation}
Set $Q=Q_{\mathcal E}$. By the internal-end argument in the proof of
Proposition~\ref{prop:diagonal}, $Q$ is connected and separable, so its bimodules in $\mathcal L$ form a multifusion category. Its tensor unit $Q$ is simple because
$\End_{Q\text{-bimod}}(Q)\subseteq\Hom_{\mathcal L}(\one,Q)=\kk$
under evaluation at the unit. Thus the bimodules form a fusion category. In this category $Q_{\mathcal L}$ is
a bimodule by multiplication. The bimodule dimension normalization
\cite[Exercise~7.16.9(i)]{EGNO} gives
\begin{equation}\label{eq:app-bimodule-dimension-20261102}
 \frac{\FPdim(Z_i)}{\FPdim(Z_\nu^{\mathcal B})}
 =\frac{\FPdim(Q_{\mathcal L})}{\FPdim(Q)}
 =\FPdim_{Q\text{-}\mathrm{Bimod}(\mathcal L)}(Q_{\mathcal L}).
\end{equation}
A nonzero object of a fusion category has positive algebraic-integer
Frobenius--Perron dimension. Under weak integrality both class
dimensions are integers, so their quotient is rational and hence an
integer. In that case every downstairs character extends by
Proposition~\ref{prop:extent}, which gives the inclusion of class graphs.
\end{proof}

\begin{corollary}\label{cor:restriction-consequences}
Let $\C$ be a weakly integral braided fusion category over $\kk$. Let
$\mathcal B\subseteq\C$ be a full fusion subcategory with the inherited
braiding. For an algebra character $\nu$ of $K(\mathcal B)$, the quotients
$\FPdim(Z_\mu^{\C})/\FPdim(Z_\nu^{\mathcal B})$ indexed by its extensions $\mu$ to $K(\C)$
form a partition of $\FPdim(\C)/\FPdim(\mathcal B)$.
\end{corollary}
\begin{proof}
Combine the individual divisibility~\eqref{eq:braided-divisibility}
of Theorem~\ref{thm:braided-divisibility} with
\eqref{eq:im-fibre-normalized-mass-20261102}, the normalized form of
Proposition~\ref{prop:extent}.
\end{proof}

\begin{proposition}[Tannakian restriction]\label{prop:tannakian-divisibility}
Let $\C$ be a weakly integral braided fusion category over $\kk$. Let
$\Rep(G)\subseteq\C$ be a Tannakian fusion subcategory for a finite
group $G$, and let $Z_i\in Z(\C)$ be the class object corresponding to
an algebra character $\mu_i$ of $K(\C)$. If $\mu_i$ restricts to
evaluation at $g\in G$ on $K(\Rep(G))$, then $|g^G|\mid \FPdim(Z_{i})$. If $\C$ is also nondegenerate, with its canonical modular structure, and
$\mu_i$ is the Verlinde character of $X\in\Irr(\C)$, then $|g^G|^2\mid {\FPdim(X)}^2$.
\end{proposition}
\begin{proof}
In $\Rep(G)$ a character indexed by $[g]$ has formal codegree
$|C_G(g)|$ by column orthogonality~\eqref{eq:traceexpansion} for $\Rep(G)$. Its class has
dimension $|g^G|$. The first assertion now follows from
the individual divisibility~\eqref{eq:braided-divisibility} of
Theorem~\ref{thm:braided-divisibility}. For the square assertion suppose $\C$ is modular
and let $X$ be the simple indexing the character. In the
de-equivariantization by $\Rep(G)$ choose a constituent $V$ of its
forgetful image, of degree $h$. The action conjugates degrees, so the
inertia subgroup $G_V$ is contained in $C_G(h)$. For a charge
$U\in\Rep(G)$ the normalized monodromy trace on each degree-$h$
constituent evaluates its character at $h$ by the crossed braiding
\cite[Section~4.4]{DGNO}. Thus $|h^G|=|g^G|$ (the inverse degree
convention gives $[g^{-1}]$ and the same size). For the projective
inertia representation $\rho$ labelling $X$, the equivariantization
dimension formula of \cite[Corollary~2.13]{BN} and the inclusion
$G_V\subseteq C_G(h)$ give
\begin{equation}\label{eq:app-inertia-square-20261102}
 \begin{aligned}
 \FPdim(X)&=[G:G_V]\,\dim(\rho)\,\FPdim(V),\\
 \frac{\FPdim(X)}{|g^G|}&=[C_G(h):G_V]\,\dim(\rho)\,\FPdim(V).
 \end{aligned}
\end{equation}\eqlabelalias{eq:equivariantization-simple-dimension}{eq:app-inertia-square-20261102}
The right side of the second equality is an algebraic integer. Its square is rational since
$\FPdim(X)^2$ is an integer, proving $|g^G|^2\mid\FPdim(X)^2$.
\end{proof}

\subsection{Projection, balancing and dimension formulas}
For the forgetful functor $F$ and its right adjoint $I$, rigidity and
adjunction give the projection formula
\begin{equation}\label{eq:projection-right-adjoint}
 I(F(U)\otimes V)\simeq U\otimes I(V),
 \qquad U\in Z(\C),\quad V\in\C.
\end{equation}

In a pseudo-unitary braided fusion category $\C$ with its canonical ribbon
structure, write $\theta_X$ for the twist of a simple object $X$. The
balancing identity gives the monodromy scalar on each simple constituent
$T$ of $X\otimes Y$. Its trace gives the unnormalized $S$-entry.
\begin{equation}\label{eq:balancing-monodromy-trace}
 \left.c_{Y,X}c_{X,Y}\right|_T
    =\frac{\theta_T}{\theta_X\theta_Y}\,\operatorname{id}_T,
 \qquad
 S_{X,Y}=\sum_{T\in\Irr(\C)}N_{X,Y}^T
       \frac{\theta_T}{\theta_X\theta_Y}\FPdim(T).
\end{equation}

Let $\mathcal A,\mathcal B$ be fusion subcategories of a braided fusion
category $\mathcal L$. Write primes for centralizers in $\mathcal L$ and
$\vee$ for the join. By \cite[Corollary~3.11 and Theorem~3.14]{DGNO},
the double-centralizer identity and the centralizer dimension formula take the form
\begin{equation}\label{eq:muger-centralizer-identities}
 \mathcal A''=\mathcal A\vee Z_2(\mathcal L),\qquad
 \FPdim(\mathcal A)\FPdim(\mathcal A')
 =\FPdim(\mathcal L)\FPdim(\mathcal A\cap Z_2(\mathcal L)).
\end{equation}
The join dimension formula \cite[Lemma~3.38]{DGNO} gives
\begin{equation}\label{eq:join-fp-dimension}
 \FPdim(\mathcal A\vee\mathcal B)
 =\frac{\FPdim(\mathcal A)\FPdim(\mathcal B)}
        {\FPdim(\mathcal A\cap\mathcal B)}.
\end{equation}

For the trivial action of $C_n$ on $\kk^\times$, the cyclic cohomology
formula \cite[Chap.~IV, \S7, Thm.~7.1]{MacLaneHomology} gives
\begin{equation}\label{eq:cyclic-even-cohomology}
 H^{2j}(C_n,\kk^\times)\simeq\kk^\times/(\kk^\times)^n=0
 \qquad (j\ge1).
\end{equation}

\section{The regular product and the It\^o--Michler splitting}\label{sec:product}\label{app:class-regular-product-20261102}
\begin{theorem}[Regular product]\label{thm:productintro}
Let $\C$ be a weakly integral braided fusion category over $\kk$. For the
forgetful functor $F\colon Z(\C)\to\C$ and its right adjoint $I$,
write $I(\one)=\bigoplus_{i\in J}Z_i$ and $Y_i=F(Z_i)$. Let
$R_{\C_{\ad}}=\sum_{X\in\Irr(\C_{\ad})}\FPdim(X)[X]$ denote the
regular element of the adjoint subcategory. In $K(\C)$ one has
\begin{equation}\label{eq:productintro}
 \prod_{i\in J}[Y_i]
 =\frac{\prod_i \FPdim(Z_{i})}{\FPdim(\C_{\ad})}R_{\C_{\ad}},\qquad
 \frac{\prod_i \FPdim(Z_{i})}{\FPdim(\C_{\ad})}\in\mathbb Z_{>0}.
\end{equation}\eqlabelalias{explicit:80658146e8:eq0012}{eq:productintro}
This identity involves the product of \emph{all} class images. Its
two-factor analogue~\eqref{eq:cross-regular-product} in
Proposition~\ref{prop:class-graph-vanishing-criterion-20261027} holds under the cross-component vanishing hypothesis.
\end{theorem}
\begin{proof}
Fix $i$ with $\FPdim(Z_i)>1$. We first show that $\mu_j(Y_i)=0$ for
some $j$. The argument uses the Galois-product method of
\cite[Theorem~4.17]{BP}, with the class dimensions themselves as weights.
Since $\tau([Y_i])=[Y_i:\one]=1$, the trace expansion
\eqref{eq:app-class-idempotents-20261102} gives
\begin{equation}\label{eq:class-weighted-mean-20261102}
 \sum_{j\in J}\frac{\FPdim(Z_j)}{\FPdim(\C)}\mu_j(Y_i)=1,
 \qquad
 \sum_{j\in J}\frac{\FPdim(Z_j)}{\FPdim(\C)}=1 .
\end{equation}
All the values in the first mean are nonnegative by
\eqref{eq:app-positive-projection-20261102} in the proof of Proposition~\ref{prop:positivity}. The Galois group of a
normal number field containing the character values permutes the
$\mu_j$. It preserves their formal codegrees because
$f_j=\sum_X\mu_j(X)\mu_j(X^*)$ is an integer. It therefore preserves
their weights $\FPdim(Z_j)$. If all the values were positive, the product
\begin{equation}\label{eq:class-positive-norm-20261102}
 \prod_{j\in J}\mu_j(Y_i)^{\FPdim(Z_j)}
\end{equation}
would be a positive Galois-invariant algebraic integer, hence an integer
at least one. The weighted geometric mean would then be at least one.
By \eqref{eq:class-weighted-mean-20261102} the weighted arithmetic mean
is one, so equality would force every value to be one. The dimension
character takes the value $\FPdim(Y_i)=\FPdim(Z_i)>1$, a contradiction.

By symmetry of the $S$-matrix and \eqref{eq:app-master-20261102},
\begin{equation}\label{eq:class-image-reciprocity-20261102}
 \FPdim(Z_i)\mu_i(Y_j)=\FPdim(Z_j)\mu_j(Y_i).
\end{equation}\eqlabelalias{eq:reciprocity}{eq:class-image-reciprocity-20261102}
Thus for every $i$ with $\FPdim(Z_i)>1$ the row $\mu_i$ vanishes on
some class image. If $\FPdim(Z_i)=1$, then $Y_i\simeq\one$, and the
reciprocity identity~\eqref{eq:reciprocity} gives $\mu_i(Y_j)=\FPdim(Y_j)$ for every $j$.
As $|\mu_i(X)|\le\FPdim(X)$ for every simple $X$, attaining
the full dimension of an object forces equality with the dimension
character on each constituent. The simples on which $\mu_i$ equals
$\FPdim$ are closed under products, constituents, and duality.
Since the $Y_j$ generate $\C_{\ad}$, the restriction of this row to
$K(\C_{\ad})$ is the dimension character. Conversely, any character
agreeing there with $\FPdim$ satisfies
$\mu_i(XX^*)=\FPdim(X)^2$ for every simple $X$, hence
$f_i=\FPdim(\C)$ and $\FPdim(Z_i)=1$.

It follows that $\prod_j [Y_j]$ evaluates to $\prod_j\FPdim(Z_j)$
on the rows restricting to $\FPdim$ on $\C_{\ad}$ and to zero on all
other rows. The regular element $R_{\C_{\ad}}$ evaluates to
$\FPdim(\C_{\ad})$ on the first set of rows and to zero on the
second. Characters separate elements of $K(\C)$, which proves
\eqref{eq:productintro}. Its coefficient at the unit is
the positive integer
$\bigl[\bigotimes_jY_j:\one\bigr]$.
\end{proof}

\subsection{The Sylow factor and its complement}\label{sec:splitting}
\begin{proof}[Proof of Theorem~\ref{thm:main}]
By the regular-product identity~\eqref{eq:productintro} of Theorem~\ref{thm:productintro},
$\FPdim(\C_{\ad})$ divides $\prod_i\FPdim(Z_i)$, proving
\textup{(i)}$\Rightarrow$\textup{(iii)}. Under \textup{(ii)}
the canonical algebra factors over the Deligne product. The class
dimensions of a braided pointed category are one, and each class
dimension in $\mathcal B$ divides $N$ since the quotient is a
rational formal codegree and hence an integer. This gives
\textup{(ii)}$\Rightarrow$\textup{(i)}. We prove \textup{(iii)}$\Rightarrow$\textup{(ii)}.

Let $P$ be the Sylow $p$-quotient of the abelian group $U(\C)$ and
let $\mathcal B$ be the neutral component of the induced faithful
$P$-grading of $\C$. The hypothesis gives
\begin{equation}\label{eq:im-adjoint-neutral-20261102}
 |P|=p^a,\qquad
 \FPdim(\mathcal B)=\frac{|U(\C)|}{|P|}\FPdim(\C_{\ad})=N .
\end{equation}
Put $\mathcal S=Z_2(\mathcal B)$. Its centralizer
$\mathcal S'_{\C}$ contains $\mathcal B$. Every fusion subcategory
containing $\mathcal B$ is a union of full grading components. Indeed,
if it contains a simple $X$, then for any simple $Y$ in the same
degree it contains $YX^*\in\mathcal B$ and hence the summand $Y$ of
$(YX^*)X$. Its index is therefore a power of $p$. On the other hand,
the centralizer dimension formula~\eqref{eq:muger-centralizer-identities} gives
\begin{equation}\label{eq:im-adjoint-index-20261102}
 \frac{\FPdim(\C)}{\FPdim(\mathcal S'_{\C})}
 =\frac{\FPdim(\mathcal S)}
        {\FPdim(\mathcal S\cap Z_2(\C))}.
\end{equation}
The right side is an integer dividing $\FPdim(\mathcal S)$, which
divides $N$. Thus the index is one and $\mathcal S\subseteq Z_2(\C)$.
For $\mathcal Q=\mathcal B'_{\C}$ we now have
$\mathcal Q\cap\mathcal B=\mathcal S=
\mathcal B\cap Z_2(\C)$. Applying~\eqref{eq:muger-centralizer-identities} to $\mathcal B$ gives
$\FPdim(\mathcal Q)=|P|\FPdim(\mathcal S)$, so the restriction of
the $P$-grading to $\mathcal Q$ is faithful.

Write $\mathcal S=\Rep(H,z)$ and $L=H/\langle z\rangle$, allowing
$z=1$. The even part $\Rep(L)$ is transparent in $\mathcal Q$.
Its de-equivariantization $\mathcal R=\mathcal Q_L$ inherits a
faithful $P$-grading preserved degree by degree by $L$
\cite[Proposition~4.30(i)--(ii)]{DGNO}. The neutral component is
$\mathcal S_L$, which is $\Vec$ for $z=1$ and transparent $\sVec$
otherwise. Every component has the dimension of the neutral component.
In the first case each component consists of one invertible simple.
In the second case $p$ is odd, since $|H|=\FPdim(\mathcal S)$ is
prime to $p$. The transparent fermion acts freely on simples:
the balancing identity~\eqref{eq:balancing-monodromy-trace} gives $\theta_{fX}=-\theta_X$. Each component of dimension
two thus consists of a pair of invertibles. In both cases
$\mathcal R$ is pointed. The degree homomorphism from its abelian
group of invertibles induces an isomorphism from its Sylow
$p$-subgroup $A_p$ onto $P$, because its kernel has order
$|\langle z\rangle|$. Hence $L$ fixes every element of $A_p$.

We check that each $x\in A_p$ lifts to an invertible of
$\mathcal Q\simeq\mathcal R^L$. Choosing isomorphisms $s(x)\to x$
for $s\in L$ gives the scalar obstruction
$\alpha_x\in H^2(L,\kk^\times)$ to an equivariant structure.
Tensoring the chosen isomorphisms using the action tensorators gives
$\alpha_{xy}=\alpha_x\alpha_y$. Since $\alpha_{\one}=1$, the order of
$\alpha_x$ divides the order of $x$. It also divides $|L|$, since the
order of a finite group annihilates its positive-degree cohomology with
arbitrary module coefficients
\cite[Chap.~XII, \S2, Prop.~2.5, p.~236]{CE}. See also \cite[Chap.~IV, \S5, Prop.~5.3, p.~117]{MacLaneHomology}.
Coprimality therefore gives $\alpha_x=1$. The group of invertibles of $\mathcal Q$ lying over
$A_p$ consequently maps onto $A_p$ with kernel
$\Hom(L,\kk^\times)$. Its Sylow $p$-subgroup maps isomorphically
onto $A_p$. The corresponding pointed subcategory $\mathcal P$ of
$\mathcal Q$ therefore has dimension $|P|$ and has an invertible
in every degree. It centralizes $\mathcal B$ and meets $\mathcal B$
in $\Vec$. Tensoring with its invertible in a given degree identifies
$\mathcal B$ with that component. Multiplication thus gives the
braided equivalence~\eqref{explicit:80658146e8:eq0011}.

Since $\FPdim(\mathcal B)=N$ is prime to $p$, the invertibles of
$\mathcal P$ form the unique Sylow $p$-subgroup of $\operatorname{Inv}(\C)$.
By construction, $\mathcal B$ consists of the universal components
whose degrees have order prime to $p$. Finally, a prime outside
$\pi(\FPdim(\C))$ divides neither a class dimension nor
$\FPdim(\C_{\ad})$. This proves~\eqref{eq:unified-vertices-20261102}.
\end{proof}

\begin{corollary}\label{cor:im-class-vertices-20261102}
Let $\C$ be a weakly integral braided fusion category over $\kk$. Let
$\C_{\pt}$ and $\C_{\ad}$ denote its pointed and adjoint subcategories,
respectively, and $\Dcl(\C)$ its class prime graph. The prime divisors of $\FPdim(\C)/\FPdim(\C_{\pt})$ are vertices of
$\Dcl(\C)$, and in particular divide $\FPdim(\C_{\ad})$.
If $\C$ is also nondegenerate, then
\begin{equation}\label{eq:im-modular-class-vertices-20261102}
 \bigcup_{X\in\Irr(\C)}\pi\bigl(\FPdim(X)^2\bigr)
   =\pi\bigl(\FPdim(\C_{\ad})\bigr).
\end{equation}
\end{corollary}
\begin{proof}
If $p$ is absent, then its full $p$-part splits off in $\C_{\pt}$ by
Theorem~\ref{thm:main}. The assertion for modular categories follows
because their classes have dimensions $\FPdim(X)^2$
by Lemma~\ref{lem:class-square-minimal-center-20261026}\textup{(1)}.\end{proof}

\begin{corollary}[Prime support of conjugacy-class dimensions]
\label{cor:im-class-support-20261024}
Let $\C$ be a weakly integral braided fusion category over $\kk$. Let
$\operatorname{Vert}\Dcl(\C)$ be the set of primes dividing the
Frobenius--Perron dimension of at least one conjugacy class of $\C$,
and let $v_p$ denote the $p$-adic valuation. There is a braided equivalence
\begin{equation}\label{explicit-alpha:a45450d4ea:eq0082}
\begin{gathered}
 \C\simeq\mathcal P\boxtimes\E,\\
 \FPdim(\mathcal P)=\prod_{\substack{p\mid\FPdim(\C)\\p\notin \operatorname{Vert}\Dcl(\C)}}p^{v_p(\FPdim(\C))},\qquad
 \FPdim(\E)=\prod_{p\in \operatorname{Vert}\Dcl(\C)}p^{v_p(\FPdim(\C))}.
\end{gathered}
\end{equation}
where $\mathcal P$ is pointed and empty products are interpreted as $1$. The following consequences hold.
\begin{enumerate}[label=\textup{(\roman*)},leftmargin=*]
\item If $|\operatorname{Vert}\Dcl(\C)|\le2$, then $\C$ is solvable.
\item If $|\operatorname{Vert}\Dcl(\C)|\le1$, then $\C$ is nilpotent.
\item If $\operatorname{Vert}\Dcl(\C)=\varnothing$, then $\C$ is pointed.
\end{enumerate}
\end{corollary}
\begin{proof}
Apply Theorem~\ref{thm:main} successively to the primes dividing ${\FPdim(\C)}$
which do not belong to $\operatorname{Vert}\Dcl(\C)$. At each step the class dimensions of the
category before splitting are the class dimensions of the remaining
factor, each repeated the rank of the pointed factor by Proposition~\ref{prop:class-graph-product-20261026}. Hence prime
avoidance is preserved for all subsequent steps. The product of the
pointed factors and the remaining category have the dimensions in
\eqref{explicit-alpha:a45450d4ea:eq0082}.

If $|\operatorname{Vert}\Dcl(\C)|\le2$, then ${\FPdim(\E)}$ has the form $p^a q^b$, allowing zero
exponents. Theorem~1.6 of \cite{ENO} shows that $\E$ is solvable.
The group of invertible objects of the braided pointed category
$\mathcal P$ is abelian, so $\mathcal P$ is solvable by
\cite[Proposition~4.5(ii)]{ENO}. Solvability is preserved by
Deligne products \cite[Proposition~4.5(i)]{ENO}, proving
\textup{(i)}.

If $|\operatorname{Vert}\Dcl(\C)|=1$, then $\E$ has prime-power dimension and is nilpotent
\cite[Example~4.5(1)]{GN}. A pointed category has trivial adjoint
subcategory, and the adjoint construction commutes with Deligne
products. It follows that $\mathcal P\boxtimes\E$ is nilpotent.
If $\operatorname{Vert}\Dcl(\C)=\varnothing$, then ${\FPdim(\E)}=1$ and $\E=\Vec$, proving both
the remaining case of \textup{(ii)} and \textup{(iii)}.
\end{proof}

\subsection{Symmetric and modular specializations}\label{sec:im-checks-20261024}
For $\C=\Rep(G)$ with its symmetric braiding, $Y_i$ is the
permutation representation on a conjugacy class $K_i$. Thus
the value $\mu_j(Y_i)$ counts the elements of $K_i$ commuting with a
representative of the class indexed by $j$. The product in
\eqref{eq:productintro} of Theorem~\ref{thm:productintro} is the permutation representation on
$\prod_i K_i$ under diagonal conjugation. Any tuple in this product
generates $G$. Indeed, the subgroup $H$ generated by its entries meets
every conjugacy class, so every group element fixes a coset in $G/H$.
The average number of fixed cosets is one by transitivity. The identity
fixes $[G:H]$ cosets, forcing $H=G$. The stabilizer of each tuple is
therefore $Z(G)$. Every orbit affords a copy of the regular
representation of $G/Z(G)$, and its multiplicity is
\begin{equation}\label{eq:group-orbits-product-20261102}
 \frac{\prod_i |K_i|}{|G/Z(G)|}.
\end{equation}
This recovers \eqref{eq:productintro} since
$\Rep(G)_{\ad}=\Rep(G/Z(G))$. Indeed, the common kernel of the
representations $X\otimes X^*$ consists of the elements acting by a
scalar on every irreducible representation, precisely $Z(G)$.
Theorem~\ref{thm:main} specializes to the central-Sylow criterion
\eqref{explicit-alpha:a45450d4ea:eq0002}.

In a weakly integral modular category with $\FPdim(\C)=p^aN$ and
$(p,N)=1$, Theorem~\ref{thm:main} and the universal grading dimension
formula~\eqref{eq:gradingdim} give
\begin{equation}\label{explicit-alpha:a45450d4ea:eq0003}
 p\nmid\FPdim(X)^2\text{ for all } X\in\Irr(\C)
 \quad\Longleftrightarrow\quad p^a\mid |U(\C)|.
\end{equation}
The reverse implication, with $\FPdim(X)^2$, holds for every weakly
integral braided category by \cite[Theorem~1.6]{JLMS26}.
For integral modular categories the equivalence also holds with
$\FPdim(X)$, since an integer and its square have the same prime
divisors. Outside the integral case, the forward implication of the
unsquared statements in \cite[Proposition~1.3]{JIMJ25} and
\cite[Corollary~1.7]{JLMS26} fails if divisibility is taken in the
ring of algebraic integers. To see the distinction, put
$m=\FPdim(X)^2\in\mathbb Z_{>0}$. Then $\FPdim(X)/p$ is an algebraic
integer if and only if $p^2\mid m$. Indeed, its square $m/p^2$ is
rational, so integrality forces $m/p^2\in\mathbb Z$. Conversely, the
square root of an integer is an algebraic integer. Thus the two
divisibility conditions differ precisely when $v_p(m)=1$.
For example, an Ising category has simple dimensions $1,1,\sqrt{2}$.
Since $2\mid 2$ but $4\nmid 2$, the prime $2$ divides the square of the
noninvertible simple's dimension and divides none of the simple dimensions
in the ring of algebraic integers. Yet $\FPdim(\C)=4$ and $|U(\C)|=2$.
The class belonging to the noninvertible simple has dimension $2$,
as required by \eqref{explicit-alpha:a45450d4ea:eq0003}.

For modular $\C$, the equivalence
$\Z(\C)\simeq\C\boxtimes\C^{\rev}$ identifies the summands
of $\mathscr A$ with $X\boxtimes X^*$. The underlying object is
$X\otimes X^*$, and $\mu([X\otimes X^*])=|\mu([X])|^2$.
Duality permutes $\Irr(\C)$, so in a weakly integral modular category
\begin{equation}\label{eq:modular-product-BP-20261102}
\begin{aligned}
 \prod_i\frac{[Y_i]}{\FPdim(Z_i)}
 &=\prod_{X\in\Irr(\C)}
       \frac{[X\otimes X^*]}{\FPdim(X)^2}
 =\left(\prod_{X\in\Irr(\C)}\frac{[X]}{\FPdim(X)}\right)^2\\
 &=\frac{R_{\C_{\ad}}}{\FPdim(\C_{\ad})}.
\end{aligned}
\end{equation}
The last equality follows from \cite[Theorem~1.16]{BP}. Its Burnside hypothesis
holds for every weakly integral modular category. Theorem~\ref{thm:productintro}
extends this identity to degenerate braidings. The product in the regular-product identity~\eqref{eq:productintro} uses multiplication in the Grothendieck ring.
The related identity of \cite[Theorem~1.2]{Harada24} uses multiplication
of normalized class sums in the algebra of central elements.

Weak integrality matters even in the modular case. In the Fibonacci
category, write $X\otimes X\simeq\one\oplus X$ and
$\varphi=\FPdim(X)=(1+\sqrt5)/2$. The nontrivial class has underlying
object $\one\oplus X$. Its two character values are $1+\varphi$ and
$1-\varphi^{-1}$, both positive. Here $\C_{\ad}=\C$, so the
normalized regular element has value zero at the second character by
the defining property of the regular element. The normalized product of the underlying
classes has positive value there. Thus pseudo-unitarity suffices for
Proposition~\ref{prop:positivity}, while the regular-product identity~\eqref{eq:productintro}
of Theorem~\ref{thm:productintro} can fail beyond weak integrality.

\section{Prime supports, products and gradings}\label{sec:unified-prime-supports}

\subsection{The group case and Deligne products}
For a finite group $G$, the conjugacy classes of $\Rep(G)$ have dimensions
$|g^G|$.  Thus $\Delta_{\mathrm{cl}}(\Rep(G))$ is the usual prime graph of
conjugacy-class sizes, denoted by $\Delta'(G)$ in \cite{NPR}.  The character-degree graph uses the dimensions of the irreducible representations of $G$.

\begin{theorem}[Classical finite-group results]
\label{thm:group-class-graph-classical-20261027}
Let $G$ be a finite group and let $\Rep(G)$ be its category of
finite-dimensional complex representations with the usual symmetric
braiding. Its class graph $\Delta_{\mathrm{cl}}(\Rep(G))$ has at
most two connected components.  It has two components precisely when
$G/Z(G)$ is a Frobenius group and the inverse images in $G$ of its kernel
and complement are abelian.  Both components are then complete.  A
connected class graph has diameter at most three. If $G$ is nonsolvable,
it is connected and has diameter at most two.
\end{theorem}

For the graph on primes, all four assertions appear together in
\cite[Theorem~2.5]{NPR}. The component bound also follows from the
bound for the graph on noncentral classes in \cite{BHM} by prime--class
incidence. For the characterization of disconnectedness and the
diameter results on the graph on primes, see \cite{Alfandary,CasoloDolfi}. For example, the nonsolvable group $A_5$ has
class sizes $1,12,12,15,20$, and its class graph is complete on
$\{2,3,5\}$.

\begin{remark}[Absent primes and possible extensions]
A prime dividing $|G|$ is absent from the class graph exactly when its
Sylow subgroup is central \cite[Remark~2.4]{NPR}.
The categorical analogue is Theorem~\ref{thm:main}, and the corresponding
two-component bound is Theorem~\ref{thm:general-two-components-20261102}.

The completeness of the two components and the existence of a uniform
diameter bound for connected graphs remain open for general braided
fusion categories.
\end{remark}

\begin{proposition}[Deligne products]
\label{prop:class-graph-product-20261026}
Let $\C$ and $\D$ be weakly integral braided fusion categories over
$\kk$, and equip $\C\boxtimes\D$ with its product braiding. Its class
dimensions are precisely the products $\FPdim(Z_{i}^{\C})\FPdim(Z_{j}^{\D})$,
with the corresponding class objects indexed by pairs $(i,j)$.
A pointed factor does not change the class graph.  If both factors are
non-pointed, the class graph of their product is connected and has diameter
at most two.
\end{proposition}
\begin{proof}
The center, forgetful functor, and induction of the unit factor over a
Deligne product.  This gives the assertion on class dimensions.  A braided
pointed category has all its class dimensions equal to one.  Conversely,
if every class dimension of a weakly integral braided category equals one,
then $\FPdim(\C)=\operatorname{rank}(\C)$, so every simple object is
invertible.

Write $P$ and $Q$ for the nonempty sets of class primes of the two
non-pointed factors.  Any distinct $p\in P$ and $q\in Q$ are adjacent in
the product graph.  Two primes lying in the same factor can be joined
through any prime of the other factor; if that prime equals an endpoint,
the path has length one.  This proves the diameter assertion even when
$P$ and $Q$ overlap.
\end{proof}

\begin{proposition}[Braided nilpotent categories]
\label{prop:class-graph-nilpotent-20261026}
Let $\C$ be a nilpotent braided fusion category over $\kk$. Such a
category is weakly integral. Its class graph is
complete.  More precisely, a single class dimension is divisible by every
vertex of the graph.
\end{proposition}
\begin{proof}
Use the Sylow decomposition for braided nilpotent categories,
$\C\simeq\boxtimes_p\C_p$, where $\FPdim(\C_p)$ is a power of $p$;
see \cite[Theorem~1.1]{DGNO-nilpotent}.
Every integral class dimension of $\C_p$ divides $\FPdim(\C_p)$: its
quotient is a formal codegree, hence a rational algebraic integer.
Thus the nontrivial class dimensions of $\C_p$ are powers of $p$.
For each $p$ which occurs as a vertex choose a nontrivial class in
$\C_p$, and choose the unit class in all remaining factors.  The product
class supplied by Proposition~\ref{prop:class-graph-product-20261026}
has all the required prime divisors.
\end{proof}

\subsection{Class dimensions under equivariantization and grading}
In an equivariantization, the acting group permutes characters of the
category being equivariantized.  In a grading, the \emph{dual} of the
grading group acts on characters of the graded category by degree twists.

\begin{lemma}[The identity fiber of an equivariantization]
\label{lem:class-equivariantization-20261026}
Let $G$ be a finite group with a categorical action by braided tensor
autoequivalences on a weakly integral braided fusion category $\D$ over
$\kk$. Equip the equivariantization $\C=\D^G$ with its induced braiding.
Let $F\colon\C\to\D$ be the forgetful functor and $F_*$ the induced map
of complexified Grothendieck rings. For an algebra character $\nu$ of
$K(\D)$ put $\mu=\nu\circ F_*$, write $f_\nu$ and $f_\mu$ for their
respective formal codegrees, and let $G_\nu$ be the stabilizer of $\nu$
under the action
\begin{equation}\label{explicit:80658146e8:eq0013}
 (g\cdot\nu)([X])=\nu([g^{-1}(X)]).
\end{equation}
Then
\begin{equation}\label{explicit:80658146e8:eq0014}
 f_\mu=|G_\nu|f_\nu,
 \qquad
 \boxed{\FPdim(Z_{\mu}^{\C})=[G:G_\nu]\FPdim(Z_{\nu}^{\D}).}
\end{equation}
The resulting characters $\mu$ are indexed by $G$-orbits of characters
of $K(\D)$.  They exhaust the restriction fiber over the dimension
character of the canonical subcategory $\Rep(G)\subseteq\C$.
\end{lemma}
\begin{proof}
We work with complexified Grothendieck algebras.  Their characters satisfy
the orthogonality relations~\eqref{eq:traceexpansion}, and
$f_\nu=\sum_{X\in\Irr(\D)}|\nu([X])|^2$.
The simple objects of an equivariantization are indexed by an orbit of
$X\in\Irr(\D)$ and an irreducible projective representation $\pi$ of
the inertia group $G_X$.  Their forgetful images are
\begin{equation}\label{explicit:80658146e8:eq0015}
 F(S_{X,\pi})=\dim(\pi)\bigoplus_{Y\in G\cdot X}Y,
 \qquad
 \sum_\pi\dim(\pi)^2=|G_X|;
\end{equation}
see \cite[Corollary~2.13]{BN}.  Therefore
\begin{align}
 f_\mu
 &=\sum_{G\cdot X}|G_X|
       \left|\sum_{Y\in G\cdot X}\nu([Y])\right|^2\label{explicit:graphs:aligned-row-1}\\
 &=\frac1{|G|}\sum_{X\in\Irr(\D)}
       \left|\sum_{g\in G}\nu([g(X)])\right|^2
 =|G_\nu|f_\nu.
\label{explicit:graphs:aligned-row-2}\end{align}
For the last equality, expand the square and use character orthogonality~\eqref{eq:traceexpansion}. Then
the pair $(g,h)$ contributes precisely when the two corresponding
translates of $\nu$ coincide.  Since
$\FPdim(\C)=|G|\FPdim(\D)$, the formal-codegree formula~\eqref{explicit:80658146e8:eq0003}
gives the class dimension in~\eqref{explicit:80658146e8:eq0014}.

The image of $K(\C)\to K(\D)$ is $K(\D)^G$: it contains every orbit
sum, by the forgetful-image formula~\eqref{explicit:80658146e8:eq0015}.  Characters of this
invariant algebra correspond to $G$-orbits of characters of $K(\D)$.
All the resulting $\mu$ restrict to the dimension character of
$\Rep(G)$.  Their class dimensions sum to
\begin{equation}\label{explicit:80658146e8:eq0016}
 \sum_{G\cdot\nu}[G:G_\nu]\FPdim(Z_{\nu}^{\D})
 =\sum_\nu \FPdim(Z_{\nu}^{\D})=\FPdim(\D).
\end{equation}
The class-fiber identity~\eqref{eq:extent-general} of Proposition~\ref{prop:extent}
gives this same sum for the entire restriction
fiber.  Positivity of class dimensions proves exhaustion.
\end{proof}

\begin{lemma}[A faithful abelian grading]
\label{lem:class-dimensions-grading-20261027}
Let $\C$ be a weakly integral braided fusion category over $\kk$ with
a faithful grading $\C=\bigoplus_{a\in A}\C_a$ by a finite abelian
group $A$. Equip the neutral component $\mathcal B=\C_e$ with the
inherited braiding.  The character group
$\widehat A=\operatorname{Hom}(A,\mathbb C^\times)$ acts on the characters
of $K(\C)$ by
\begin{equation}\label{explicit:80658146e8:eq0017}
 (\chi\cdot\mu)([X])=\chi(a)\mu([X]),
 \qquad X\in\Irr(\C_a).
\end{equation}
For every character $\mu$,
\begin{equation}\label{explicit:80658146e8:eq0018}
 \boxed{
 \FPdim(Z_{\mu}^{\C})
 =\bigl|\operatorname{Stab}_{\widehat A}(\mu)\bigr|
       \FPdim(Z_{\mu|_{\mathcal B}}^{\mathcal B}).}
\end{equation}
The restriction fibers are exactly the $\widehat A$-orbits.
\end{lemma}
\begin{proof}
Multiplication of homogeneous classes adds their degrees, so degree
twisting is an algebra automorphism and carries algebra characters to
algebra characters.  Its invariant algebra is precisely $K(\mathcal B)$.
For a finite-dimensional commutative semisimple algebra, the characters
of the invariant algebra identify with the orbits on the character set.
This proves the statement about fibers.

Formal codegrees are constant on such an orbit, since $|\chi(a)|=1$.
Put $\nu=\mu|_{\mathcal B}$. The order of the stabilizer is
$|\operatorname{Stab}_{\widehat A}(\mu)|$.
The fiber has $|A|/|\operatorname{Stab}_{\widehat A}(\mu)|$ elements, all of the same
class dimension.  Faithfulness gives
$\FPdim(\C)=|A|\FPdim(\mathcal B)$, so the class-fiber identity~\eqref{eq:extent-general}
of Proposition~\ref{prop:extent} reads
\begin{equation}\label{explicit:80658146e8:eq0019}
 \frac{|A|}{|\operatorname{Stab}_{\widehat A}(\mu)|}\FPdim(Z_{\mu}^{\C})=|A|\FPdim(Z_{\nu}^{\mathcal B}).
\end{equation}
Cancelling $|A|$ proves the formula.
\end{proof}

\begin{example}[The grading group and its dual in the Ising example]
For an Ising category,
\begin{equation}\label{explicit:80658146e8:eq0020}
 \C_e=\langle\mathbf1,\psi\rangle,
 \qquad \C_a=\langle\sigma\rangle,
 \qquad A=\mathbb Z/2\mathbb Z.
\end{equation}
The nontrivial element of $\widehat A$ changes the sign of
$\mu(\sigma)$ and fixes $\mu(\psi)$.  It exchanges the two characters
with values $\mu(\sigma)=\pm\sqrt2$, whose class dimensions are one,
and fixes the character with $\mu(\sigma)=0$, whose class dimension is
two.  All class dimensions of the neutral pointed category are one.
This realizes both stabilizer sizes in the grading formula~\eqref{explicit:80658146e8:eq0018}
of Lemma~\ref{lem:class-dimensions-grading-20261027}.
\end{example}

\begin{corollary}[The odd-prime part of the weakly integral problem]
\label{cor:class-graph-odd-primes-20261027}
Let $\C$ be a weakly integral braided fusion category over $\kk$. Let
$\C_{\mathrm{int}}$ be the full fusion subcategory generated by the simple objects of integer
Frobenius--Perron dimension, equipped with the inherited braiding. Then
\begin{equation}\label{explicit:80658146e8:eq0021}
 \Delta_{\mathrm{cl}}(\C)\setminus\{2\}
 =\Delta_{\mathrm{cl}}(\C_{\mathrm{int}})\setminus\{2\}
\end{equation}
as induced graphs.  Every edge of the latter category's full class graph
also persists in the graph of $\C$.
If $\Delta_{\mathrm{cl}}(\C)$ has more components than
$\Delta_{\mathrm{cl}}(\C_{\mathrm{int}})$, then $2$ is a new isolated
vertex and the number of components increases by one.
\end{corollary}
\begin{proof}
The canonical dimension grading of a weakly integral category has
elementary abelian $2$-group $A$ and neutral component
$\C_{\mathrm{int}}$; see \cite{GN}.  Thus the grading formula~\eqref{explicit:80658146e8:eq0018} of
Lemma~\ref{lem:class-dimensions-grading-20261027} gives
\begin{equation}\label{explicit:80658146e8:eq0022}
 \FPdim(Z_{\mu}^{\C})=2^{e(\mu)}\FPdim(Z_{\mu|_{\C_{\mathrm{int}}}}^{\C_{\mathrm{int}}}).
\end{equation}
Every character of the neutral component occurs as a restriction.  Hence
the odd prime supports, and the simultaneous occurrence of pairs of odd
primes, agree in the two categories.  Equation~\eqref{explicit:80658146e8:eq0022} also preserves all old edges.  Only the vertex $2$ can be newly introduced, proving the final
assertion.
\end{proof}

\begin{remark}[The vertex $2$ and the integral subcategory]
When $2$ is absent from the class graph of $\C_{\mathrm{int}}$, the
class It\^o--Michler splitting~\eqref{explicit:80658146e8:eq0011} of
Theorem~\ref{thm:main} gives

\begin{equation}\label{explicit:80658146e8:eq0023}
 \C_{\mathrm{int}}\simeq\mathcal P_2\boxtimes\mathcal E,
 \qquad \mathcal P_2\text{ pointed of }2\text{-power dimension},
 \quad \FPdim(\mathcal E)\text{ odd}.
\end{equation}
If the integral part has two components supported on odd primes and $2$
occurs in $\Dcl(\C)$, Theorem~\ref{thm:general-two-components-20261102}
forces $2$ to have an odd neighbor. Thus some even class dimension of
$\C$ has an odd prime divisor.
\end{remark}

\subsection{Prime-power class dimensions and solvability}
\begin{theorem}[Prime-power conjugacy-class dimensions]
\label{thm:class-prime-powers-solvable-20261026}
Let $\C$ be a weakly integral braided fusion category over $\kk$. Suppose that the
Frobenius--Perron dimension of each conjugacy-class object of $\C$ is either one or a power of a prime, where the prime may
depend on the class.  Then $\C$ is solvable.
\end{theorem}
\begin{proof}
We argue by induction on $\FPdim(\C)$. Suppose first that every
Tannakian subcategory of $\C$ is pointed. Then
Theorem~\ref{thm:class-graph-pointed-tannakians-20261027} gives at
most two connected components.

The prime-power hypothesis makes the graph discrete, so its
vertex set has at most two elements.  Split off all absent primes by the
class It\^o--Michler splitting~\eqref{explicit:80658146e8:eq0011} of
Theorem~\ref{thm:main}.  The remaining braided factor has dimension
$p^a q^b$, allowing either exponent to be zero, and is solvable by the
categorical Burnside theorem \cite[Theorem~1.6]{ENO}.  The split-off braided pointed
factor is solvable as well.

Otherwise choose a non-pointed Tannakian subcategory
$\Rep(G)\subseteq\C$.  By Proposition~\ref{prop:extent}, each character of $K(\Rep(G))$ extends to
$K(\C)$. Tannakian restriction (Proposition~\ref{prop:tannakian-divisibility})
then gives the divisibility~\eqref{eq:braided-divisibility} for every conjugacy-class size of $G$.  Every nontrivial class size of $G$ is therefore a prime power.  Its class graph is discrete and has at most two vertices by
Theorem~\ref{thm:group-class-graph-classical-20261027}.  Primes absent
from that graph have central Sylow subgroups.  Splitting off those
central factors reduces $G$ to a group of order $p^a q^b$, so $G$ is
solvable by Burnside's theorem.

Let $\mathcal B=(\Rep(G))'_{\C}$ and let
$\mathcal T=\C_G$ be the de-equivariantization.  Its neutral component
$\mathcal T_e$ is braided and
\begin{equation}\label{explicit:80658146e8:eq0024}
 \mathcal B\simeq(\mathcal T_e)^G.
\end{equation}
The category $\mathcal T_e$ is weakly integral, and the induced
$G$-action on it is braided.
By Theorem~\ref{thm:braided-divisibility}, every class dimension of
$\mathcal B$ divides a class dimension of $\C$, so $\mathcal B$ retains the
prime-power property.  By
the identity-fiber formula~\eqref{explicit:80658146e8:eq0014} of
Lemma~\ref{lem:class-equivariantization-20261026}, every class dimension
of $\mathcal T_e$ divides a class dimension of $\mathcal B$.  Hence
$\mathcal T_e$ has the same property and is solvable by induction, since
\begin{equation}\label{explicit:80658146e8:eq0025}
 \FPdim(\mathcal T_e)=\frac{\FPdim(\mathcal B)}{|G|}
 <\FPdim(\C).
\end{equation}
The faithful support of the crossed grading of $\mathcal T$ is a subgroup
of the solvable group $G$.  Therefore $\mathcal T$ is a solvable-group
extension of the solvable category $\mathcal T_e$, and
$\C\simeq\mathcal T^G$ is solvable as well.  Here we use closure of
solvable fusion categories under extensions and equivariantizations by
solvable finite groups \cite{ENO}.
\end{proof}

The general component bound of Theorem~\ref{thm:general-two-components-20261102}
gives another proof of Theorem~\ref{thm:class-prime-powers-solvable-20261026}.
The prime-power hypothesis makes the graph discrete, so at most two
primes occur. The assertion then follows from Corollary~\ref{cor:im-class-support-20261024}\textup{(i)} on prime supports.

\begin{remark}[Comparison with simple-object dimensions]
If the Frobenius--Perron dimension of every simple object is one or
a prime power, \cite[Theorem~1.3]{NPR} gives weak group-theoreticity.
It gives solvability if at most two primes occur or every prime that
occurs exceeds seven.

Here the prime-power assumption concerns conjugacy-class dimensions.
For integral modular categories the two prime-power hypotheses agree,
since the class dimensions are exactly the squares of the simple
dimensions. Outside the modular case the dimensions may have different
prime supports.
\end{remark}

\section{Coprime class dimensions and Tannakian subcategories}\label{sec:unified-tannakian}

\subsection{Cross-component vanishing and the regular-object identity}
Write $\mathscr A=I_{\C}(\mathbf1)$, $Z_i\prec\mathscr A$,
$Y_i=F(Z_i)$.  All characters below are algebra
characters of the complexified Grothendieck ring.  The proposed vanishing
statement is
\begin{equation}\label{eq:cross-component-vanishing-target}
 \mu_i(Y_j)=0\quad\text{if }\FPdim(Z_{i}),\FPdim(Z_{j})>1\text{ belong to different
 components of }\Delta_{\mathrm{cl}}(\C).
\end{equation}
The prime divisors of each nontrivial class dimension lie in a unique
component. Whether \eqref{eq:cross-component-vanishing-target} holds
for every weakly integral braided fusion category is
Problem~\ref{prob:cross-component-vanishing-20261027}.

\begin{proposition}[A sufficient vanishing criterion]
\label{prop:class-graph-vanishing-criterion-20261027}
Let $\C$ be a weakly integral braided fusion category over $\kk$. Let
$Z_i\in Z(\C)$ be its class objects, $Y_i$ their forgetful images, and
$\mu_i$ the corresponding algebra characters of $K(\C)$. Suppose that
$\mu_i(Y_j)=0$ whenever $\FPdim(Z_i),\FPdim(Z_j)>1$ have prime
supports in distinct connected components of $\Dcl(\C)$, as
in~\eqref{eq:cross-component-vanishing-target}. Then
$\Delta_{\mathrm{cl}}(\C)$ has at most two connected components.  More precisely, for every such pair,
\begin{equation}\label{eq:cross-regular-product}
 [Y_i][Y_j]=\frac{\FPdim(Z_{i})\FPdim(Z_{j})}{\FPdim(\C_{\ad})}R_{\C_{\mathrm{ad}}},\qquad
 R_{\C_{\mathrm{ad}}}
 =\sum_{X\in\Irr(\C_{\mathrm{ad}})}\FPdim(X)[X].
\end{equation}
In particular $\FPdim(\C_{\ad})\mid \FPdim(Z_{i})\FPdim(Z_{j})$.
\end{proposition}
\begin{proof}
As shown in the proof of Theorem~\ref{thm:productintro}, the class-size-one
characters are exactly the characters whose restriction to $K(\C_{\mathrm{ad}})$
is $\FPdim$.  Consequently
\begin{equation}\label{explicit:80658146e8:eq0026}
 \mu_s(R_{\C_{\mathrm{ad}}})=
 \begin{cases}\FPdim(\C_{\ad}),&\FPdim(Z_{s})=1,\\0,&\FPdim(Z_{s})>1.\end{cases}
\end{equation}
If $\FPdim(Z_{s})=1$, then $Y_i,Y_j\in\C_{\mathrm{ad}}$ give
$\mu_s(Y_iY_j)=\FPdim(Z_{i})\FPdim(Z_{j})$.  If $\FPdim(Z_{s})>1$, its component is different from at
least one of the components of $\FPdim(Z_{i}),\FPdim(Z_{j})$, so the assumed vanishing gives
$\mu_s(Y_iY_j)=0$.  Thus both sides of
\eqref{eq:cross-regular-product} have the same value under every character.
The complexified Grothendieck algebra is commutative semisimple, so its
characters separate its elements.  This proves the equality.

The coefficient of $\mathbf1$ on the right is $\FPdim(Z_{i})\FPdim(Z_{j})/\FPdim(\C_{\ad})$, whereas
the corresponding coefficient on the left is an integer.  Hence
$\FPdim(\C_{\ad})\mid \FPdim(Z_{i})\FPdim(Z_{j})$.  If three components existed, choose class objects $Z_a,Z_b,Z_c$ whose nontrivial dimensions lie
in the three components. These dimensions are pairwise coprime, and the
preceding divisibility shows that $\FPdim(\C_{\ad})$ divides each of
$\FPdim(Z_a)\FPdim(Z_b)$, $\FPdim(Z_a)\FPdim(Z_c)$, and
$\FPdim(Z_b)\FPdim(Z_c)$.  Hence $\FPdim(\C_{\ad})=1$.  This makes
$\C$ pointed, contradicting the existence of these nontrivial class
dimensions.
\end{proof}

\begin{lemma}[The coprime dichotomy]
\label{lem:coprime-class-vanishing-dichotomy-20261027}
Let $\C$ be a weakly integral braided fusion category over $\kk$. Let
$Z_i,Z_j$ be class objects in $Z(\C)$, let $\mu_i$ be the algebra
character of $K(\C)$ corresponding to $Z_i$, and let $Y_j$ be the
forgetful image of $Z_j$. If $\gcd(\FPdim(Z_i),\FPdim(Z_j))=1$, then
\begin{equation}\label{explicit:80658146e8:eq0027}
 \mu_i(Y_j)\in\{0,\FPdim(Z_{j})\}.
\end{equation}
The second alternative is equivalent to $Z_i$ and $Z_j$ centralizing each
other in $Z(\C)$.
\end{lemma}
\begin{proof}
Use the unnormalized $S$-matrix of $Z(\C)$.  The class-row identity~\eqref{explicit:80658146e8:eq0007}
and positivity~\eqref{eq:positive} from Proposition~\ref{prop:positivity} give
\begin{equation}\label{explicit:80658146e8:eq0028}
 S_{Z_i,Z_j}=\FPdim(Z_{i})\mu_i(Y_j)\geq0.
\end{equation}
Both $S_{Z_i,Z_j}/\FPdim(Z_{i})$ and $S_{Z_i,Z_j}/\FPdim(Z_{j})$ are algebraic integers.
Bezout's identity therefore shows that
$\alpha=S_{Z_i,Z_j}/(\FPdim(Z_{i})\FPdim(Z_{j}))$ is an algebraic integer.
The objects $Z_i,Z_j$ are integral, so every constituent of their tensor
product is integral.  Their twists are one, so the balancing formula~\eqref{eq:balancing-monodromy-trace} expresses
$\alpha$ as
\begin{equation}\label{explicit:80658146e8:eq0029}
 \alpha=\frac{1}{\FPdim(Z_{i})\FPdim(Z_{j})}
 \sum_W N_{Z_i,Z_j}^{W}\FPdim(W)\theta_W.
\end{equation}
The coefficients are nonnegative integers with total $\FPdim(Z_{i})\FPdim(Z_{j})$.
Since twists are roots of unity, every Galois conjugate of $\alpha$ has
absolute value at most one.  Kronecker's theorem implies that $\alpha$ is
zero or a root of unity.  Its nonnegative real value leaves only zero or
one.  In the latter case equality in the pseudo-unitary $S$-matrix bound
is precisely the centralization criterion.
\end{proof}

\subsection{The common Tannakian category and additional summands}
For a centralizing pair of class objects of coprime dimensions,
the common constituents of their adjoint products determine a
Tannakian subcategory. The construction applies for every value of
$\dim\operatorname{Hom}(Z_iZ_j,\mathscr A)$.

\begin{proposition}[The common Tannakian category and its image]
\label{prop:common-tannakian-localization-20261027}
Let $\C$ be a weakly integral braided fusion category over $\kk$. Let
$F\colon Z(\C)\to\C$ be the forgetful functor, and let
$Z_i,Z_j\in Z(\C)$ be conjugacy-class objects that centralize each
other in $Z(\C)$ and have coprime Frobenius--Perron dimensions. Put
$Y_i=F(Z_i)$ and $Y_j=F(Z_j)$. In $Z(\C)$ set
\begin{equation}\label{explicit:80658146e8:eq0030}
 \mathcal D_i=\langle Z_iZ_i^*\rangle,\qquad
 \mathcal D_j=\langle Z_jZ_j^*\rangle,
\end{equation}
and let $\mathcal E$ be generated by all simple objects occurring in both
$Z_iZ_i^*$ and $Z_jZ_j^*$.  Then $\mathcal E$ is Tannakian.  The full fusion subcategory $\mathcal K=F(\mathcal E)$
generated by its forgetful image is Tannakian in $\C$ for the inherited
braiding. The restriction $F|_{\mathcal E}$ is braided, and
$\mathcal E$ centralizes the canonical lift $\mathcal K^+$ through the
given braiding of $\C$.

If $Z_iZ_j$ is reducible, then $\mathcal E$ is non-pointed.  If, in addition,
every Tannakian subcategory of $\C$ is pointed, then
\begin{equation}\label{explicit:80658146e8:eq0031}
 \langle Y_i,Y_j\rangle\subsetneq\C.
\end{equation}
In this case $\mathcal E$ fails to centralize $X^+$ for some simple
$X\notin\langle Y_i,Y_j\rangle$.
\end{proposition}
\begin{proof}
The categories $\mathcal D_i,\mathcal D_j$ centralize one another.  A common
simple constituent $U$ belongs to $\mathcal D_j$, so it centralizes $Z_i$;
also $\operatorname{Hom}(U Z_i,Z_i)\ne0$.  The balancing identity~\eqref{eq:balancing-monodromy-trace} then gives
$\theta_U=1$.  The category generated by these objects is symmetric,
being contained in $\mathcal D_i\cap\mathcal D_j$, and its generating
objects have twist one.  Thus it is Tannakian.

Set $\mathcal L_i=\langle Y_i\rangle$ and
$\mathcal L_j=\langle Y_j\rangle$.  The internal-end identity~\eqref{eq:offdiagonal} of
Proposition~\ref{prop:diagonal}, written as
\begin{equation}\label{explicit:80658146e8:eq0032}
 [Y_i:X]=\dim\operatorname{Hom}(X^+Z_i,Z_i),
\end{equation}
implies $\mathcal L_i^+\subseteq\mathcal D_i$, and similarly for $j$.
Therefore $\mathcal L_i$ and $\mathcal L_j$ centralize.  Moreover
$\mathcal K\subseteq\mathcal L_i\cap\mathcal L_j$, so $\mathcal K$ is
symmetric.  Since
$\mathcal E\subseteq\mathcal D_i$ and
$\mathcal K^+\subseteq\mathcal L_j^+\subseteq\mathcal D_j$, the category
$\mathcal E$ centralizes $\mathcal K^+$.
For $U,V\in\mathcal E$, this gives
\begin{equation}\label{explicit:80658146e8:eq0033}
 F(c^{Z(\C)}_{U,V})
 =(c^{\C}_{F(V),F(U)})^{-1}
 =c^{\C}_{F(U),F(V)},
\end{equation}
the last equality using symmetry of $\mathcal K$.  Consequently the
dominant functor $F|_{\mathcal E}:\mathcal E\to\mathcal K$ is braided,
and $\mathcal K$ is Tannakian.

Rigidity gives
\begin{equation}\label{explicit:80658146e8:eq0034}
 \dim\operatorname{End}(Z_iZ_j)
 =\dim\operatorname{Hom}(Z_iZ_i^*,Z_jZ_j^*).
\end{equation}
If the product is reducible, these two internal-end objects have a common
nontrivial constituent.  A common invertible constituent stabilizes both
$Z_i$ and $Z_j$. Each stabilizer has order dividing the corresponding
squared dimension by \cite[Remark~9.4]{NPR}. Since these squares are
coprime, the stabilizers intersect trivially. This is the stabilizer
argument used in the simplicity criterion of \cite[Lemma~9.3]{NPR}.
Hence a noninvertible common constituent exists, and $\mathcal E$ is
non-pointed.

Under the additional hypothesis, $\mathcal K$ is pointed.  If
$\langle Y_i,Y_j\rangle=\C$, then $\mathcal E$ centralizes $\C^+$ and
therefore lies in $(\C^+)'=\C^-$.  The functor $F$ is fully faithful on
$\C^-$; it would identify the non-pointed $\mathcal E$ with the pointed
$\mathcal K$, a contradiction.  The same argument shows that some simple
object outside $\langle Y_i,Y_j\rangle$ detects $\mathcal E$ by monodromy.
\end{proof}

\begin{proposition}[A sufficient condition for cross-component vanishing]
\label{prop:class-vanishing-pointed-tannakian-center-20261027}
Let $\C$ be a weakly integral braided fusion category over $\kk$. Suppose every Tannakian
fusion subcategory of the Drinfeld center $Z(\C)$ is pointed. Let
$Z_i,Z_j$ be conjugacy-class objects in $Z(\C)$ with Frobenius--Perron
dimensions greater than one
whose prime supports lie in distinct components of $\Dcl(\C)$.
Then the algebra character $\mu_i$ corresponding to $Z_i$ vanishes
on the forgetful image $Y_j$ of $Z_j$, as
in~\eqref{eq:cross-component-vanishing-target}.  In particular
$\Delta_{\mathrm{cl}}(\C)$ has at most two components.
\end{proposition}
\begin{proof}
Take class dimensions in different components.  They are coprime.  If the
evaluation were nonzero, the coprime dichotomy~\eqref{explicit:80658146e8:eq0027} of
Lemma~\ref{lem:coprime-class-vanishing-dichotomy-20261027} would imply centralization.
Proposition~\ref{prop:common-tannakian-localization-20261027} then shows
that the product must be simple, since a reducible product would generate
a non-pointed Tannakian category in the center.  But a simple product
$Z_iZ_j$ is itself a class object: its forgetful image contains
$\mathbf1$, and adjunction places it in $\mathscr A$.  Its dimension
$\FPdim(Z_{i})\FPdim(Z_{j})$ connects the two components, a contradiction.  Now apply
Proposition~\ref{prop:class-graph-vanishing-criterion-20261027}.
\end{proof}

\begin{lemma}[Common constituents of the class images]
\label{lem:class-images-common-constituents-20261027}
Let $\C$ be a weakly integral braided fusion category over $\kk$. Let
$F\colon Z(\C)\to\C$ be the forgetful functor with right adjoint $I$.
Let $Z_i,Z_j$ be simple summands of $\mathscr A=I(\one)$, and
put $Y_i=F(Z_i)$ and $Y_j=F(Z_j)$. If $Z_i,Z_j$ centralize each
other in $Z(\C)$ and $\gcd(\FPdim(Z_i),\FPdim(Z_j))=1$, then
\begin{equation}\label{explicit:80658146e8:eq0035}
 [Y_iY_j:\mathbf1]>1
 \quad\Longrightarrow\quad
 \C\text{ contains a non-pointed Tannakian subcategory}.
\end{equation}
Consequently, if every Tannakian subcategory of $\C$ is pointed, then
$\dim\Hom_{Z(\C)}(Z_iZ_j,\mathscr A)=1$.
\end{lemma}
\begin{proof}
The objects $Y_i,Y_j$ are self-dual and contain the unit once by
Proposition~\ref{prop:positivity}. Their generated subcategories centralize
by the proof of Proposition~\ref{prop:common-tannakian-localization-20261027}.  A common
invertible constituent $g$ would stabilize both $Z_i,Z_j$, so coprimeness
forces $g=\mathbf1$.  Thus $[Y_iY_j:\mathbf1]>1$ gives a common
noninvertible simple $X$.
By the internal-end identity~\eqref{eq:offdiagonal} of
Proposition~\ref{prop:diagonal}, every common simple $X$ satisfies
$X^+\prec Z_iZ_i^*,Z_jZ_j^*$.  It centralizes $Z_i$, and
$\operatorname{Hom}(X^+Z_i,Z_i)\ne0$, so $\theta_X=1$ by the balancing identity~\eqref{eq:balancing-monodromy-trace}.
The common constituents therefore generate a Tannakian subcategory of
the symmetric intersection $\langle Y_i\rangle\cap\langle Y_j\rangle$;
it is non-pointed.  Finally, adjunction identifies
$\operatorname{Hom}(Z_iZ_j,\mathscr A)$ with
$\operatorname{Hom}(Y_iY_j,\mathbf1)$.
\end{proof}

\begin{lemma}[Constraints on the additional summands]
\label{lem:additional-class-summands-20261027}
Let $\C$ be a weakly integral braided fusion category over $\kk$. For the
forgetful functor $F\colon Z(\C)\to\C$ and its right adjoint $I$, let
$\mathscr A=I(\one)=\bigoplus_a Z_a$ and $Y_a=F(Z_a)$. Suppose
two of these class objects $Z_i,Z_j$ centralize each other in $Z(\C)$, have coprime
Frobenius--Perron dimensions, and satisfy
$\dim\Hom_{Z(\C)}(Z_iZ_j,\mathscr A)=1$. Then $Z_iZ_j$ has a unique
class constituent $Z_k$, of multiplicity one, and
\begin{equation}\label{explicit:80658146e8:eq0036}
 Z_iZ_j=Z_k\oplus R,\qquad
 R=\bigoplus_{W\not\prec\mathscr A}m_WW.
\end{equation}
Here $m_W$ is the multiplicity of the simple $W\in Z(\C)$ in
$Z_iZ_j$. For the canonical ribbon twist on the center, every
additional constituent satisfies
\begin{equation}\label{explicit:80658146e8:eq0037}
 \theta_W=1,\qquad [F(W):\mathbf1]=0.
\end{equation}
Furthermore,
\begin{align}
 Y_iY_j&=Y_k\oplus\bigoplus_W m_W F(W),\label{explicit:graphs:aligned-row-3}\\
 \sum_Wm_W\FPdim(W)&=\FPdim(Z_{i})\FPdim(Z_{j})-\FPdim(Z_{k}),\label{explicit:graphs:aligned-row-4}\\
 1+\sum_Wm_W^2
 &=\dim\operatorname{Hom}(Z_iZ_i^*,Z_jZ_j^*).
\label{explicit:graphs:aligned-row-5}\end{align}
\end{lemma}
\begin{proof}
Multiplicity-freeness of $\mathscr A$ and the stated Hom dimension give
the unique class constituent.  The balancing identity~\eqref{eq:balancing-monodromy-trace} gives trivial twist on the whole
product, hence on each constituent.  Adjunction gives
$[F(W):\mathbf1]=0$ for nonclass constituents.  Applying $F$, taking
dimensions, and computing the dimension of the endomorphism algebra of
the product give the remaining formulas.
\end{proof}

\begin{lemma}[A negative character contribution from the extra summands]
\label{lem:extra-summands-negative-character-20261027}
Let $\C$ be a weakly integral braided fusion category over $\kk$. For the
forgetful functor $F\colon Z(\C)\to\C$ with right adjoint $I$, write
$I(\one)=\bigoplus_a Z_a$, $Y_a=F(Z_a)$, and let $\mu_a$ be the
algebra character of $K(\C)$ corresponding to $Z_a$. Suppose that
two of these class objects $Z_i,Z_j$ centralize each other in $Z(\C)$, have coprime
Frobenius--Perron dimensions, and that
$\dim\Hom_{Z(\C)}(Z_iZ_j,I(\one))=1$. Decompose
$Z_iZ_j=Z_k\oplus R$ as in
Lemma~\ref{lem:additional-class-summands-20261027}, with $Z_k$ its
unique class constituent. If $R\ne0$, then some $\mu_s$ satisfies
\begin{equation}\label{explicit:80658146e8:eq0038}
 \mu_s(Y_k)>\mu_s(Y_i)\mu_s(Y_j)\geq0.
\end{equation}
If this witness satisfies $\gcd(\FPdim(Z_{s}),\FPdim(Z_{i})\FPdim(Z_{j})\FPdim(Z_{k}))=1$, then
\begin{equation}\label{explicit:80658146e8:eq0039}
 \mu_s(F(R))=-\FPdim(Z_{k}),\qquad 2\FPdim(Z_{k})\leq \FPdim(Z_{i})\FPdim(Z_{j}).
\end{equation}
\end{lemma}
\begin{proof}
The actual object $V=F(R)$ has zero unit multiplicity and positive
dimension.  It is self-dual, since
$[V]=[Y_i][Y_j]-[Y_k]$ and all class images are self-dual.
The trace expansion~\eqref{eq:traceexpansion} and the formal-codegree
formula~\eqref{explicit:80658146e8:eq0003} give
\begin{equation}\label{explicit:80658146e8:eq0040}
 0=[V:\mathbf1]=\frac1{\FPdim(\C)}\sum_s \FPdim(Z_{s})\mu_s(V).
\end{equation}
All evaluations are real, all weights are positive, and the
$\FPdim$ character contributes positively.  Thus some evaluation is
negative, proving the first assertion by the positivity of class images
in Proposition~\ref{prop:positivity}.
Under the additional coprimeness, the dichotomy~\eqref{explicit:80658146e8:eq0027}
of Lemma~\ref{lem:coprime-class-vanishing-dichotomy-20261027} forces
$\mu_s(Y_k)=\FPdim(Z_{k})$ and
$\mu_s(Y_i)\mu_s(Y_j)=0$, because the only other possible product is
$\FPdim(Z_{i})\FPdim(Z_{j})>\FPdim(Z_{k})$.  Finally
$|\mu_s(V)|\leq\FPdim(V)=\FPdim(Z_{i})\FPdim(Z_{j})-\FPdim(Z_{k})$ gives the inequality.
\end{proof}

\subsection{The additional summands in the group case}
The pure-flux formula~\eqref{eq:centralizer-fiber-class-product} is the single-double-coset case of
\cite[Theorem~4.1 and Example~4.3]{BN} for the untwisted double.
Coprimality makes the diagonal conjugation action transitive and gives
the following direct description of all multiplicities.
\begin{lemma}[The centralizer-fiber formula]
\label{lem:coprime-group-class-product-fiber-20261027}
Let $G$ be a finite group and let $\Rep(G)$ be its category of
finite-dimensional complex representations with its symmetric braiding.
Use the canonical ribbon structure on $Z(\Rep(G))\simeq\Rep(D(G))$.
Let $a,b\in G$ have coprime conjugacy-class sizes. Write $Z_{a,\rho}$
for the simple object indexed by the class of $a$ and an irreducible
representation $\rho$ of $C_G(a)$.  Put
$K=C_G(a)\cap C_G(b)$ and $L=C_G(ab)$.  Then
\begin{equation}\label{eq:centralizer-fiber-class-product}
 Z_{a,\mathbf1}\otimes Z_{b,\mathbf1}
 \simeq\bigoplus_{\rho\in\Irr(L)}
       \dim(\rho^K) Z_{ab,\rho}.
\end{equation}
There is exactly one canonical class summand, $Z_{ab,\mathbf1}$, of
multiplicity one.  The product is simple if and only if $L=K$.
If $a,b$ commute, the two class objects centralize and every summand in
\eqref{eq:centralizer-fiber-class-product} has twist one.
\end{lemma}
\begin{proof}
The index $[G:K]$ is divisible by both coprime indices
$[G:C_G(a)]$, $[G:C_G(b)]$ and is at most their product.  It therefore
equals that product.  The diagonal $G$-action on $a^G\times b^G$ is
transitive.  Its multiplication map has image $(ab)^G$ and fiber $L/K$
over $ab$.  The pure-class tensor product is the corresponding equivariant
vector bundle, with fiber representation $\operatorname{Ind}_K^L\mathbf1$.
Frobenius reciprocity gives the decomposition~\eqref{eq:centralizer-fiber-class-product}.
The trivial representation occurs once, and the permutation representation
is irreducible exactly when its dimension is one, namely when $L=K$.

If $a,b$ commute, transitivity shows that every member of $a^G$ commutes
with every member of $b^G$.  This is centralization for these pure-class
objects.  Moreover $ab\in K\cap Z(L)$.  It acts trivially on $L/K$, hence
on every occurring $\rho$, which is exactly the trivial-twist condition.
\end{proof}

\begin{example}[Extra summands with a graph connection elsewhere]
\label{ex:coprime-centralizing-class-extra-summands-20261027}
For an odd prime $p$, let
\begin{equation}\label{explicit:80658146e8:eq0041}
 V=\mathbb F_p^2,\qquad
 H=\left\{\begin{pmatrix}t&s\\0&1\end{pmatrix}:
 t\in\mathbb F_p^\times,\ s\in\mathbb F_p\right\},\qquad
 G=V\rtimes H.
\end{equation}
This is the linear action of \cite[Example~3.1]{DGPS}, written
using column vectors. We compute the tensor product of the two
associated class objects.

For $a=(1,0)$ and $b=(0,1)$,
\begin{equation}\label{explicit:80658146e8:eq0042}
 a^G=\{(u,0):u\ne0\},\qquad
 b^G=\{(v,1):v\in\mathbb F_p\}.
\end{equation}
Their sizes are $p-1,p$, and the two classes commute because they lie in
the abelian normal subgroup $V$.  Their product has support $b^G$.
The fiber over $b$ consists of
\begin{equation}\label{explicit:80658146e8:eq0043}
 ((u,0),(-u,1)),\qquad u\in\mathbb F_p^\times.
\end{equation}
The quotient $C_G(b)/V\simeq\mathbb F_p^\times$ acts regularly on it.
The centralizer-fiber formula~\eqref{eq:centralizer-fiber-class-product}
of Lemma~\ref{lem:coprime-group-class-product-fiber-20261027} therefore gives
\begin{equation}\label{eq:affine-extra-summands}
 Z_{a,\mathbf1}\otimes Z_{b,\mathbf1}
 \simeq\bigoplus_{\chi\in\widehat{\mathbb F_p^\times}}Z_{b,\chi},
\end{equation}
where each $\chi$ is inflated to $C_G(b)$.  Every summand has dimension
$p$ and twist one.  Exactly the trivial character yields a canonical
class object; the other $p-2$ summands have no unit in their forgetful
images.

For $p=3$, this is the decomposition $2\cdot3=3+3$ in the center of
the representation category of a group of order $54$.
The element
\begin{equation}\label{explicit:80658146e8:eq0044}
 u=\begin{pmatrix}1&1\\0&1\end{pmatrix}\in H\subseteq G
\end{equation}
has $|C_H(u)|=p$ and $|C_V(u)|=p$, so its class size is $p(p-1)$.
This supplies a graph connection outside the product under consideration.
More generally, the multiplicities for all odd primes $p$ are
\begin{equation}\label{eq:affine-multiplicities-20261102}
\begin{array}{c|ccccc}
 \text{class size} & 1 & p-1 & p & p(p-1) & p^2 \\
 \hline
 \text{number of classes} & 1 & 1 & p-1 & p & p(p-2).
\end{array}
\end{equation}
Indeed, an element $(u,v,r,h)$, $r\ne0$, commutes with $(x,y,t,s)$
exactly when $(t-1)u=(r-1)x+hy-sv$ and $(t-1)h=(r-1)s$.
For $t\ne1$ the free parameters are $r,v$, giving centralizer order
$p(p-1)$. For $t=1$, $s\ne0$ one has $r=1$ and $sv=hy$, giving
order $p^2$. For $t=1$, $s=0$, the centralizer has order $p^2(p-1)$
if $y\ne0$ and order $p^3$ if $y=0$, $x\ne0$. The remaining
element is the identity. Dividing the number of elements in each
case by the corresponding class size gives the table.
\end{example}

\begin{remark}[Reducible products of coprime class objects]
In Example~\ref{ex:coprime-centralizing-class-extra-summands-20261027} the two class objects centralize each other and
have coprime dimensions. Their tensor product has a unique canonical
class constituent and all its summands have twist one. Nevertheless
the product is reducible, and every summand has dimension $p$. The
ambient category $\Rep(G)$ is itself a non-pointed Tannakian category.
\end{remark}

\begin{example}[The canonical constituents need not exhibit a local connection]
\label{ex:class-product-no-local-bridge-20261027}
Let $\C=\Rep(D(A_5))$, and take the electric irreducibles $X,Y$ of
dimensions $3,4$.  They centralize, and the ordinary $A_5$ fusion rule is
\begin{equation}\label{explicit:80658146e8:eq0045}
 X\otimes Y\simeq X'\oplus Y\oplus T,
 \qquad \FPdim(X')=3,\quad\FPdim(T)=5.
\end{equation}
Under $Z(\C)\simeq\C\boxtimes\C^{\mathrm{rev}}$, the corresponding
class objects are $X\boxtimes X^*$ and $Y\boxtimes Y^*$.  Their dimensions
are $9,16$, and the canonical constituents of their tensor product have
dimensions $9,16,25$.  None of these dimensions joins the prime $3$ of
the first class to the prime $2$ of the second.  Nevertheless $\C$ has
a simple pure-flux object of dimension $12$, coming from a five-cycle
class in $A_5$, and hence has a class dimension $12^2=144$.
The primes $2$ and $3$ are therefore joined by this class.
\end{example}

\subsection{Pointed forgetful images and the two-component theorem}
\label{sec:pointed-tannakian-bound}

When every Tannakian subcategory of $\C$ is pointed, a failure of
cross-component vanishing produces a solvable Tannakian subcategory
of $Z(\C)$ with pointed forgetful image. Its representation degrees
occur among the class dimensions of $\C$ and bound the number of
components.

We retain $F:Z(\C)\to\C$, its right adjoint $I$, the canonical algebra
$\mathscr A=I(\mathbf1)=\bigoplus_i Z_i$, and $Y_i=F(Z_i)$.
For $X\in\C$, write $X^+=\iota_+(X)$ for its lift determined by the
given braiding.

\begin{lemma}[Pointed images are translates of class objects]
\label{lem:pointed-image-class-translation-20261027}
Let $\C$ be a braided fusion category over $\kk$, let
$F:Z(\C)\to\C$ be the forgetful functor, and let $I$ be its right
adjoint. Suppose $W\in\Irr(Z(\C))$ and $F(W)$ is a direct sum of
invertible objects. Write $g^+$ for the lift given by the braiding.
For every invertible constituent $g\prec F(W)$, there exists, up to
isomorphism, a unique simple summand $Z_\ell$ of $I(\one)$ such that
\begin{equation}\label{explicit:80658146e8:eq0046}
 W\simeq g^+\otimes Z_\ell.
\end{equation}
Consequently $\FPdim(W)=\FPdim(Z_{\ell})$, and $Y_\ell=F(Z_\ell)$ is pointed.
In particular, let $\mathcal P\subseteq Z(\C)$ be a full fusion
subcategory tensor equivalent to $\Rep(G)$ for a finite group $G$.
If $F(\mathcal P)$ is pointed, every ordinary irreducible degree of $G$
is a class dimension of $\C$.
\end{lemma}

\begin{proof}
Adjunction and the projection formula~\eqref{eq:projection-right-adjoint} give
\begin{equation}\label{explicit:80658146e8:eq0047}
 \operatorname{Hom}_{\C}(F(W),g)
 \simeq\operatorname{Hom}_{Z(\C)}(W,I(g)),
 \qquad
 I(g)\simeq g^+\otimes\mathscr A.
\end{equation}
Thus $(g^+)^*\otimes W$ is a simple summand of $I(\one)$.
Its isomorphism class is uniquely determined by $g$ and $W$.
Applying $F$ gives $F(W)\simeq g\otimes Y_\ell$.
The last assertion follows by applying this argument to each simple
object of $\mathcal P$.
\end{proof}

\begin{remark}
For the stabilizer description, assume further that $\C$ is
pseudo-unitary and equip it with its canonical positive spherical
structure. Put
\begin{equation}\label{explicit:80658146e8:eq0048}
 G_{\C}[Z_\ell]
 =\{a\in\Irr(\C_{\mathrm{pt}}):a^+\otimes Z_\ell\simeq Z_\ell\}.
\end{equation}
For each invertible $a$, the internal-end identity~\eqref{eq:offdiagonal} of Proposition~\ref{prop:diagonal} gives
$[Y_\ell:a]=\dim\Hom_{Z(\C)}(a^+\otimes Z_\ell,Z_\ell)$.
By Schur's lemma this is one for $a\in G_{\C}[Z_\ell]$ and zero
otherwise. Hence, if $Y_\ell$ is pointed, then
\begin{equation}\label{explicit:80658146e8:eq0049}
 Y_\ell=\bigoplus_{a\in G_{\C}[Z_\ell]}a,
 \qquad \FPdim(Z_{\ell})=|G_{\C}[Z_\ell]|.
\end{equation}
If $F(W)$ lies in a pointed fusion subcategory $\mathcal K$ and
$g\prec F(W)$, then $Y_\ell\in\mathcal K$ as well.  In particular,
$\FPdim(W)$ divides $\FPdim(\mathcal K)$.
These conclusions concern pointed forgetful images; no analogous
divisibility $\FPdim(Z_{\ell})\mid\FPdim(W)$ is asserted for an arbitrary
central object in the same $\C^+$-module block.
\end{remark}

\begin{proposition}[A solvable Tannakian subcategory with pointed image]
\label{prop:solvable-tannakian-class-graph-20261027}
Let $\C$ be a weakly integral braided fusion category over $\kk$,
and let $F:Z(\C)\to\C$ be the forgetful functor. Suppose that
$Z(\C)$ contains a Tannakian subcategory
\begin{equation}\label{explicit:80658146e8:eq0050}
 \mathcal P\simeq\Rep(G)
\end{equation}
for a finite nonabelian solvable group $G$, and that $F(U)\in\C_{\pt}$
for every $U\in\mathcal P$.  Then $\Delta_{\mathrm{cl}}(\C)$ has at most two connected
components.
\end{proposition}

\begin{proof}
We first establish the class-object version of the Gallagher argument.
Let $\mathcal Q\simeq\Rep(L)$ be a fusion subcategory of $\mathcal P$,
and let $Z_i$ be a class object such that $\gcd(\FPdim(Z_{i}),|L|)=1$.
De-equivariantizing $Z(\C)$ by $\mathcal Q$ gives a fusion category
$\mathcal M$ with an $L$-action and $Z(\C)\simeq\mathcal M^L$.
Taking dimensions in the forgetful-image formula~\eqref{explicit:80658146e8:eq0015} gives
\begin{equation}\label{explicit:80658146e8:eq0051}
 \FPdim(Z_{i})=[L:L_V]\dim(\pi)\FPdim(V),
\end{equation}
where $V\in\Irr(\mathcal M)$, $L_V$ is its inertia subgroup, and
$\pi$ is an irreducible projective representation of $L_V$.
Because $\FPdim(Z_{i})$ is an integer, $\FPdim(V)$ is a rational algebraic
integer and is therefore an integer.  Every prime dividing
$[L:L_V]\dim(\pi)$ divides $|L|$.  Coprimeness now forces
$[L:L_V]=\dim(\pi)=1$.

It follows that the de-equivariantization of $Z_i$ is simple.
Consequently $U\otimes Z_i$ is simple for every
$U\in\Irr(\mathcal Q)$: equivalently, the corresponding projective
representation $\pi$ is one-dimensional, and tensoring it with an
ordinary irreducible representation preserves irreducibility.
This is the argument of \cite[Corollary~6.2 and Lemma~7.2]{NPR}
applied to the integral dimension $\FPdim(V)$.

By Lemma~\ref{lem:pointed-image-class-translation-20261027}, write
$U\simeq g^+\otimes Z_\ell$.  Then
\begin{equation}\label{explicit:80658146e8:eq0052}
 Z_\ell\otimes Z_i
 \simeq(g^+)^{-1}\otimes U\otimes Z_i
\end{equation}
is simple.  Its forgetful image contains $\mathbf1$, since both
$Y_\ell$ and $Y_i$ contain $\mathbf1$.  Adjunction therefore shows
that it is a constituent of $\mathscr A$, and hence is a class object.
We have proved the implication
\begin{equation}
\begin{aligned}
 &\gcd(\FPdim(Z_i),|L|)=1\\
 &\quad\Longrightarrow\quad
 \dim(U)\FPdim(Z_i)\text{ is a class dimension of }\C.
\end{aligned}
 \label{eq:class-gallagher-product}
\end{equation}
This holds for every $U\in\Irr(L)$.

Choose $N\lhd G$ maximal among normal subgroups for which $L=G/N$
is nonabelian.  Then all proper quotients of $L$ are abelian, and
$\Rep(L)\subseteq\mathcal P$ still has pointed forgetful image.
The solvable-group classification used in
\cite[proof of Theorem~7.7]{NPR} gives two possibilities:
either $L$ is a nonabelian $p$-group, or $L$ is a Frobenius group
with elementary abelian kernel of order $r^n$ and abelian complement $K$. Its order is $r^n|K|$, and $|K|$ is its unique nonlinear irreducible degree.

In the first case choose a nonlinear degree $p^a$.
It is a class dimension by
Lemma~\ref{lem:pointed-image-class-translation-20261027}.
Every class dimension divisible by $p$ connects all its prime
divisors to $p$.  Every class dimension coprime to $p$ does so through
the class dimension $p^a \FPdim(Z_{i})$ from
\eqref{eq:class-gallagher-product}.  Thus the graph is connected.

In the Frobenius case, $|K|$ is a class dimension, so $\pi(|K|)$ is a
clique.  If $r\mid \FPdim(Z_{i})$, all prime divisors of $\FPdim(Z_{i})$ lie in the
component containing $r$.  If $r\nmid \FPdim(Z_{i})$ but
$\gcd(\FPdim(Z_{i}),|K|)>1$, they connect to $\pi(|K|)$.  In the remaining case
$\gcd(\FPdim(Z_{i}),|L|)=1$, and $|K| \FPdim(Z_{i})$ supplies the same connection by
\eqref{eq:class-gallagher-product}.  There are at most two components;
the component containing $r$ need not exist or be distinct from the
one containing $\pi(|K|)$.
\end{proof}

\begin{proposition}[Solvable enlargement of a pointed-image obstruction]
\label{prop:pointed-image-solvable-envelope-20261027}
Let $\C$ be a braided fusion category over $\kk$, let
$F:Z(\C)\to\C$ be the forgetful functor, and let
$\mathcal E\subseteq Z(\C)$ be a Tannakian fusion subcategory.
Suppose its image fusion subcategory $\mathcal K=F(\mathcal E)$ is
pointed Tannakian for the braiding of $\C$. Write
$\mathcal K^+=\iota_+(\mathcal K)$ for its lift by that braiding, and
suppose that $\mathcal E$ centralizes $\mathcal K^+$ in $Z(\C)$.
Then
\begin{equation}\label{explicit:80658146e8:eq0053}
 \mathcal P=\mathcal E\vee\mathcal K^+
 \simeq\Rep(\Gamma)
\end{equation}
is Tannakian for a solvable finite group $\Gamma$.
More precisely, writing $\mathcal K\simeq\Rep(H)$ with $H$ abelian,
one has $\Gamma\simeq N\rtimes H$, where $N$ is nilpotent of class
at most two.  If $\mathcal E$ is non-pointed, then $\Gamma$ is
nonabelian.
\end{proposition}

\begin{proof}
The join of the two centralizing Tannakian subcategories
$\mathcal E$ and $\mathcal K^+$ is Tannakian, and
$F(\mathcal P)=\mathcal K$.
The restriction $F|_{\mathcal P}$ is braided.  Indeed, for
$U\in\mathcal E$ and $V\in\mathcal K$, centralization with $V^+$
identifies the half-braiding of $U$ on $V$ with
$c_{V,F(U)}^{-1}=c_{F(U),V}$, where the equality uses symmetry of
$\mathcal K$.  This proves compatibility on $\mathcal E$ and on
crossings with $\mathcal K^+$; compatibility on $\mathcal K^+$
itself holds by definition.  It therefore holds on their join.

The dominant symmetric functor
\begin{equation}\label{explicit:80658146e8:eq0054}
 F|_{\mathcal P}:\Rep(\Gamma)\longrightarrow\Rep(H)
\end{equation}
has the symmetric section $\mathcal K^+\hookrightarrow\mathcal P$.
Tannaka duality gives homomorphisms
\begin{equation}\label{explicit:80658146e8:eq0055}
 H\xrightarrow{\iota}\Gamma\xrightarrow{q}H,
 \qquad q\iota=\operatorname{id}_H,
\end{equation}
where $\iota$ is injective and $q$ is surjective.  Hence
$\Gamma=N\rtimes H$, with $N=\ker q$.

Let $B=\operatorname{Fun}(H)$ be the regular algebra in
$\mathcal K^+$ and set $\D=\C_H$, using the compatible module
convention for de-equivariantization. The equivalence of centers of module categories
\cite{Schauenburg}, in the form of \cite[Theorem~3.20]{DMNO}, gives
\begin{equation}
 Z(\D)\simeq Z(\C)^0_B
 \simeq\bigl((\mathcal K^+)'\bigr)_H.
 \label{eq:envelope-center-deequivariantization}
\end{equation}
Because $\mathcal P\subseteq(\mathcal K^+)'$, this contains
\begin{equation}\label{explicit:80658146e8:eq0056}
 \mathcal P_H\simeq\Rep(N).
\end{equation}

Under the forgetful functor $Z(\D)\to\D$, this entire subcategory
has trivial underlying objects.  Explicitly, its generating free
modules have the form $U\otimes B$, $U\in\mathcal P$, and their
underlying $\D$-objects are
\begin{equation}\label{explicit:80658146e8:eq0057}
 F(U)\otimes B\simeq
 \mathbf1_{\D}^{\oplus\FPdim(U)},
\end{equation}
since $F(U)\in\Rep(H)$.  Every object of $\mathcal P_H$ is a
direct summand of such a free module.  This also verifies the
needed compatibility of \eqref{eq:envelope-center-deequivariantization}
with forgetting, rather than merely using an abstract equivalence
of centers.

For any fusion category $\D$, the fusion subcategory of $Z(\D)$
consisting of objects with trivial underlying object is
$\Rep(U(\D))$, where $U(\D)$ is the universal grading group.
One can see this directly: the half-braiding of
$V\otimes\mathbf1$ with a simple object $X$ specifies an invertible
linear map on $V$; the hexagon identifies these maps according to
the defining relations of $U(\D)$.  Thus the inclusion
$\Rep(N)\subseteq\Rep(U(\D))$ makes $N$ a quotient of $U(\D)$.

The category $\D$ is braided $H$-crossed by
\cite[Section~4.4]{DGNO}.  Its $H$-grading factors through a
homomorphism
\begin{equation}\label{explicit:80658146e8:eq0058}
 \partial:U(\D)\longrightarrow H.
\end{equation}
No surjectivity onto $H$ is needed.  If $u\in\ker\partial$, choose
a simple $X$ of universal degree $u$.  Its $H$-degree is neutral,
so the crossed braiding provides $X\otimes Y\simeq Y\otimes X$
for every simple $Y$.  Comparing universal degrees proves that
$u$ commutes with every element of $U(\D)$.  Therefore
\begin{equation}\label{explicit:80658146e8:eq0059}
 \ker\partial\subseteq Z(U(\D)),\qquad
 [U(\D),U(\D)]\subseteq\ker\partial,
\end{equation}
where the second inclusion uses that $H$ is abelian.
It follows that $U(\D)$, and hence $N$, is nilpotent of class at
most two.  The extension $N\rtimes H$ is solvable.  Finally,
$\mathcal E$ non-pointed implies that $\mathcal P$ is non-pointed,
so $\Gamma$ is nonabelian.
\end{proof}

\begin{theorem}[The bound with pointed Tannakian subcategories]
\label{thm:class-graph-pointed-tannakians-20261027}
Let $\C$ be a weakly integral braided fusion category over $\kk$.  If every
Tannakian subcategory of $\C$ is pointed, then
\begin{equation}\label{explicit:80658146e8:eq0060}
 \Delta_{\mathrm{cl}}(\C)
 \text{ has at most two connected components}.
\end{equation}
\end{theorem}

\begin{proof}
Suppose, to the contrary, that the graph has at least three
components.  If all evaluations between classes in different
components vanish, then
Proposition~\ref{prop:class-graph-vanishing-criterion-20261027}
already gives a contradiction.  Thus there exist class objects
$Z_i,Z_j$ whose nontrivial dimensions belong to different
components and for which $\mu_i(Y_j)\ne0$.
Their dimensions are coprime.  By
Lemma~\ref{lem:coprime-class-vanishing-dichotomy-20261027},
$Z_i$ and $Z_j$ centralize each other.

Their product is not simple: a simple product would have a unit
in its forgetful image and hence would be a class object of
dimension $\FPdim(Z_{i})\FPdim(Z_{j})$, connecting the two components.
Consequently
\begin{equation}\label{explicit:80658146e8:eq0061}
 \dim\operatorname{Hom}(Z_iZ_i^*,Z_jZ_j^*)
 =\dim\operatorname{End}(Z_iZ_j)>1.
\end{equation}
The common-constituent construction in
Proposition~\ref{prop:common-tannakian-localization-20261027}
therefore produces a non-pointed Tannakian category
$\mathcal E\subseteq Z(\C)$, with
$\mathcal K=F(\mathcal E)$ Tannakian in $\C$, such that
$\mathcal E$ centralizes $\mathcal K^+$.
By hypothesis $\mathcal K$ is pointed.

Proposition~\ref{prop:pointed-image-solvable-envelope-20261027}
then yields a nonabelian solvable Tannakian subcategory
$\mathcal P=\mathcal E\vee\mathcal K^+$ with pointed forgetful
image.  Applying
Proposition~\ref{prop:solvable-tannakian-class-graph-20261027}
gives at most two components, a contradiction.
\end{proof}

\begin{remark}[The role of the pointed-image hypothesis]
The hypothesis concerns Tannakian subcategories of $\C$. The proof uses
a possibly non-pointed Tannakian subcategory $\mathcal E\subseteq Z(\C)$
whose image $F(\mathcal E)$ is pointed. When its image is non-pointed,
the general component bound instead follows by equivariantization in
Theorem~\ref{thm:general-two-components-20261102}. Determining whether
cross-component vanishing holds in this case remains a separate problem.
\end{remark}

\section{The modular and slightly degenerate cases}\label{sec:unified-modular}

\subsection{Modular categories and the group-theoretic core}
The arguments in this section concern the prime graph of the integer
squares of simple dimensions.  For a modular category this is precisely the
class graph, even when the category is not integral.
For integral braided categories, the coprime-dimension monodromy
criterion is \cite[Lemma~7.1]{ENO}. The following version uses
coprime integer squares in a weakly integral modular category.

\begin{lemma}[Coprime squares and monodromy]
\label{lem:coprime-squares-monodromy-20261027}
Let $\mathcal M$ be a weakly integral modular fusion category over
$\kk$, and let $S$ be its unnormalized $S$-matrix for the
canonical positive spherical structure. Let $X,Y\in\Irr(\mathcal M)$.  If
\begin{equation}\label{explicit:80658146e8:eq0062}
 \gcd\bigl(\FPdim(X)^2,\FPdim(Y)^2\bigr)=1,
\end{equation}
then either $S_{X,Y}=0$ or $X$ and $Y$ projectively centralize each other.
\end{lemma}
\begin{proof}
Write $\alpha=S_{X,Y}/({\FPdim(X)}{\FPdim(Y)})$, using the unnormalized $S$-matrix.  The two
numbers $S_{X,Y}/{\FPdim(X)}$ and $S_{X,Y}/{\FPdim(Y)}$ are algebraic integers, since they
are values of Verlinde characters.  Thus ${\FPdim(X)}^2\alpha^2$ and
${\FPdim(Y)}^2\alpha^2$ are algebraic integers.  Bezout's identity implies that
$\alpha^2$, and hence $\alpha$, is an algebraic integer.

The balancing formula~\eqref{eq:balancing-monodromy-trace} expresses $\alpha$ as a weighted sum of roots of
unity, with positive weights summing to one.  The square of the dimension
of every simple object is an integer.  Consequently any Galois
automorphism changes its dimension only by a sign.  Applying a Galois
automorphism to this expression therefore shows that every
conjugate of $\alpha$ has absolute value at most one.  Kronecker's theorem
implies that $\alpha$ is zero or a root of unity.  In the second case,
equality in the triangle inequality for this weighted sum
forces all monodromy eigenvalues on $X\otimes Y$ to be equal.  This is
projective centralization.
\end{proof}

\begin{lemma}[Squares realized as class dimensions]
\label{lem:class-square-minimal-center-20261026}
Let $\C$ be a weakly integral braided fusion category over $\kk$,
equipped with its canonical positive spherical structure.
\begin{enumerate}
\item If $\C$ is modular, its class dimensions are exactly
$\FPdim(X)^2$, for $X\in\Irr(\C)$, with the corresponding class image
$X\otimes X^*$.
\item If $\C$ is slightly degenerate, then every number $\FPdim(X)^2$,
$X\in\Irr(\C)$, occurs among the class dimensions of $\C$.
The characters with $\mu([f])=1$, where $f$ is the transparent fermion,
give one such class per fermion orbit. The other characters can yield
different dimensions, as in Example~\ref{ex:fermion-ising-cube-nonsquare-class}.
\end{enumerate}
\end{lemma}
\begin{proof}
In the modular case,
$Z(\C)\simeq\C\boxtimes\C^{\mathrm{rev}}$ and the canonical algebra is
$\bigoplus_X X\boxtimes X^*$, which proves the assertion.

For the second assertion, choose a minimal modular extension
$\mathcal M\supseteq\C$, whose existence follows from \cite[Corollary~1.3]{JFR}.
Let $f$ be the transparent fermion of $\C$.  Then
$\C=\langle f\rangle'_{\mathcal M}$ and
$\C'_{\mathcal M}=\langle f\rangle$.  The modular quotient formula~\eqref{eq:modular-quotient} of Theorem~\ref{thm:modular-quotient} for
the character obtained from $X\in\Irr(\C)$ gives
\begin{equation}\label{explicit:80658146e8:eq0063}
 \FPdim(Z_{\mu_X}^{\C})=
 \frac{{\FPdim(X)}^2}{\FPdim((X\otimes X^*)\cap\langle f\rangle)}.
\end{equation}
The fermion acts freely on the simple objects of $\C$: if $fX\simeq X$,
the balancing identity~\eqref{eq:balancing-monodromy-trace} and the trivial monodromy of $f$ with $X$ would give
$\theta_X=-\theta_X$.  Hence the denominator is one. Proposition~\ref{prop:fermion-sign-class-dimensions}
below identifies exactly the characters with $\mu([f])=1$ by the reduced
character table.
\end{proof}

\begin{proposition}[The group-theoretic core in a modular category]
\label{prop:class-graph-tannakian-core-20261027}
Let $\mathcal M$ be a weakly integral modular fusion category over
$\kk$. Suppose that $\mathcal M$ contains a non-pointed Tannakian
fusion subcategory braided equivalent to $\Rep(G)$ for a finite group $G$.  Then
$\Delta_{\mathrm{cl}}(\mathcal M)$ has at most two connected components.
If $G$ is nonsolvable, this graph is connected.
\end{proposition}
\begin{proof}
Remove the central Sylow subgroups of $G$.  Their product is a central
Hall subgroup, and $G=P\times H$, where $P$ is abelian and $H$ has no
central Sylow subgroup.  The subgroup $H$ is nonabelian.  The ordinary
conjugacy-class graph of $H$ has vertex set $\pi(H)$ by the central Sylow
criterion~\eqref{explicit-alpha:a45450d4ea:eq0002}. By
Theorem~\ref{thm:group-class-graph-classical-20261027}, it has at most two
components and is connected if $H$ is nonsolvable. See also \cite[Section~2]{NPR}.

We regard $\Rep(H)$ as a subcategory of $\mathcal M$ and
de-equivariantize.  The resulting braided $H$-crossed category has
faithful $H$-grading.  If a simple object $V$ has degree $g$, then its
inertia subgroup $H_V$ is contained in $C_H(g)$.  By the equivariantization dimension formula~\eqref{eq:equivariantization-simple-dimension},
a simple equivariant object $W$ lying over $V$ has dimension
\begin{equation}\label{explicit:80658146e8:eq0064}
 {\FPdim(W)}=[H:H_V]\dim(\rho){\FPdim(V)}.
\end{equation}
The de-equivariantization is weakly integral, so ${\FPdim(V)}^2$ is an integer.
It follows that $|g^H|^2$ divides ${\FPdim(W)}^2$.  Each graded component is
nonzero, so the ordinary class graph of $H$ embeds in the class graph of
$\mathcal M$.

Let $q$ be any vertex outside $\pi(H)$, and choose a simple $X$ such that
$q\mid {\FPdim(X)}^2$.  If ${\FPdim(X)}^2$ has a prime divisor in $\pi(H)$, then $q$ is
already connected to the embedded graph.  Otherwise
$\gcd({\FPdim(X)}^2,|H|)=1$.  The equivariantization dimension formula~\eqref{eq:equivariantization-simple-dimension}, after
squaring, shows that the inertia index and the projective representation
degree are both one.  Thus the image of $X$ under de-equivariantization
is simple.  Choose any nonlinear $U\in\Irr(H)$.  The categorical Gallagher theorem \cite[Corollary~6.2]{NPR} implies that $U\otimes X$ is simple.  Its square
dimension ${\FPdim(U)}^2{\FPdim(X)}^2$ connects $q$ to a prime divisor of ${\FPdim(U)}$ in
$\pi(H)$.  Every additional vertex therefore attaches to the embedded
group graph, proving both claims.
\end{proof}

In the integral modular case, \cite[Theorem~1.1 and
Proposition~10.4]{NPR} give the bound of three in general and the bound
of two when the category is solvable or every Tannakian subcategory
is pointed. The following
theorem gives the bound of two for the graph of squared dimensions
throughout the weakly integral modular case.

\begin{theorem}[The modular two-component bound]
\label{thm:class-graph-modular-two-20261027}
Let $\mathcal M$ be a weakly integral modular fusion category over
$\kk$. Then
$\Delta_{\mathrm{cl}}(\mathcal M)$ has at most two connected components.
\end{theorem}
\begin{proof}
If $\mathcal M$ contains a non-pointed Tannakian subcategory, apply
Proposition~\ref{prop:class-graph-tannakian-core-20261027}.
Suppose that every Tannakian subcategory of $\mathcal M$ is
pointed. Let $X,Y$ label classes in different components. Their
squared dimensions are coprime. If $S_{X,Y}\ne0$, then
Lemma~\ref{lem:coprime-squares-monodromy-20261027} gives projective
centralization. The stabilizers of $X^*$ and $Y$ have coprime orders,
since those orders divide $\FPdim(X)^2$ and $\FPdim(Y)^2$, respectively.
It follows that $X\otimes Y$ is simple
\cite[Lemma~9.3 and Remark~9.4]{NPR}. Its squared dimension would
join the two components. Thus $S_{X,Y}=0$. For the Verlinde character
$\mu_X$ and the class image $Y\otimes Y^*$ this gives
\begin{equation}\label{eq:modular-cross-component-class-value-20261102}
 \mu_X(Y\otimes Y^*)=
 \frac{|S_{X,Y}|^2}{\FPdim(X)^2}=0.
\end{equation}
Proposition~\ref{prop:class-graph-vanishing-criterion-20261027} now
bounds the number of components. For integral categories the
simple-product or zero-entry alternative used here is
\cite[Corollary~9.5]{NPR}.
\end{proof}

\begin{example}\label{ex:class-graph-double-s3-20261102}
Take $G=S_3=C_3\rtimes C_2$ in the Frobenius-group family of
\cite[Example~8.6]{NPR}. In $\mathcal M=\Rep(D(S_3))$, the centralizers
of the three flux representatives are $S_3$, $C_2$, and $C_3$. The corresponding simple
dimensions, with multiplicity, are $(1,1,2)$, $(3,3)$, and $(2,2,2)$.
By Lemma~\ref{lem:class-square-minimal-center-20261026}\textup{(1)}, the class dimensions have distinct values $1,4,9$, so the graph has
two isolated vertices.  On the other
hand, $Z(\mathcal M)\simeq\mathcal M\boxtimes\mathcal M^{\mathrm{rev}}$
has class dimensions including $4\cdot9=36$, so its class graph is
connected. More generally, the graph of the center of an integral
modular category is complete, as observed in \cite[Example~10.1]{NPR}.
\end{example}

\subsection{Slightly degenerate categories and an Ising modification}
Let $\C$ be a weakly integral slightly degenerate braided fusion category. With its canonical
positive spherical structure it admits a pseudo-unitary minimal modular
extension $\mathcal M$ by \cite[Corollary~1.3]{JFR}. We have
\begin{equation}\label{explicit:80658146e8:eq0065}
 \FPdim(\mathcal M)=2\FPdim(\C),\qquad
 \C=\langle f\rangle'_{\mathcal M},
\end{equation}
where $f$ is the fermion.  Write
$\mathcal M=\mathcal M_0\oplus\mathcal M_1$ for the grading detected by
monodromy with $f$, so that $\mathcal M_0=\C$.  
By Proposition~\ref{prop:extent}, every character of $\C$ is the restriction of a Verlinde character of
$\mathcal M$, and the modular quotient formula~\eqref{eq:modular-quotient} of Theorem~\ref{thm:modular-quotient} reads
\begin{equation}
 \FPdim(Z_{\mu_W}^{\C})=\frac{\FPdim(W)^2}{1+\delta_{fW,W}}.
 \label{eq:fermion-class-quotient-20261027}
\end{equation}
For the restricted character $\mu_W$, the value $\mu_W([f])$ is the
monodromy sign of $W$ with $f$. Thus the odd sector supplies the
negative fiber of Proposition~\ref{prop:fermion-sign-class-dimensions};
its class dimensions need not be simple squares in $\C$.

\begin{lemma}[Parity in the odd sector]
\label{lem:fermion-parity-20261027}
Let $\C$ be a weakly integral slightly degenerate braided fusion
category over $\kk$, and let $f$ be the fermion generating $Z_2(\C)$. Let $\mathcal M$ be a minimal modular extension
of $\C$. It is weakly integral, with $\C=\langle f\rangle'_{\mathcal M}$ and
$\FPdim(\mathcal M)=2\FPdim(\C)$. Write
$\mathcal M=\mathcal M_0\oplus\mathcal M_1$ for the faithful
$C_2$-grading defined by monodromy with $f$, with $\mathcal M_0=\C$.
Suppose $W\in\Irr(\mathcal M_1)$ satisfies $f\otimes W\simeq W$ and
$v_2(\FPdim(W)^2)=1$.  Then $\FPdim(V)^2$ is even for every
$V\in\Irr(\mathcal M_1)$.
\end{lemma}
\begin{proof}
Assume that $\FPdim(V)^2$ is odd.  The product $W\otimes V$ belongs to $\C$ and
is fixed by tensoring with $f$.  Since $f$ acts freely on the simples of
$\C$, its constituents occur in pairs $X,fX$ of equal multiplicity and
equal dimension.

All constituents of $W\otimes V$ belong to one component of the
canonical dimension grading of the weakly integral category
$\mathcal M$.  Thus their dimensions are integer multiples of the same
square root $\sqrt{s}$, where $s$ is squarefree.  Write
${\FPdim(W)}{\FPdim(V)}=b\sqrt{s}$.  The paired decomposition implies that $b$ is even.
Thus $4\mid b^2s=\FPdim(W)^2\FPdim(V)^2$, contradicting
$v_2(\FPdim(W)^2\FPdim(V)^2)=1$.
\end{proof}

\begin{proposition}[Modular realization of the class graph]
\label{prop:class-graph-modular-realization-20261027}
Let $\C$ be a weakly integral slightly degenerate braided fusion
category over $\kk$. There is
a weakly integral modular fusion category $\mathcal N$ containing
$\C$ as a full braided fusion subcategory such that
\begin{equation}\label{explicit:80658146e8:eq0066}
 \Delta_{\mathrm{cl}}(\mathcal N)
 =\Delta_{\mathrm{cl}}(\C).
\end{equation}
Moreover $\mathcal N$ can be chosen to have dimension
$2\FPdim(\C)$.
\end{proposition}
\begin{proof}
Start with the minimal extension $\mathcal M$ in
\eqref{explicit:80658146e8:eq0065}. For each $f$-fixed simple $W$, the fermion specialization~\eqref{eq:fermion-class-quotient-20261027}
of the modular quotient formula~\eqref{eq:modular-quotient} of Theorem~\ref{thm:modular-quotient} gives the integer
class dimension $\FPdim(W)^2/2$. Thus $\FPdim(W)^2$ is even. If it is divisible by four for
every fixed $W$, dividing these squares by two preserves all their prime
divisors. The free orbits contribute the same squares in both categories,
so $\mathcal N=\mathcal M$ has the required graph.

Suppose now that $v_2(\FPdim(W)^2)=1$ for some fixed $W$.
Choose a positive Ising modular category $\mathcal I$ with fermion
$\psi$ and noninvertible simple $\sigma$. Write
$\mathcal A=\mathcal M\boxtimes\mathcal I$ and $b=f\boxtimes\psi$.
The product of the two fermionic self-braidings is trivial, so
$\langle b\rangle\simeq\Rep(C_2)$ is Tannakian. Take its connected
\'etale regular algebra and the category of local modules
\begin{equation}\label{explicit:80658146e8:eq0067}
 B=\one\oplus b,\qquad \mathcal N=\mathcal A^0_B.
\end{equation}
This is the Ising modification of \cite[Theorem~3.13]{BRW} in the
unitary case. Here \cite[Corollaries~3.30 and~3.32]{DMNO} apply to
the nondegenerate category $\mathcal A$ and give nondegeneracy of
$\mathcal N$ and
\begin{equation}\label{eq:ising-local-extension-dimension-20261102}
 \FPdim(\mathcal N)
 =\frac{\FPdim(\mathcal M)\FPdim(\mathcal I)}{\FPdim(B)^2}
 =\frac{4\FPdim(\mathcal M)}{2^2}=2\FPdim(\C).
\end{equation}
Hence $\mathcal N$ is weakly integral and is modular for its canonical
positive spherical structure.

Let $F_B\colon\mathcal A\to\mathcal A_B$, $F_B(T)=T\otimes B$, be the
free right $B$-module functor \cite[(15)]{DMNO}. By the locality criterion
\cite[Lemma~3.15]{DMNO}, $F_B(T)$ is local precisely when $T$ centralizes $b$,
since $B=\one\oplus b$. Thus $F_B$ restricts to a braided tensor functor
$\langle b\rangle'_{\mathcal A}\to\mathcal N=\mathcal A_B^0$.
Every simple local module is a summand of $F_B(T)$ for some simple $T$
in this centralizer \cite[Example~3.14]{DMNO}. For $X,Y\in\C$, the objects
$X\boxtimes\one$ and $Y\boxtimes\one$ centralize $b$. Adjunction gives
\begin{equation}\label{eq:ising-free-embedding-hom-20261102}
 \Hom_{\mathcal N}(F_B(X\boxtimes\one),F_B(Y\boxtimes\one))
 \simeq \Hom_{\mathcal A}(X\boxtimes\one,(Y\boxtimes\one)\otimes B)
 \simeq \Hom_{\C}(X,Y),
\end{equation}
since $\Hom_{\mathcal I}(\one,\psi)=0$. By~\eqref{eq:ising-free-embedding-hom-20261102}, the functor
$X\mapsto F_B(X\boxtimes\one)$ embeds $\C$ as a full braided fusion
subcategory of $\mathcal N$.

Let $U_B\colon\mathcal A_B\to\mathcal A$ be the forgetful functor.
The underlying object $U_B(F_B(T))=T\otimes B$ has dimension
$2\FPdim_{\mathcal A}(T)$, since $\FPdim_{\mathcal A}(B)=2$.
In $\mathcal A_B$ one has
$\FPdim_{\mathcal A_B}(F_B(T))=\FPdim_{\mathcal A}(T)$ because
the tensor functor $F_B$ preserves Frobenius--Perron dimensions
\cite[Proposition~3.3.13(1)]{EGNO}. For $T$ centralizing $b$, this is
also $\FPdim_{\mathcal N}(F_B(T))$, since $\mathcal N$ is a fusion
subcategory of $\mathcal A_B$.

The simples centralizing $b$ are $X\boxtimes\one$, $X\boxtimes\psi$
with $X\in\Irr(\C)$, and $V\boxtimes\sigma$ with
$V\in\Irr(\mathcal M_1)$, by the monodromy signs in the two factors.
Since $F_B(T\otimes b)\simeq F_B(T)$, the first two families give the
embedded simples of $\C$. For $T=V\boxtimes\sigma$ the same adjunction
gives $\End_{\mathcal N}(F_B(T))\simeq
\Hom_{\mathcal A}(T,T\oplus(T\otimes b))$ as vector spaces.
If $fV\not\simeq V$, this space has dimension one. Thus $F_B(T)$ is
simple in $\mathcal N$ and
$\FPdim_{\mathcal N}(F_B(T))=\FPdim_{\mathcal A}(T)=\sqrt{2}\,\FPdim_{\mathcal M}(V)$.

If $fV\simeq V$, choose an isomorphism $T\otimes b\to T$. By the cyclic cohomology formula~\eqref{eq:cyclic-even-cohomology},
$H^2(C_2,\kk^\times)=0$, so the isomorphism can be rescaled so that
applying it twice agrees with multiplication $b\otimes b\to\one$.
The two choices of sign then give inequivalent simple
$B$-modules on $T$, each occurring once in $F_B(T)$ by adjunction.
They exhaust $F_B(T)$, whose endomorphism space has dimension two.
Twisting both module actions by the algebra automorphism of $B$
acting by $-1$ on $b$ preserves the relative tensor product over $B$
and the locality condition. It defines a tensor autoequivalence of
$\mathcal N$ interchanging these two summands, which therefore have
equal Frobenius--Perron dimension in $\mathcal N$. Each summand has
dimension $\FPdim_{\mathcal N}(F_B(T))/2=\FPdim_{\mathcal A}(T)/2
=\FPdim_{\mathcal M}(V)/\sqrt{2}$.

Combining these calculations with \eqref{eq:fermion-class-quotient-20261027}
gives all the class dimensions in $\C$ and all the simple squares in
$\mathcal N$, up to repetition.
\begin{equation}\label{explicit:80658146e8:eq0068}
\begin{array}{c|c|c}
\text{objects in }\mathcal M&
 \text{class dimension in }\C&
 \text{simple square in }\mathcal N\\ \hline
 X\in\Irr(\C)&\FPdim(X)^2&\FPdim(X)^2\\
 W\in\Irr(\mathcal M_1),\ fW\simeq W&\FPdim(W)^2/2&\FPdim(W)^2/2\\
 \{V,fV\}\subset\Irr(\mathcal M_1),\ fV\not\simeq V&\FPdim(V)^2&2\FPdim(V)^2
\end{array}
\end{equation}
Lemma~\ref{lem:fermion-parity-20261027} gives $2\mid\FPdim(V)^2$ in
the last row. Thus each entry in the middle column and the corresponding
entry in the last column have precisely the same prime divisors. By Lemma~\ref{lem:class-square-minimal-center-20261026}\textup{(1)}, the simple squares in the last column are the class dimensions of $\mathcal N$. The two columns therefore give the same vertices and edges.
\end{proof}

\begin{theorem}[The slightly degenerate two-component bound]
\label{thm:class-graph-slightly-degenerate-20261027}
Let $\C$ be a weakly integral slightly degenerate braided fusion
category over $\kk$. Then
$\Delta_{\mathrm{cl}}(\C)$ has at most two connected components.
If $\C$ contains a Tannakian fusion subcategory braided equivalent to
$\Rep(G)$ for a finite nonsolvable group $G$, this graph is connected.
\end{theorem}
\begin{proof}
Use Proposition~\ref{prop:class-graph-modular-realization-20261027}
and Theorem~\ref{thm:class-graph-modular-two-20261027}.
A Tannakian subcategory of $\C$ remains a Tannakian subcategory of the
chosen $\mathcal N$, so the connectedness assertion follows from
Proposition~\ref{prop:class-graph-tannakian-core-20261027}.
\end{proof}

\subsection{Centralizer orbit formulas and cyclic Tannakian center}
\begin{lemma}[Class dimensions in a Tannakian centralizer]
\label{lem:class-dimensions-centralizer-orbits-20261027}
Let $\mathcal M$ be a weakly integral modular fusion category over
$\kk$, with its canonical positive spherical structure. Let
$\Rep(G)\subseteq\mathcal M$ be a Tannakian fusion subcategory for a
finite group $G$, and put $\C=\Rep(G)'_{\mathcal M}$, its M\"uger
centralizer, equipped with the inherited braiding. This is a weakly
integral fusion subcategory of $\mathcal M$.
The de-equivariantization $\mathcal V=\mathcal M_G$ by this copy of $\Rep(G)$
is a weakly integral braided $G$-crossed fusion category. For $W\in\Irr(\mathcal M)$ denote by $\mu_W$
the restriction to $K(\C)$ of its Verlinde character
$[X]\mapsto S_{X,W}/\FPdim(W)$. Under $\mathcal M\simeq\mathcal V^G$,
choose a representative $V\in\Irr(\mathcal V)$ of the orbit
underlying $W$ and the corresponding irreducible projective representation
$\rho$ of its inertia group $G_V$. Set
\begin{equation}\label{explicit:80658146e8:eq0069}
 H=G_V,\qquad t=[G:H].
\end{equation}
Then
\begin{align}
 \FPdim(W)^2&=t^2\dim(\rho)^2\FPdim(V)^2,
 \label{explicit:80658146e8:eq0070}\\
 \FPdim((W\otimes W^*)\cap\Rep(G))&=t\dim(\rho)^2,
 \label{explicit:graphs:centralizer-internal-end}\\
 \FPdim(Z_{\mu_W}^{\C})&=t\,\FPdim(V)^2.
 \label{explicit:graphs:centralizer-class-dimension}
\end{align}
Every class dimension of $\C$ is obtained in this way.
\end{lemma}
\begin{proof}
The first identity follows by squaring the equivariantization dimension
formula~\eqref{eq:equivariantization-simple-dimension}.  For $U\in\Irr(G)$, the multiplicity of $U$ in $W\otimes
W^*$ is
\begin{equation}\label{explicit:80658146e8:eq0071}
 \dim\operatorname{Hom}_H
       (\operatorname{Res}^G_H U\otimes\rho,\rho).
\end{equation}
Since $\rho\otimes\rho^*$ is an ordinary $H$-representation, summing
with weights ${\FPdim(U)}$ gives
\begin{equation}\label{explicit:80658146e8:eq0072}
 \sum_U {\FPdim(U)}\,
 \dim\operatorname{Hom}_H
 (\operatorname{Res}^G_HU,\rho\otimes\rho^*)
 =\dim\operatorname{Ind}_H^G(\rho\otimes\rho^*)=t\dim(\rho)^2.
\end{equation}
The modular quotient formula~\eqref{eq:modular-quotient} of Theorem~\ref{thm:modular-quotient} and the double-centralizer identity~\eqref{eq:muger-centralizer-identities}
$\C'_{\mathcal M}=\Rep(G)$ now yield
$\FPdim(Z_{\mu_W}^{\C})={\FPdim(W)}^2/(t\dim(\rho)^2)=t{\FPdim(V)}^2$.  Finally, by Proposition~\ref{prop:extent}, every character of the
commutative semisimple algebra $K(\C)\otimes\mathbb C$ extends to a
character of $K(\mathcal M)\otimes\mathbb C$, and the latter are
Verlinde characters.
\end{proof}

\begin{theorem}[Cyclic Tannakian M\"uger center]
\label{thm:class-graph-cyclic-center-20261027}
Let $\C$ be a weakly integral braided fusion category over $\kk$
whose M\"uger center is Tannakian and braided equivalent to
$\Rep(C_n)$, where $C_n$ is cyclic of order $n\ge1$. Then there is a
weakly integral modular extension $\mathcal M$ containing $\C$ as a
full braided fusion subcategory such that
\begin{equation}\label{explicit:80658146e8:eq0073}
 \Delta_{\mathrm{cl}}(\C)=
 \Delta_{\mathrm{cl}}(\mathcal M).
\end{equation}
In particular, $\Delta_{\mathrm{cl}}(\C)$ has at most two components.
\end{theorem}
\begin{proof}
If $n=1$, take $\mathcal M=\C$. Suppose $n>1$. For a Tannakian
M\"uger center $\Rep(G)$, the obstruction to a minimal
nondegenerate extension lies in $H^4(G,\mathbb C^\times)$; see
\cite[Theorem~4.9]{GVR}. The cyclic cohomology formula~\eqref{eq:cyclic-even-cohomology} with $j=2$ gives $H^4(C_n,\mathbb C^\times)=0$, so a minimal extension $\mathcal M$ exists.
Its dimension is $n\FPdim(\C)$, so it is weakly integral, and
$\C=\Rep(C_n)'_{\mathcal M}$.

Apply the orbit formulas~\eqref{explicit:80658146e8:eq0070} and~\eqref{explicit:graphs:centralizer-class-dimension} of Lemma~\ref{lem:class-dimensions-centralizer-orbits-20261027}.
Every inertia subgroup is cyclic, and its projective irreducible
representations have dimension one by~\eqref{eq:cyclic-even-cohomology} with $j=1$.  Thus the corresponding
dimensions are
\begin{equation}\label{explicit:80658146e8:eq0074}
 {\FPdim(W)}^2=t^2{\FPdim(V)}^2,
 \qquad \FPdim(Z_{\mu_W}^{\C})=t{\FPdim(V)}^2.
\end{equation}
The number ${\FPdim(V)}^2$ is an integer because the de-equivariantization is
weakly integral.  These two numbers have identical prime supports.
All class dimensions of $\C$ occur, so the graphs agree.  Apply
Theorem~\ref{thm:class-graph-modular-two-20261027}.
\end{proof}

\begin{remark}[A prime-support test for a chosen extension]
For a noncyclic Tannakian center the fourth-cohomology obstruction may
be nonzero. Suppose a minimal modular extension exists. By Lemma~\ref{lem:class-dimensions-centralizer-orbits-20261027}, the dimensions in~\eqref{explicit:graphs:centralizer-class-dimension} and~\eqref{explicit:80658146e8:eq0070}, namely $t\FPdim(V)^2$ and $t^2\dim(\rho)^2\FPdim(V)^2$, have
the same prime support precisely when
\begin{equation}\label{eq:modular-extension-inertia-primes-20261102}
 \pi(\dim(\rho))\subseteq\pi\bigl(t\FPdim(V)^2\bigr).
\end{equation}
Thus the square of the projective inertia degree $\dim(\rho)^2$ is the
only factor that can introduce an additional prime in the ambient
simple square. If \eqref{eq:modular-extension-inertia-primes-20261102}
holds for every equivariant label $(V,\rho)$, the two graphs agree,
since every downstairs class and every upstairs simple square occurs in these orbit formulas. The general bound in Theorem~\ref{thm:general-two-components-20261102}
uses central de-equivariantization and also covers super-Tannakian
centers.
\end{remark}

\subsection{Class-image descent}
\subsubsection{Restricting class objects to a smaller braided category}

Restriction to a braided subcategory containing $Y_i$ detects when
the class image becomes transparent.

\begin{lemma}[Descent and collapse of a class]
\label{lem:class-restriction-collapse-20261027}
Let $\C$ be a weakly integral braided fusion category over $\kk$
with its canonical positive spherical structure, and let
$\mathcal L\subseteq\C$ be a full fusion subcategory with the induced
braiding. For a simple class summand $Z_i$ of $I_{\C}(\one)$ in
$Z(\C)$, suppose its forgetful image $Y_i=F(Z_i)$ belongs to
$\mathcal L$. Write $\mu_i$ for the corresponding character of $K(\C)$
and put $\nu_i=\mu_i|_{K(\mathcal L)}$.
Restricting the half-braiding of $Z_i$ to $\mathcal L$ produces an object
of $Z(\mathcal L)$ with exactly one canonical class constituent, of
multiplicity one, namely the class labeled by $\nu_i$. If $\FPdim(Z_{\nu_i}^{\mathcal L})=1$, then
\begin{equation}\label{explicit:80658146e8:eq0075}
 Y_i\in Z_2(\mathcal L).
\end{equation}
If $Z_j$ is another class summand with $Y_j=F(Z_j)\in\mathcal L$ and
$Z_i,Z_j$ centralize each other in $Z(\C)$, their descended
canonical class constituents centralize as well.
\end{lemma}
\begin{proof}
Restrict half-braidings on the full subcategory of $Z(\C)$ whose
underlying objects lie in $\mathcal L$. This gives a braided tensor
functor $R$ to $Z(\mathcal L)$. Adjunction and $[Y_i:\one]=1$ show
that $R(Z_i)$ contains exactly one canonical class constituent, once.
For $X\in\Irr(\mathcal L)$, preservation of the canonical positive
spherical structure gives
\begin{equation}\label{explicit:80658146e8:eq0076}
 \frac{S_{R(Z_i),X^+}}{\FPdim(Z_{i})}=\nu_i(X).
\end{equation}
Writing $R(Z_i)=\bigoplus_W n_WW$ expresses this character as the
positive weighted average of the restrictions of the Verlinde characters
of the $W$, with weights $n_W{\FPdim(W)}/\FPdim(Z_{i})$. Distinct algebra characters are
linearly independent, so every character in this average equals $\nu_i$.
This identifies the canonical constituent. Braided functors preserve
centralization, which also passes to constituents.

For the collapse assertion apply the modular quotient formula~\eqref{eq:modular-quotient} of Theorem~\ref{thm:modular-quotient} to
$\mathcal L^+\subseteq Z(\C)$ and the simple $Z_i$:
\begin{equation}\label{explicit:80658146e8:eq0077}
 \FPdim\bigl((Z_iZ_i^*)\cap(\mathcal L^+)'\bigr)
 =\frac{\FPdim(Z_{i})^2}{\FPdim(Z_{\nu_i}^{\mathcal L})}.
\end{equation}
If $\FPdim(Z_{\nu_i}^{\mathcal L})=1$, every constituent of $Z_iZ_i^*$ centralizes $\mathcal L^+$.
For every $X\prec Y_i$, the internal-end identity~\eqref{eq:offdiagonal} of Proposition~\ref{prop:diagonal} places $X^+$ in
$Z_iZ_i^*$; hence $X$ centralizes $\mathcal L$.
\end{proof}

\begin{remark}[The residual case in induction]
Suppose that the stronger cross-component vanishing assertion holds in
every proper fusion subcategory of $\C$. Let $Z_i,Z_j$ be a
counterexample, and let $\mathcal L$ be a proper fusion subcategory
containing $Y_i,Y_j$. Theorem~\ref{thm:braided-divisibility}, applied in the form~\eqref{eq:braided-divisibility}, gives
$\FPdim(Z_{\nu_i}^{\mathcal L})\mid \FPdim(Z_{i})$ and $\FPdim(Z_{\nu_j}^{\mathcal L})\mid \FPdim(Z_{j})$. If both $\FPdim(Z_{\nu_i}^{\mathcal L}),\FPdim(Z_{\nu_j}^{\mathcal L})>1$, their prime
supports remain in different components, while Lemma~\ref{lem:class-restriction-collapse-20261027} gives centralizing descended classes. This contradicts the inductive
hypothesis. Thus one of $\FPdim(Z_{\nu_i}^{\mathcal L}),\FPdim(Z_{\nu_j}^{\mathcal L})$ is one, and the corresponding entire
class image is transparent in $\mathcal L$. 
\end{remark}

\begin{proposition}[Arithmetic of a symmetric class image]
\label{prop:symmetric-class-image-stabilizer-20261027}
Let $\C$ be a weakly integral braided fusion category over $\kk$
whose Tannakian fusion subcategories are all pointed. Let $Z_i$ be
a simple class summand of $I_{\C}(\one)$ in $Z(\C)$ with forgetful
image $Y_i=F(Z_i)$. Assume that the fusion subcategory
$\langle Y_i\rangle\subseteq\C$ is symmetric for its induced braiding.
For an invertible $g\in\C$ write $g^+$ for its lift using that braiding. Let
\begin{equation}\label{explicit:80658146e8:eq0078}
 G_i=\{g\in\Irr(\C_{\pt}):g^+Z_i\simeq Z_i\}.
\end{equation}
Then precisely one of the following holds:
\begin{enumerate}[label=\textup{(\roman*)}]
\item $Y_i$ is pointed, $Y_i=\bigoplus_{g\in G_i}g$, and $\FPdim(Z_{i})=|G_i|$;
\item $Y_i$ is not pointed, $\FPdim(Z_{i})=2|G_i|$, and $4\mid \FPdim(Z_{i})$.
\end{enumerate}
In particular an odd-dimensional symmetric class image is pointed.
\end{proposition}
\begin{proof}
The internal-end identity~\eqref{eq:offdiagonal} of Proposition~\ref{prop:diagonal} and Schur's lemma give
$[Y_i:g]=1$ for $g\in G_i$ and zero for other invertible $g$.
This proves the pointed case.

In the non-pointed case write the symmetric category generated by $Y_i$
as $\Rep(G,z)$, where $z$ is a central involution giving the parity.
Its even Tannakian subcategory is pointed, so $G/\langle z\rangle$
is abelian. non-pointedness implies $[G,G]=\langle z\rangle$.
All linear representations are even, and every odd irreducible has
even degree: its determinant is trivial on the commutator $z$, whereas
$z$ acts on it as $-1$.

As a connected separable algebra in the underlying category $\Rep(G)$,
$Y_i$ is $\operatorname{Ind}_H^G\operatorname{End}(V)$ for a subgroup
$H$ and an irreducible projective $H$-representation $V$; see
\cite[Section~3.4]{OstrikModules}. The trace of $z$ on this algebra is
either zero or its full dimension. Indeed, if $z\notin H$ it permutes
the matrix blocks without fixed blocks. If $z\in H$, projective
commutation gives a character $a:H\to\kk^\times$ with
$\rho(t)\rho(z)\rho(t)^{-1}=a(t)\rho(z)$. If $a$ is nontrivial then
$\operatorname{Tr}\rho(z)=0$; otherwise Schur's lemma makes $\rho(z)$
scalar. Conjugation on $\operatorname{End}(V)$ therefore has trace zero
or full dimension, and the same holds on the induced algebra.
The full-dimension alternative would make every constituent even and
thus invertible, contrary to the assumption. Hence the even and odd
parts each have dimension $\FPdim(Z_{i})/2$. The even part consists exactly of
the invertible constituents, each once, giving $\FPdim(Z_{i})/2=|G_i|$. Its odd
part is a sum of even-dimensional simples, so $|G_i|$ is even.
\end{proof}

\begin{corollary}[Odd pointed class images]
\label{cor:odd-pointed-class-tannakian-20261027}
Let $\C$ be a weakly integral braided fusion category over $\kk$
with its canonical positive spherical structure. Let $Z_i$ be a simple
class summand of $I_{\C}(\one)$ in $Z(\C)$ whose forgetful image
$Y_i=F(Z_i)$ is a direct sum of invertibles in $\C$. If
$\FPdim(Z_i)$ is odd, then $Z_i$ centralizes itself in $Z(\C)$.
Consequently $\langle Z_i\rangle\subseteq Z(\C)$ is Tannakian and
has pointed forgetful image.
\end{corollary}
\begin{proof}
By Proposition~\ref{prop:symmetric-class-image-stabilizer-20261027},
$Y_i=\bigoplus_{g\in G_i}g$. The balancing identity~\eqref{eq:balancing-monodromy-trace} for $g^+Z_i\simeq Z_i$
gives $\mu_i(g)=\theta_g^{-1}$. Since $\mu_i$ is multiplicative, the
twist restricts to a character on $G_i$, so the corresponding pointed
subcategory is symmetric. A symmetric pointed category of odd order
has trivial twist. Thus $\mu_i(Y_i)=\FPdim(Z_{i})$, proving self-centralization.
The category generated by the twist-one self-centralizing simple
$Z_i$ is Tannakian. Its image is generated by the pointed object $Y_i$.
\end{proof}

\subsection{Disconnected graphs, Tannakian subcategories, and modular solvability}
\label{sec:modular-solvability-completion-20261027}

A disconnected class graph forces a nontrivial Tannakian subcategory.
In the modular case this gives a reduction by a cyclic group of prime
order. Throughout this section a disconnected graph has at least two
nonempty connected components.

\begin{proposition}[Disconnectedness forces a Tannakian subcategory]
\label{prop:disconnected-class-graph-tannakian-20261027}
Let $\C$ be a weakly integral braided fusion category over $\kk$.  If $\C$ has no
nontrivial Tannakian subcategory, then $\Dcl(\C)$ is empty or connected.
Equivalently, if $\Dcl(\C)$ has at least two nonempty connected
components, then $\C$ contains a nontrivial Tannakian fusion subcategory.
\end{proposition}
\begin{proof}
We first treat the case in which $\C$ is modular, and suppose for a
contradiction that its class graph is disconnected.

If noninvertible simples $X,Y$ have square dimensions belonging to
different components, then $S_{X,Y}=0$.  Indeed, their square dimensions
are coprime.  Lemma~\ref{lem:coprime-squares-monodromy-20261027} says that
a nonzero $S_{X,Y}$ would make them projectively centralize one another.
The orders of their invertible stabilizers divide ${\FPdim(X)}^2$ and ${\FPdim(Y)}^2$,
respectively, by \cite[Remark~9.4]{NPR}, so the stabilizers have trivial intersection.  By
\cite[Lemma~9.3]{NPR}, $X\otimes Y$ would be simple.  Its square dimension
${\FPdim(X)}^2{\FPdim(Y)}^2$ would join the two components, a contradiction.  We next remove a nondegenerate pointed factor.  The symmetric center of
$\C_{\pt}$ is either $\Vec$ or $\sVec$: any larger symmetric category has
a nontrivial even Tannakian subcategory.  If $\C_{\pt}$ is nondegenerate,
split it off.  Otherwise the decomposition of a slightly degenerate
pointed category in \cite[Corollary~A.19]{DGNO} writes it as a
nondegenerate pointed category times $\sVec$, and we split off the first
factor. By Proposition~\ref{prop:class-graph-product-20261026}, the pointed factor does not change the class graph.  Thus we may
assume
\begin{equation}\label{explicit:80658146e8:eq0079}
 \C_{\pt}=\Vec\quad\text{or}\quad\C_{\pt}=\sVec,
 \qquad \FPdim(\C_{\pt})\in\{1,2\}.
\end{equation}
The remaining modular factor still has no nontrivial Tannakian
subcategory.

If $\C_{\ad}$ is pointed, then $\C$ is braided nilpotent and its graph is
complete by Proposition~\ref{prop:class-graph-nilpotent-20261026}.  Hence
assume that $\C_{\ad}$ is non-pointed.  Choose any noninvertible
$Y\in\Irr(\C_{\ad})$ and choose a simple $X$ whose square dimension lies
in another component.  The category $\C_{\ad}$ is integral in the weakly
integral setting, so ${\FPdim(Y)}$ is an integer at least two.
We use the Verlinde computation of \cite[Remark~10.3]{NPR}. For integral
modular categories with trivial pointed part, it yields the connectedness
criterion of \cite[Proposition~10.5]{NPR}. Here $X$ may have nonintegral
dimension, and the pointed part may be $\sVec$. The Verlinde formula gives
\begin{equation}\label{explicit:80658146e8:eq0080}
 N_{X,X^*}^{Y}
 =\frac1{\FPdim(\C)}\sum_{T\in\Irr(\C)}
       \frac{S_{X,T}S_{X^*,T}S_{Y^*,T}}{{\FPdim(T)}}.
\end{equation}
Every summand with $T$ noninvertible vanishes: either its component differs
from that of $X$, or it differs from that of $Y$.  For invertible $T$,
$S_{X,T}S_{X^*,T}={\FPdim(X)}^2$, while $S_{Y^*,T}={\FPdim(Y)}$, because
$\C_{\ad}=(\C_{\pt})'$.  Therefore
\begin{equation}
 N_{X,X^*}^{Y}=\frac{\FPdim(\C_{\pt}) {\FPdim(X)}^2{\FPdim(Y)}}{{\FPdim(\C)}}.
 \label{eq:small-pointed-part-verlinde-20261027}
\end{equation}
Both ${\FPdim(X)}^2$ and ${\FPdim(Y)}^2$ are class dimensions and divide ${\FPdim(\C)}$.  Since they
are coprime, ${\FPdim(\C)}=t {\FPdim(X)}^2{\FPdim(Y)}^2$ for a positive integer $t$.  Consequently
\begin{equation}\label{explicit:80658146e8:eq0081}
 N_{X,X^*}^{Y}=\frac{\FPdim(\C_{\pt})}{t {\FPdim(Y)}}.
\end{equation}
For $\FPdim(\C_{\pt})=1$ this is impossible.  For $\FPdim(\C_{\pt})=2$, integrality forces ${\FPdim(Y)}=2$.
Since $Y$ was arbitrary, every noninvertible simple of $\C_{\ad}$ has
dimension two.

Let $f$ be the fermion when $\FPdim(\C_{\pt})=2$.  It belongs to $\C_{\ad}$ and is
transparent there.  It acts freely on its simples: $fY\simeq Y$ would
contradict the balancing identity~\eqref{eq:balancing-monodromy-trace}, because $\theta_f=-1$.  Thus $Y\otimes Y^*$ contains
the unit once and does not contain $f$.  Its remaining dimension is
$4-1=3$, whereas every remaining constituent has dimension two.  This
contradiction finishes the modular case.

Now let $\C$ be arbitrary with no nontrivial Tannakian subcategory.  Its
M\"uger center is $\Vec$ or $\sVec$.  Only the second case remains.  Choose
the graph-preserving minimal modular extension $\mathcal N$ supplied by
Proposition~\ref{prop:class-graph-modular-realization-20261027}, and write
$\C=\langle f\rangle'_{\mathcal N}$.  If $\mathcal N$ has no nontrivial
Tannakian subcategory, the modular case applies.  Otherwise choose one,
say $\mathcal T$.  The faithful $C_2$-grading of $\mathcal N$ with neutral
component $\C$ restricts to a grading of $\mathcal T$ with neutral
component $\mathcal T\cap\C=\Vec$.  Hence
$\mathcal T=\langle g\rangle\simeq\Rep(C_2)$, where $g$ lies outside
$\C$.  Its monodromy with $f$ is $-1$.

The pointed category $\mathcal P=\langle f,g\rangle$ has dimension four
and is nondegenerate: its generators have orders two, twists $-1,1$, and
mutual monodromy $-1$.  Thus
\begin{equation}\label{explicit:80658146e8:eq0082}
 \mathcal N=\mathcal P\boxtimes\mathcal B,
 \qquad
 \C=\langle f\rangle\boxtimes\mathcal B
       \simeq\sVec\boxtimes\mathcal B.
\end{equation}
Here $\mathcal B$ is modular and has no nontrivial Tannakian subcategory,
since it is a subcategory of $\C$.  The modular case applies to $\mathcal B$, and Proposition~\ref{prop:class-graph-product-20261026} identifies its graph with that of $\C$.
\end{proof}

\begin{lemma}[Free invertible actions and class divisibility]
\label{lem:free-invertible-action-class-divisibility-20261027}
Let $\mathcal A$ be a weakly integral braided fusion category over
$\kk$. For each algebra character $\mu$ of $K(\mathcal A)$, consider the
corresponding class object $Z_\mu^{\mathcal A}$ in $Z(\mathcal A)$.
If a subgroup
$H\subseteq\operatorname{Inv}(\mathcal A)$ acts freely on
$\Irr(\mathcal A)$ by tensor product, then
\begin{equation}\label{explicit:80658146e8:eq0083}
 \FPdim(Z_{\mu}^{\mathcal A})\mid\frac{\FPdim(\mathcal A)}{|H|}
 \qquad\text{for every character }\mu.
\end{equation}
In particular, let a cyclic group $C_p$ of prime order act by braided
tensor autoequivalences on a weakly integral braided fusion category
$\mathcal B$ over $\kk$. If this categorical action fixes every
simple isomorphism class of $\mathcal B$, then every character $\mu$ of
the complexified Grothendieck algebra $K(\mathcal B^{C_p})$ satisfies
\begin{equation}\label{explicit:80658146e8:eq0084}
 \FPdim(Z_{\mu}^{\mathcal B^{C_p}})\mid\FPdim(\mathcal B).
\end{equation}
\end{lemma}
\begin{proof}
For an invertible $a$, the number $\mu(a)$ is a root of unity.  Hence
$\mu(aX)\mu((aX)^*)=\mu(X)\mu(X^*)$.  Choose one representative in each
free $H$-orbit.  The formal codegree formula~\eqref{eq:traceexpansion} becomes
\begin{equation}\label{explicit:80658146e8:eq0085}
 f_\mu=|H|\sum_{X\in\Irr(\mathcal A)/H}\mu(X)\mu(X^*).
\end{equation}
Thus $f_\mu/|H|$ is an algebraic integer.  By~\eqref{eq:master}, it is the rational number
$\FPdim(\mathcal A)/(|H|\FPdim(Z_{\mu}^{\mathcal A}))$, so it is an integer.

For the final assertion, each simple of $\mathcal B$ has inertia group
$C_p$.  Its projective equivariant representations are all
one-dimensional, since $H^2(C_p,\kk^\times)=0$ by~\eqref{eq:cyclic-even-cohomology}.  Tensoring with the
canonical charges permutes its $p$ equivariant lifts freely.  Apply the
first assertion to these charges.
\end{proof}

\begin{lemma}[Fixing the noninvertibles in a modular category]
\label{lem:fixed-noninvertibles-modular-action-20261027}
Let $\mathcal D$ be a non-pointed weakly integral modular fusion category
over $\kk$. If a braided tensor autoequivalence of $\mathcal D$
fixes the isomorphism class of every noninvertible simple object,
then it fixes the isomorphism class of every invertible simple object.
\end{lemma}
\begin{proof}
Let $T$ be the autoequivalence, let $a$ be invertible, and let $X$ be
noninvertible.  The canonical positive spherical structure, and therefore
the $S$-matrix, is preserved by $T$.  Consequently
\begin{equation}\label{explicit:80658146e8:eq0086}
 S_{T(a),X}=S_{a,T^{-1}(X)}=S_{a,X}.
\end{equation}
The invertible object $T(a)a^{-1}$ thus centralizes every noninvertible
simple.  These simples generate $\mathcal D$: fix one noninvertible $X$;
for every invertible $b$, both $bX$ and $X^*$ are noninvertible, and their
product contains $b$.  Modularity forces $T(a)a^{-1}=\one$.
\end{proof}

\begin{proposition}[Components survive prime gauging]
\label{prop:class-graph-prime-gauging-components-20261027}
Let $\mathcal M$ be a weakly integral modular fusion category over
$\kk$ and let $\mathcal E\simeq\Rep(C_p)$ be a Tannakian fusion
subcategory of $\mathcal M$, with $p$ prime. De-equivariantize
$\mathcal M$ by $\mathcal E$, and let $\mathcal D=(\mathcal M_{C_p})_e$
be the weakly integral modular neutral component of the faithful $C_p$-grading on the
resulting braided $C_p$-crossed fusion category.  If
$\mathcal D$ is non-pointed, then every component of $\Dcl(\mathcal M)$
contains a vertex of $\Dcl(\mathcal D)$.  In particular,
\begin{equation}\label{explicit:80658146e8:eq0087}
 \#\pi_0\Dcl(\mathcal M)\le\#\pi_0\Dcl(\mathcal D).
\end{equation}
\end{proposition}
\begin{proof}
Set $\Pi=\operatorname{Vert}\Dcl(\mathcal D)$.  The iterated class It\^o--Michler decomposition~\eqref{explicit-alpha:a45450d4ea:eq0082} of Corollary~\ref{cor:im-class-support-20261024} gives
\begin{equation}\label{explicit:80658146e8:eq0088}
 \mathcal D=\mathcal P\boxtimes\mathcal B,
 \qquad\mathcal P\text{ pointed},\qquad
 \pi(\FPdim\mathcal B)=\Pi,
 \qquad\gcd(\FPdim\mathcal P,\FPdim\mathcal B)=1.
\end{equation}
The factor $\mathcal B$ is non-pointed.  These factors are preserved by the
$C_p$-action.  Indeed, $\mathcal P$ is the product of the full pointed
Sylow factors for the absent primes, hence is canonical; since
$\mathcal D$ is modular, its complementary factor is its centralizer.
Inside
$\mathcal E'_{\mathcal M}=\mathcal D^{C_p}$ put
\begin{equation}\label{explicit:80658146e8:eq0089}
 \mathcal A=\mathcal B^{C_p},\qquad
 \mathcal Q=\mathcal P^{C_p}.
\end{equation}
They centralize one another.  Their dimensions multiply to
$p^2\FPdim(\mathcal D)=\FPdim(\mathcal M)$, so modularity gives
\begin{equation}\label{explicit:80658146e8:eq0090}
 \mathcal A'_{\mathcal M}=\mathcal Q,
 \qquad \mathcal A\cap\mathcal Q=\mathcal E.
\end{equation}
The intersection assertion follows by forgetting equivariance: an object
belonging to both factors has underlying object in
$\mathcal P\cap\mathcal B=\Vec$.

Let $\Sigma$ be the union of the components of $\Dcl(\mathcal M)$ which
meet $\Pi$.  Every downstairs edge persists upstairs.  More explicitly,
the squares of neutral equivariant simples have the form
$t^2{\FPdim(X)}^2$, where $X\in\Irr(\mathcal D)$ and $t\in\{1,p\}$; projective
representations of cyclic inertia subgroups have degree one.

We claim that all neutral simple squares and all class dimensions of
$\mathcal A$ have prime support in $\Sigma$.  If $p\in\Sigma$, this
follows from the equivariantization dimension formula~\eqref{eq:equivariantization-simple-dimension} and
$\FPdim(Z_{\mu}^{\mathcal A})\mid p\FPdim(\mathcal B)$.  Suppose instead that
$p\notin\Sigma$.  A noninvertible $\mathcal D$-simple cannot have an
orbit of length $p$, since its equivariant square would connect $p$ to
$\Pi$.  Lemma~\ref{lem:fixed-noninvertibles-modular-action-20261027} then shows that the action fixes every
simple of $\mathcal D$.  Neutral simple dimensions are unchanged, and
the free-action divisibility~\eqref{explicit:80658146e8:eq0084} of Lemma~\ref{lem:free-invertible-action-class-divisibility-20261027} gives $\FPdim(Z_{\mu}^{\mathcal A})\mid\FPdim(\mathcal B)$.  This proves the claim.

Suppose a component of $\Dcl(\mathcal M)$ misses $\Pi$, and choose a
simple $W$ whose nontrivial square dimension belongs to that component.
The claim makes $W$ nonneutral. By the braided class divisibility~\eqref{eq:braided-divisibility} of Theorem~\ref{thm:braided-divisibility},
\begin{equation}\label{explicit:80658146e8:eq0091}
 \FPdim(Z_{\mu_W|_{\mathcal A}}^{\mathcal A})\mid {\FPdim(W)}^2.
\end{equation}
The two sides have disjoint prime support unless the left side is one.
The modular quotient formula~\eqref{eq:modular-quotient} of Theorem~\ref{thm:modular-quotient} therefore yields
\begin{equation}\label{explicit:80658146e8:eq0092}
 \FPdim((W\otimes W^*)\cap\mathcal A')={\FPdim(W)}^2,
 \qquad W\otimes W^*\in\mathcal Q.
\end{equation}

Choose a simple $Z\in\mathcal A$ over a noninvertible simple of
$\mathcal B$.  Its square dimension is nontrivial and supported in
$\Sigma$.  Every common simple constituent of $W^*W$ and $ZZ^*$ lies in
$\mathcal Q\cap\mathcal A=\mathcal E$.  No nontrivial charge $a$ can
stabilize $W$: if $aW\simeq W$, the balancing identity~\eqref{eq:balancing-monodromy-trace} and $\theta_a=1$ would make
$W$ centralize $a$, hence all of the prime-order category $\mathcal E$.
That would make $W$ neutral.  Thus the unit is the only common
constituent, with multiplicity one.  Rigidity gives
\begin{equation}\label{explicit:80658146e8:eq0093}
 \dim\End(W\otimes Z)=\dim\Hom(W^*W,ZZ^*)=1.
\end{equation}
The simple $W\otimes Z$ has square dimension ${\FPdim(W)}^2{\FPdim(Z)}^2$, joining the
allegedly distinct components.  This contradiction proves the first
assertion.  Because downstairs edges persist, the induced map from
downstairs components to upstairs components is surjective, proving the
component inequality.
\end{proof}

\begin{theorem}[Disconnected modular class graphs imply solvability]
\label{thm:disconnected-class-graph-modular-solvable-20261027}
Let $\mathcal M$ be a weakly integral modular fusion category over
$\kk$. If $\Dcl(\mathcal M)$ has at least two nonempty connected
components, then $\mathcal M$ is solvable.
\end{theorem}
\begin{proof}
We argue by induction on the integer $\FPdim(\mathcal M)$.
Proposition~\ref{prop:disconnected-class-graph-tannakian-20261027}
supplies a nontrivial Tannakian subcategory $\Rep(G)$.
Proposition~\ref{prop:class-graph-tannakian-core-20261027} excludes a
nonsolvable $G$.  A nontrivial solvable group has a quotient $C_p$, so
we obtain $\mathcal E=\Rep(C_p)\subseteq\mathcal M$.

De-equivariantize and put $\mathcal D=(\mathcal M_{C_p})_e$.  This is
modular, weakly integral, and has dimension
$\FPdim(\mathcal M)/p^2$.  If it is pointed, it is solvable.  If it is
non-pointed, Proposition~\ref{prop:class-graph-prime-gauging-components-20261027}
shows that its graph is disconnected, since the upstairs graph is
disconnected.  The induction hypothesis then makes $\mathcal D$
solvable.  In both cases, $\mathcal M_{C_p}$ is a $C_p$-extension of a
solvable category, and
$\mathcal M=(\mathcal M_{C_p})^{C_p}$ is solvable.  We use closure under
extensions and equivariantizations by solvable groups
\cite[Proposition~4.5]{ENO}.
\end{proof}

\begin{corollary}[Slightly degenerate and cyclic-center cases]
\label{cor:disconnected-class-graph-minimal-center-solvable-20261027}%
\label{thm:disconnected-class-graph-wgt-solvable-20261027}
Let $\C$ be a weakly integral braided fusion category over $\kk$.
Suppose that either $\C$ is slightly degenerate or its M\"uger center is
braided equivalent to the Tannakian category $\Rep(C_n)$ for a cyclic
group of order $n\ge1$.
If $\Dcl(\C)$ has at least two nonempty connected components, then $\C$
is solvable.
\end{corollary}
\begin{proof}
Proposition~\ref{prop:class-graph-modular-realization-20261027} and
Theorem~\ref{thm:class-graph-cyclic-center-20261027}, respectively,
provide a weakly integral modular extension with the same class graph.
Theorem~\ref{thm:disconnected-class-graph-modular-solvable-20261027} makes that extension solvable.  Its fusion
subcategory $\C$ is therefore solvable, by
\cite[Proposition~4.5]{ENO}.
\end{proof}

In particular, Corollary~\ref{cor:disconnected-class-graph-minimal-center-solvable-20261027} applies to weakly group-theoretical braided
categories with any of the stated M\"uger centers, since they are weakly
integral.

\subsection{Class dimensions from the reduced ring and the signed quotient}
\label{sec:fermionic-class-dimensions}

Let $\C$ be weakly integral and slightly degenerate, with transparent
fermion $f$. Choose a set $\Pi_0$ of representatives for the pairs
$\{X,fX\}$ in $\Irr(\C)$, with $\one\in\Pi_0$. For $\epsilon=\pm1$ write
\begin{equation}\label{eq:fermion-signed-quotients}
 R_\epsilon=K_0(\C)/([f]-\epsilon[\one]),\qquad
 x_X^\epsilon=[X]\bmod([f]-\epsilon[\one]) .
\end{equation}
The quotient $R_+$ is the reduced fusion ring used in
\cite[Proposition~2.4]{BDPReduced}. With the positive spherical structure,
its character table is given by the reduced $S$-matrix
\cite[Proposition~3.3]{YuReduced}. The classes detected by these
characters are also described by the coset dimension formula
\cite[Corollary~4.19 and Proposition~4.21]{BurciuCosets}.
The quotient $R_-$ retains the characters taking value $-1$ at $f$.
Its structure constants in the orbit basis can be negative.

\begin{proposition}[Class dimensions for the two fermion characters]
\label{prop:fermion-sign-class-dimensions}
Let $\C,f,\Pi_0$ and $R_\epsilon$ be as in
\eqref{eq:fermion-signed-quotients}. An algebra character $\lambda$ of
$R_\epsilon\otimes_{\mathbb Z}\kk$ lifts to a unique character $\mu$ of
$K(\C)$ with $\mu([f])=\epsilon$. As $\epsilon$ and $\lambda$ vary, every
conjugacy class occurs exactly once, with dimension
\begin{equation}\label{eq:fermion-reduced-class-dimension}
 \FPdim(Z_\mu^\C)=
 \frac{\FPdim(\C)}
      {2\sum_{X\in\Pi_0}|\lambda(x_X^\epsilon)|^2}.
\end{equation}
For $\epsilon=+1$ these dimensions are $\FPdim(X)^2$, one for each
$X\in\Pi_0$. Each sign accounts for $|\Pi_0|$ classes and for half
the total class dimension.
\end{proposition}
\begin{proof}
The balancing identity~\eqref{eq:balancing-monodromy-trace} gives
$\theta_{fX}=-\theta_X$ for every simple $X$, so the fermion acts freely.
The quotient $R_\epsilon$ is free over $\mathbb Z$ on
$\{x_X^\epsilon:X\in\Pi_0\}$, with multiplication
\begin{equation}\label{eq:fermion-signed-multiplication}
 x_X^\epsilon x_Y^\epsilon=
 \sum_{Z\in\Pi_0}(N_{X,Y}^Z+\epsilon N_{X,Y}^{fZ})x_Z^\epsilon.
\end{equation}
Indeed, multiplication by $[f]$ exchanges the two basis vectors in
each pair. Since $[f]^2=[\one]$, each algebra character factors through
exactly one quotient. Both complexified quotients are commutative
semisimple of dimension $|\Pi_0|$. On a pair $X,fX$ the two character
values differ by the factor $\epsilon$. Thus the full formal codegree
in~\eqref{eq:traceexpansion} is twice the orbit sum of squared absolute
values. The class-dimension formula~\eqref{eq:master} gives
\eqref{eq:fermion-reduced-class-dimension}.

For $Y\in\Pi_0$, the normalized reduced column sends $x_X^+$ to
$S_{X,Y}/\FPdim(Y)$. Reduced character orthogonality
\cite[Proposition~3.3]{YuReduced} gives
\begin{equation}\label{eq:fermion-reduced-positive-codegree}
 \sum_{X\in\Pi_0}\left|\frac{S_{X,Y}}{\FPdim(Y)}\right|^2
   =\frac{\FPdim(\C)}{2\FPdim(Y)^2}.
\end{equation}
These columns exhaust the characters of $R_+\otimes_{\mathbb Z}\kk$,
which proves the assertion for $\epsilon=+1$. Finally, apply the
fiber-extent formula~\eqref{eq:extent-general} of
Proposition~\ref{prop:extent} to $\langle f\rangle\subseteq\C$.
Each of the two characters of $\langle f\rangle$ has class dimension one,
so each fiber has total class dimension $\FPdim(\C)/2$.
\end{proof}

\begin{corollary}[When the two class-dimension multisets agree]
\label{cor:fermion-equal-multisets}
In Proposition~\ref{prop:fermion-sign-class-dimensions}, the two multisets
of class dimensions, counted with multiplicity, coincide if and only if
$f\notin\C_{\ad}$. In this case the full multiset is
$\{\FPdim(X)^2:X\in\Irr(\C)\}$, with one entry per simple object. In particular,
under this condition the reduced ring $R_+$ gives the full class prime graph.
\end{corollary}
\begin{proof}
If $f\notin\C_{\ad}$, its universal degree has order two in the finite
abelian group $U(\C)$. Choose a character $\gamma$ of $U(\C)$ taking
value $-1$ there. Twisting $\mu([X])$ by $\gamma(\deg X)$ gives a
bijection from the positive fiber to the negative fiber. It preserves
formal codegrees by~\eqref{eq:traceexpansion}. Conversely, the positive
fiber contains a class of dimension one, so equality of multisets gives
such a class in the negative fiber. By the characterization
\eqref{explicit:80658146e8:eq0098} in the proof of
Lemma~\ref{lem:nontrivial-class-character-rigidity-20261027}, its
character is $\FPdim$ twisted by a character of $U(\C)$. Its value
$-1$ at $f$ forces $\deg f\ne e$.
\end{proof}

\begin{example}[A class of dimension two in an integral fermion centralizer]
\label{ex:fermion-ising-cube-nonsquare-class}
Let $\mathcal I$ be a positive Ising modular category, with fermion
$\psi$ and noninvertible simple $\sigma$, and set
\begin{equation}\label{eq:ising-cube-fermion-centralizer}
 \mathcal M=\mathcal I^{\boxtimes3},\qquad
 f=\psi\boxtimes\psi\boxtimes\psi,\qquad
 \C=\langle f\rangle'_{\mathcal M}.
\end{equation}
The double-centralizer identity~\eqref{eq:muger-centralizer-identities}
gives $Z_2(\C)=\langle f\rangle\simeq\sVec$. Monodromy with $f$ acts
on a simple triple by $(-1)^k$, where $k$ counts its $\sigma$ factors.
Hence $\C$ has eight invertibles and six simples of dimension two, so
$\FPdim(\C)=32$ and the only simple squares are $1,4$. On fermion
orbits these contribute four classes of dimension one and three of
dimension four by Proposition~\ref{prop:fermion-sign-class-dimensions}.

To see a further class directly in the fusion ring, let $\eta$ be
the Ising character with $\eta([\psi])=-1$ and $\eta([\sigma])=0$.
Restrict $\eta\boxtimes\FPdim\boxtimes\FPdim$ to $K(\C)$. Its value
at $[f]$ is $-1$. On the eight invertibles its absolute value is one.
Of the six remaining simples, only those with $\sigma$ in the last two
positions contribute, with values $2,-2$. Thus its full formal codegree
is $8+2\cdot 2^2=16$, and the class-dimension formula~\eqref{eq:master}
gives $32/16=2$.

For completeness the negative fiber has six classes of dimension two
and one of dimension four. Indeed, put $\eta$ successively in each of
the three positions and choose either sign of $\sqrt2$ for the value
at $\sigma$ of one of the other two Ising characters with value $1$ at
$\psi$. This gives six distinct restrictions of codegree $16$.
The restriction of $\eta^{\boxtimes3}$ has codegree $8$, since it
vanishes on the six noninvertibles. These seven characters exhaust
the fiber by Proposition~\ref{prop:fermion-sign-class-dimensions}. In this example the
additional dimension $2$ has the same prime support as $4$, so both
fibers give the single vertex $2$.
\end{example}

\section{Equivariantization and the component bound}\label{sec:unified-equivar}

\subsection{Coprime actions and components under equivariantization}
\label{sec:equivariant-connectedness-20261027}

For a non-pointed base, we show that equivariantization induces a
surjection on the sets of connected components of the class graphs.
All group actions in this section are categorical actions by braided
tensor autoequivalences. We call an action trivial if it is equivalent
to the action with identity functors and identity coherence maps.

\subsubsection{Rank-one components and the size-one spectral subspace}

\begin{lemma}[A rank-one universal component forces solvability]
\label{lem:rank-one-universal-component-solvable-20261027}
Let $\mathcal B$ be a weakly integral braided fusion category over
$\kk$. If a homogeneous component of its universal grading has
exactly one isomorphism class of simple objects,
then $\mathcal B$ is solvable.
\end{lemma}
\begin{proof}
Let $\mathcal A=\mathcal B_{\ad}$ and suppose
$\Irr(\mathcal B_u)=\{X\}$.  For every $Y\in\Irr(\mathcal A)$ the object
$Y\otimes X$ belongs to $\mathcal B_u$, so comparison of dimensions gives
\begin{equation}\label{explicit:80658146e8:eq0094}
 Y\otimes X\simeq {\FPdim(Y)} X,\qquad
 X\otimes X^*=\bigoplus_{Y\in\Irr(\mathcal A)}{\FPdim(Y)} Y,
 \qquad {\FPdim(X)}^2=\FPdim(\mathcal A),
\end{equation}
where $\FPdim(Y)$ is consequently an integer.  The second equality
follows from the first by rigidity.

Equip $\mathcal B$ with its canonical positive ribbon structure.
The balancing identity~\eqref{eq:balancing-monodromy-trace} shows that the monodromy of $Y$ and $X$ is the scalar
$\theta_Y^{-1}$, because $Y\otimes X$ is a multiple of $X$.  Thus $X$ and
$Y$ projectively centralize each other.  By
\cite[Proposition~3.22]{DGNO}, $X\otimes X^*$ centralizes $Y$.
Its support is all of $\mathcal A$, so $\mathcal A$ is symmetric.

Write $\mathcal A\simeq\Rep(L,z)$, where $z$ is the central involution
determining its symmetric braiding, and allow $z=1$.  Its even Tannakian
subcategory is
\begin{equation}\label{explicit:80658146e8:eq0095}
 \mathcal E=\Rep(H),\qquad H=L/\langle z\rangle.
\end{equation}
Let $\mathcal T=\langle X\rangle$.  Since $XX^*$ contains every simple of
$\mathcal A$, one has $\mathcal T_{\ad}=\mathcal A$.  The universal
grading of $\mathcal T$ is cyclic, generated by the degree of $X$.
Every object of $\mathcal E$ has twist one, so the same balancing identity~\eqref{eq:balancing-monodromy-trace}
gives actual centralization of $X$.  Hence
$\mathcal E\subseteq\mathcal Z_2(\mathcal T)$.

De-equivariantize $\mathcal T$ by $\mathcal E$.  The cyclic grading
descends because the algebra of functions on $H$ lies in the neutral
component.  Its neutral component becomes $\Vec$ if $z=1$, and $\sVec$
otherwise.  Each homogeneous component therefore has dimension $1$ or
$2$, and contains at most two simple objects.

Consider the component corresponding to $X$.  Its equivariantization
has just one simple.  The parametrization of equivariant simple objects
\cite[Corollary~2.13]{BN} implies that its simples downstairs form one $H$-orbit, and
that the inertia subgroup $K$ of one of them has exactly one irreducible
projective representation for the relevant cocycle $\alpha$.  Thus
\begin{equation}\label{explicit:80658146e8:eq0096}
 [H:K]\le2,\qquad \kk^\alpha[K]\text{ is a matrix algebra}.
\end{equation}
The theorem on groups of central type implies that $K$ is solvable
\cite[Theorem~7.3]{HowlettIsaacs}.  Since $[H:K]\le2$, the group $H$ is
solvable, and so is its central extension $L$.  Hence $\mathcal A$ is
solvable.  Finally, the universal grading group of the braided category
$\mathcal B$ is abelian; $\mathcal B$ is an extension of $\mathcal A$ by
this group, and is solvable by \cite[Proposition~4.5]{ENO}.
\end{proof}

\begin{lemma}[Rigidity on the nontrivial class characters]
\label{lem:nontrivial-class-character-rigidity-20261027}
Let $\mathcal B$ be a weakly integral braided fusion category over
$\kk$, and let $T$ be a based-ring automorphism of $K_0(\mathcal B)$.
Suppose the induced action on the algebra characters of $K(\mathcal B)$
fixes every $\mu$ such that $\FPdim(Z_\mu^{\mathcal B})>1$.  Then
every simple moved by $T$ lies in a rank-one component of the universal
grading.  Moreover, $T$ cannot move a simple within its own component.
In particular, if $\mathcal B$ is nonsolvable, then $T$ fixes every
simple and every character.
\end{lemma}
\begin{proof}
Put $U=U(\mathcal B)$ and
\begin{equation}\label{explicit:80658146e8:eq0097}
 R_u=\sum_{X\in\Irr(\mathcal B_u)}{\FPdim(X)}[X].
\end{equation}
The class-dimension-one characters are exactly
\begin{equation}\label{explicit:80658146e8:eq0098}
 \mu_\chi([X])=\chi(u){\FPdim(X)},
 \qquad X\in\mathcal B_u,\quad \chi\in\widehat U.
\end{equation}
By~\eqref{eq:master} and the formal codegree formula~\eqref{eq:traceexpansion}, $\FPdim(Z_{\mu}^{\mathcal B})=1$ means
$\sum_X|\mu([X])|^2=\FPdim(\mathcal B)$.  Since
$|\mu([X])|\le {\FPdim(X)}$, equality holds termwise.  Multiplicativity and
equality in the triangle inequality show that the phases
$\mu([X])/{\FPdim(X)}$ multiply across every nonzero fusion coefficient.
They therefore define a character of the universal grading group.
The converse is immediate.

The primitive-idempotent formula~\eqref{eq:idempotents} gives
\begin{equation}\label{explicit:80658146e8:eq0099}
 e_\chi=\frac1{\FPdim(\mathcal B)}
       \sum_{u\in U}\chi(u)^{-1}R_u.
\end{equation}
Fourier inversion identifies their linear span with
$\operatorname{span}_{\kk}\{R_u:u\in U\}$.
For every simple $X$, the element $[T(X)]-[X]$ is annihilated by all
characters $\mu$ with $\FPdim(Z_{\mu}^{\mathcal B})>1$, and therefore belongs to this span.

Each $R_u$ has strictly positive coefficients on precisely the simples
of $\mathcal B_u$.  A nonzero difference of two basis elements in the
same component cannot be a multiple of $R_u$.  If the two elements
belong to distinct components, restricting their difference to either
component forces that component to have rank one.  This proves the
first two assertions.  The final assertion follows from
Lemma~\ref{lem:rank-one-universal-component-solvable-20261027}.
\end{proof}

\subsubsection{Two coprime-action lemmas}

\begin{lemma}[Coprime actions fixing all simples]
\label{lem:coprime-simple-fixed-action-trivial-20261027}
Let $\mathcal B$ be a weakly integral braided fusion category over
$\kk$, and let a finite group $H$ act categorically on $\mathcal B$
by braided tensor autoequivalences. Suppose
\begin{equation}\label{explicit:80658146e8:eq0100}
 \gcd(|H|,\FPdim(\mathcal B))=1
\end{equation}
and every simple isomorphism class of $\mathcal B$ is $H$-fixed.
Then the categorical action is equivalent to the trivial action.
\end{lemma}
\begin{proof}
First restrict to a cyclic subgroup $L\le H$.  Since every
simple is $L$-fixed and $H^2(L,\kk^\times)=0$ by~\eqref{eq:cyclic-even-cohomology}, each simple of $\mathcal B$
has $|L|$ inequivalent equivariant lifts, permuted freely by the canonical
charges $\widehat L$.  The free-orbit codegree formula~\eqref{explicit:80658146e8:eq0085} from the proof of Lemma~\ref{lem:free-invertible-action-class-divisibility-20261027} gives, for every character $\mu$ of
$K(\mathcal B^L)$,
\begin{equation}\label{explicit:80658146e8:eq0101}
 f_\mu=|L|\sum_{X\in\Irr(\mathcal B^L)/\widehat L}
                 \mu([X])\mu([X^*]).
\end{equation}
Thus $f_\mu/|L|$ is an algebraic integer.  Since
\begin{equation}\label{explicit:80658146e8:eq0102}
 \frac{f_\mu}{|L|}
   =\frac{\FPdim(\mathcal B)}{\FPdim(Z_{\mu}^{\mathcal B^L})}
   \in\mathbb Q,
\end{equation}
it is an integer.  Every prime dividing $|L|$ is therefore absent from
$\Dcl(\mathcal B^L)$.

The iterated class It\^o--Michler decomposition~\eqref{explicit-alpha:a45450d4ea:eq0082} of Corollary~\ref{cor:im-class-support-20261024} splits off the full pointed factors for
these primes:
\begin{equation}\label{explicit:80658146e8:eq0103}
 \mathcal B^L\simeq\mathcal P\boxtimes\mathcal N,
 \qquad
 \FPdim(\mathcal P)=|L|,\quad
 \FPdim(\mathcal N)=\FPdim(\mathcal B).
\end{equation}
All canonical charges have orders dividing $|L|$.  Their components in
$\mathcal N$ must be trivial, by coprimality.  Hence the specified
subcategory $\Rep(L)$ lies in $\mathcal P$, and equality of dimensions
gives $\mathcal P=\Rep(L)$.  This splitting is over the canonical
$\Rep(L)$.  De-equivariantization shows that the restricted categorical
$L$-action is equivalent to the trivial action.

In particular, each functor belonging to an element of $H$ is
monoidally isomorphic to the identity.  After choosing such
isomorphisms, the coherence of the action is represented by a class in
\begin{equation}\label{explicit:80658146e8:eq0104}
 H^2\bigl(H,\operatorname{Aut}_{\otimes}(\operatorname{Id}_{\mathcal B})\bigr).
\end{equation}
Here
$\operatorname{Aut}_{\otimes}(\operatorname{Id}_{\mathcal B})
 \simeq\widehat{U(\mathcal B)}$.
Its order divides $\FPdim(\mathcal B)$ because the universal grading
is faithful and its neutral component has integral Frobenius--Perron
dimension. For a finite group $H$ and a finite abelian $H$-module $A$
of order coprime to $|H|$, one has $H^n(H,A)=0$ for every $n>0$
\cite[Chap.~XII, \S2, Cor.~2.7, p.~237]{CE}. Here
$\gcd(|H|,\FPdim(\mathcal B))=1$, so the result applies to
$A=\operatorname{Aut}_{\otimes}(\operatorname{Id}_{\mathcal B})$.
Taking $n=2$ makes the coherence obstruction vanish.
\end{proof}

\begin{lemma}[Coprime actions fixing the nontrivial classes]
\label{lem:coprime-nontrivial-class-action-trivial-20261027}
Let $\mathcal B$ be a non-pointed weakly integral braided fusion category
over $\kk$, and let a finite group $H$ act categorically on
$\mathcal B$ by braided tensor autoequivalences. If
\begin{equation}\label{explicit:80658146e8:eq0105}
 \gcd(|H|,\FPdim(\mathcal B))=1
\end{equation}
and $H$ fixes every character $\mu$ of $K(\mathcal B)$ with
$\FPdim(Z_\mu^{\mathcal B})>1$,
then the categorical action is equivalent to the trivial action.
\end{lemma}
\begin{proof}
Set $\mathcal A=\mathcal B_{\ad}$ and $U=U(\mathcal B)$.
Since $\mathcal B$ is non-pointed, $\mathcal A\ne\Vec$.
Lemma~\ref{lem:nontrivial-class-character-rigidity-20261027} shows that
every simple in a universal component of rank greater than one is
fixed.  In particular, all simples of $\mathcal A$ are fixed.
Lemma~\ref{lem:coprime-simple-fixed-action-trivial-20261027} trivializes
the restricted action on $\mathcal A$.

The universal grading determines a homomorphism
\begin{equation}\label{explicit:80658146e8:eq0106}
 \phi:U\longrightarrow\operatorname{BrPic}(\mathcal A),
 \qquad u\longmapsto[\mathcal B_u].
\end{equation}
Indeed, the components are invertible $\mathcal A$-bimodule categories
and multiplication induces equivalences
$\mathcal B_u\boxtimes_{\mathcal A}\mathcal B_v\simeq\mathcal B_{uv}$;
see \cite[Theorem~6.1 and Corollary~6.2]{ENO-homotopy}.
Let $K=\ker\phi$.

Since each acting functor is monoidally isomorphic to the identity on
$\mathcal A$, its restriction
$\mathcal B_u\to\mathcal B_{h(u)}$ is an equivalence of
$\mathcal A$-bimodule categories.  Hence $\phi(h(u))=\phi(u)$, and $H$
acts trivially on $U/K$.  For $k\in K$, the component $\mathcal B_k$ is
equivalent to the regular $\mathcal A$-bimodule category, so has rank
$\operatorname{rank}(\mathcal A)>1$.  The rigidity assertion in Lemma~\ref{lem:nontrivial-class-character-rigidity-20261027}
shows that $H$ fixes its simples, hence fixes $K$ pointwise.

The group $U$ is abelian.  Its automorphisms inducing the identity on
both $K$ and $U/K$ identify with
\begin{equation}\label{explicit:80658146e8:eq0107}
 \operatorname{Hom}(U/K,K),
 \qquad
 f\longmapsto \operatorname{id}_U+\iota f\pi,
\end{equation}
using additive notation, where $\iota:K\to U$ and $\pi:U\to U/K$ are
the canonical maps.  Its order has only primes dividing $|U(\mathcal B)|$.
Since $|U(\mathcal B)|\mid\FPdim(\mathcal B)$, coprimality forces the $H$-action on
$U$ to be trivial.  Every rank-one component is therefore preserved and
its unique simple is fixed.  The other simples were already fixed.
A final application of
Lemma~\ref{lem:coprime-simple-fixed-action-trivial-20261027}
trivializes the action on $\mathcal B$.
\end{proof}

\subsubsection{Components under equivariantization}

\begin{proposition}[Components under prime cyclic equivariantization]
\label{prop:class-graph-prime-equivariantization-20261027}
Let $p$ be a prime, and let $\mathcal B$ be a non-pointed weakly
integral braided fusion category over $\kk$. Suppose $C_p$ acts
categorically on $\mathcal B$ by braided tensor autoequivalences. The inclusion
$\Dcl(\mathcal B)\subseteq\Dcl(\mathcal B^{C_p})$ induces a surjection
on connected components. In particular,
\begin{equation}\label{eq:prime-equivariantization-components-20261102}
 \#\pi_0\Dcl(\mathcal B^{C_p})\leq\#\pi_0\Dcl(\mathcal B).
\end{equation}
\end{proposition}
\begin{proof}
Write $\mathcal A=\mathcal B^{C_p}$, and let $a$ generate the group of
canonical invertible charges.  Characters of $\Rep(C_p)$ are evaluations
at group elements.  We call the evaluation at the identity the identity
fiber.

We first record the arithmetic of the other fibers.  If $\mu(a)\ne1$
and $aX\simeq X$, then
$\mu([X])=\mu(a)\mu([X])$, so $\mu([X])=0$.
Every other charge orbit has length $p$.  Since
$\mu(a)\mu(a^*)=1$, the summand
$\mu([X])\mu([X^*])$ in the formal codegree formula~\eqref{eq:traceexpansion} is constant on a charge
orbit.  It follows that $f_\mu/p$ is an algebraic integer.  Therefore
\begin{equation}
 \label{eq:prime-nonidentity-fiber-divisibility-20261027}
 \FPdim(Z_{\mu}^{\mathcal A})\mid\FPdim(\mathcal B)
 \qquad\text{if }\mu|_{\Rep(C_p)}\ne\FPdim.
\end{equation}
For the identity fiber, the orbit formula~\eqref{explicit:80658146e8:eq0014} of Lemma~\ref{lem:class-equivariantization-20261026} gives
\begin{equation}\label{explicit:80658146e8:eq0108}
 \FPdim(Z_{\nu F}^{\mathcal A})
 =[C_p:(C_p)_\nu]\FPdim(Z_{\nu}^{\mathcal B}),
 \qquad [C_p:(C_p)_\nu]\in\{1,p\}.
\end{equation}
Thus the inclusion of graphs preserves every vertex and edge. We must
show that each component of $\Dcl(\mathcal A)$ contains a vertex of
$\Dcl(\mathcal B)$. The class It\^o--Michler theorem (Theorem~\ref{thm:main}), in the iterated form~\eqref{explicit-alpha:a45450d4ea:eq0082} of Corollary~\ref{cor:im-class-support-20261024}, gives
\begin{equation}\label{explicit:80658146e8:eq0109}
 \mathcal B=\mathcal P\boxtimes\mathcal D,\qquad
 \gcd(\FPdim(\mathcal P),\FPdim(\mathcal D))=1,\qquad
 \pi(\FPdim(\mathcal D))=\Pi,
\end{equation}
where $\mathcal P$ is pointed and
$\Pi=\operatorname{Vert}\Dcl(\mathcal B)
     =\operatorname{Vert}\Dcl(\mathcal D)$.
Since $\mathcal B$ is non-pointed, so is $\mathcal D$, and $\Pi$ is nonempty.
Both factors are preserved by the action.  For $\mathcal P$, this follows
because its invertibles form the full Hall subgroup of invertibles for
the primes absent from $\Pi$.  To check $\mathcal D$, project the
$\mathcal P$-coordinate of each simple image $T(X)$, $X\in\Irr(\mathcal D)$.
Tensor compatibility makes these coordinates a grading of $\mathcal D$
by a subgroup of $\operatorname{Inv}(\mathcal P)$.  Its order divides both
$|U(\mathcal D)|$ and $\FPdim(\mathcal P)$, so it is trivial.
Thus $T(\mathcal D)=\mathcal D$.  The multiplication equivalence
$\mathcal P\boxtimes\mathcal D\to\mathcal B$ is equivariant through the
tensorators of the braided action.

Put
\begin{equation}\label{explicit:80658146e8:eq0110}
 \mathcal A_P=\mathcal P^{C_p},
 \qquad \mathcal A_D=\mathcal D^{C_p}.
\end{equation}
By \eqref{explicit:80658146e8:eq0108}, the only possible new prime in
the identity fiber is $p$. Assume $p\notin\Pi$. If a character of
$\mathcal P$ moves, pair it with any character of $\mathcal D$ of class
dimension greater than one. The resulting character of $\mathcal B$
has orbit size $p$, so its class dimension upstairs is divisible by
$p$ and by a prime in $\Pi$. If instead a character of $\mathcal D$
of class dimension greater than one moves, pairing it with the
Frobenius--Perron character of $\mathcal P$ has the same effect.
In the remaining case all $\mathcal P$-characters and all
$\mathcal D$-characters of class dimension greater than one are fixed.
Since $p\nmid\FPdim(\mathcal D)$,
Lemma~\ref{lem:coprime-nontrivial-class-action-trivial-20261027}
trivializes the action on $\mathcal D$. Thus every pair of characters
is fixed and the identity fiber introduces no new prime.

It remains to treat a nonidentity charge value $\zeta$.  Write
$K(\mathcal A)_\zeta=e_\zeta K(\mathcal A)$, where
\begin{equation}\label{explicit:80658146e8:eq0111}
 e_\zeta=\frac1p\sum_{j=0}^{p-1}\zeta^{-j}[a^j].
\end{equation}
A simple lying over a free underlying $C_p$-orbit is fixed by all charges,
so its image under $e_\zeta$ vanishes.  A fixed underlying simple has
$p$ equivariant lifts by the cyclic cohomology formula~\eqref{eq:cyclic-even-cohomology}. Their images span
one line.  These lines, one for each fixed underlying simple, give a
basis of the fiber algebra.  A pair of underlying simples of
$\mathcal P$ and $\mathcal D$ is fixed if and only if both are fixed.
The product of chosen equivariant lifts from the two factors gives a
lift of the corresponding pair. Multiplication therefore takes the
tensor products of these basis vectors to a basis and induces an algebra
isomorphism
\begin{equation}\label{explicit:80658146e8:eq0112}
 K(\mathcal A_P)_\zeta\otimes_{\kk}K(\mathcal A_D)_\zeta
 \ \xrightarrow{\ \sim\ }\ K(\mathcal A)_\zeta.
\end{equation}
The units are the corresponding charge idempotents.

Let $\mu$ correspond to $(\mu_P,\mu_D)$.  Choose one equivariant lift of
each fixed underlying simple.  The charge-orbit computation gives
\begin{equation}\label{explicit:80658146e8:eq0113}
 f_\mu
 =p\!\left(\sum_{X\in\Irr(\mathcal P)^{C_p}}|\mu_P(\widetilde X)|^2\right)
     \left(\sum_{Y\in\Irr(\mathcal D)^{C_p}}|\mu_D(\widetilde Y)|^2\right)
 =\frac{f_{\mu_P}f_{\mu_D}}p.
\end{equation}
Comparison with the three global dimensions yields the exact formula
\begin{equation}
 \label{eq:prime-equivariant-fiber-product-20261027}
 \FPdim(Z_{\mu}^{\mathcal A})
 =\FPdim(Z_{\mu_P}^{\mathcal A_P})\,\FPdim(Z_{\mu_D}^{\mathcal A_D}).
\end{equation}

Every such fiber of $\mathcal A_D$ contains a class of dimension greater
than one.  Indeed, its total class mass is $\FPdim(\mathcal D)$ by the
extent formula~\eqref{eq:extent-general} of Proposition~\ref{prop:extent}, while its number of characters equals the number of
fixed simple objects of $\mathcal D$, which is at most
$\operatorname{rank}\mathcal D<\FPdim(\mathcal D)$.
Moreover,
\eqref{eq:prime-nonidentity-fiber-divisibility-20261027}
puts all its class primes in $\Pi$.
Pair any $\mu_P$ with such a $\mu_D$ in
\eqref{eq:prime-equivariant-fiber-product-20261027}.
Each prime of $\FPdim(Z_{\mu_P}^{\mathcal A_P})$ thus lies in a
component meeting $\Pi$, while every prime from $\mathcal A_D$ already
belongs to $\Pi$. This proves the assertion in every fiber.
\end{proof}

\begin{theorem}[Components under arbitrary equivariantization]
\label{thm:class-graph-connected-equivariantization-20261027}
Let $\mathcal B$ be a non-pointed weakly integral braided fusion
category over $\kk$, and let a finite group $G$ act categorically
on $\mathcal B$ by braided tensor autoequivalences. The inclusion
$\Dcl(\mathcal B)\subseteq\Dcl(\mathcal B^G)$ induces a surjection
on connected components. Consequently,
\begin{equation}\label{eq:equivariantization-components-20261102}
 \#\pi_0\Dcl(\mathcal B^G)\leq\#\pi_0\Dcl(\mathcal B).
\end{equation}
In particular, if $\Dcl(\mathcal B)$ is connected, then so is
$\Dcl(\mathcal B^G)$.
\end{theorem}
\begin{proof}
The identity-fiber formula~\eqref{explicit:80658146e8:eq0014} of Lemma~\ref{lem:class-equivariantization-20261026} gives the inclusion of graphs. For a normal subgroup $N\lhd G$, equivariantization in stages
gives the induced $G/N$-action on $\mathcal B^N$ and a braided
equivalence $(\mathcal B^N)^{G/N}\simeq\mathcal B^G$,
as in \cite[\S3]{GalindoJones}. Each intermediate equivariantization is
non-pointed because its dominant forgetful functor has non-pointed image.
For solvable $G$, iteration of
Proposition~\ref{prop:class-graph-prime-equivariantization-20261027}
along a composition series proves the assertion.

Suppose now that $G$ is centerless and nonsolvable. Its ordinary
class graph is connected, and every prime dividing $|G|$ is a vertex:
a prime absent from this graph would have central Sylow subgroup
\cite[Remark~2.4 and Theorem~2.5]{NPR}.
By the individual braided divisibility~\eqref{eq:braided-divisibility} of Theorem~\ref{thm:braided-divisibility}, the ordinary class graph
lies in a single component $\Sigma$ of $\Dcl(\mathcal B^G)$.
Consider any class of dimension greater than one lying outside $\Sigma$.
Its restriction to $\Rep(G)$ has class dimension one by the same
divisibility. It is evaluation at a central element, which must be the
identity. Therefore, by the identity-fiber formula~\eqref{explicit:80658146e8:eq0014}, there is a character $\nu$ of $K(\mathcal B)$ such that
\begin{equation}\label{explicit:80658146e8:eq0114}
 \FPdim(Z_{\mu}^{\mathcal B^G})=[G:G_\nu]\FPdim(Z_{\nu}^{\mathcal B}).
\end{equation}
The left side is coprime to $|G|$, so the orbit index is one. Every
vertex outside $\Sigma$ therefore belongs to $\Dcl(\mathcal B)$.

It remains to show that $\Sigma$ meets
$\Pi=\operatorname{Vert}\Dcl(\mathcal B)$. Otherwise $\Pi$ is disjoint
from $\pi(|G|)$, and all characters of $\mathcal B$ of class dimension
greater than one are $G$-fixed. Indeed, a nontrivial orbit index in
\eqref{explicit:80658146e8:eq0114} would join a prime of $|G|$ to a
prime in $\Pi$. As in \eqref{explicit:80658146e8:eq0109}, apply Corollary~\ref{cor:im-class-support-20261024} to split off the
pointed factor for the absent primes to obtain
\begin{equation}\label{explicit:80658146e8:eq0115}
 \mathcal B=\mathcal P\boxtimes\mathcal D,\qquad
 \pi(\FPdim(\mathcal D))=\Pi.
\end{equation}
The factors are $G$-stable by the argument following
\eqref{explicit:80658146e8:eq0109}. The category $\mathcal D$ is
non-pointed and $\gcd(|G|,\FPdim(\mathcal D))=1$. Pairing a character
of $\mathcal D$ of class dimension greater than one with the dimension
character of $\mathcal P$ shows that the former is fixed.
Lemma~\ref{lem:coprime-nontrivial-class-action-trivial-20261027}
trivializes the action on $\mathcal D$. The equivariance of the product
decomposition now gives
\begin{equation}\label{explicit:80658146e8:eq0116}
 \mathcal B^G\simeq\mathcal D\boxtimes\mathcal P^G.
\end{equation}
Both factors are non-pointed since $\mathcal P^G$ contains $\Rep(G)$.
Proposition~\ref{prop:class-graph-product-20261026} makes the product
graph connected. Its unique component contains the nonempty set $\Pi$,
contrary to the assumption on $\Sigma$.

For arbitrary nonsolvable $G$, the hypercenter $N=Z_\infty(G)$ is
nilpotent and $G/N$ is centerless and nonsolvable. Apply the solvable
case to $\mathcal B^N$ over $\mathcal B$, then the centerless case to
$(\mathcal B^N)^{G/N}$ over $\mathcal B^N$. The two surjections on
components compose to give the assertion.
\end{proof}

\subsection{All equivariant character fibers and pointed bases}
\label{sec:all-equivariant-fibers-20261027}

For the fiber over the conjugacy class of $g$, the relevant action is
the action of $C_G(g)$ on the cyclic fiber algebra. It incorporates
the coherence maps of the categorical action.

\begin{proposition}[The all-fiber equivariantization formula]
\label{prop:equivariant-class-fiber-formula-20261027}
Let $\mathcal D$ be a weakly integral braided fusion category over
$\kk$, and let a finite group $G$ act categorically on $\mathcal D$
by braided tensor autoequivalences. Put $\mathcal C=\mathcal D^G$.
Fix $g\in G$, set $L=\langle g\rangle$ and $Q=C_G(g)$.
There is a canonical $Q$-action on the $g$-fiber algebra of
$K(\mathcal D^L)$.  Characters in the conjugacy-$[g]$ fiber of
$K(\mathcal C)$ correspond to $Q$-orbits of characters $\nu$ in that
algebra.  For the corresponding character $\mu$,
\begin{equation}\label{explicit:80658146e8:eq0117}
 \boxed{
 \FPdim(Z_{\mu}^{\mathcal C})
 =|g^G|\,[Q:Q_\nu]\,\FPdim(Z_{\nu}^{\mathcal D^L}).}
\end{equation}
Here the $g$-fiber means restriction to evaluation at $g$ on
$K(\Rep(L))$, and the $[g]$-fiber means restriction to evaluation on
the conjugacy class of $g$ on $K(\Rep(G))$.
\end{proposition}
\begin{proof}
Regard the braided action as a crossed braided category with trivial
grading.  The extended Verlinde algebra of \cite{Kirillov} has the decomposition
\begin{equation}\label{explicit:80658146e8:eq0118}
 \widetilde K(\mathcal D)=\bigoplus_{x\in G}T_x(\mathcal D),
 \qquad
 T_x(\mathcal D)
 =\bigoplus_{\substack{X\in\Irr(\mathcal D)\\xX\simeq X}}
       \Hom(X,xX).
\end{equation}
The morphism spaces are understood in the usual twisted-trace
description, so the construction is independent of the chosen simple
representatives.  Tensor multiplication is zero between distinct
$x$-summands.  Each $T_x$ is commutative here, because the original
crossed grading is trivial.  The group action transports $T_x$ to
$T_{hxh^{-1}}$.

By \cite[Corollary~8.18]{Kirillov},
\begin{equation}\label{explicit:80658146e8:eq0119}
 K(\mathcal D^G)\simeq\widetilde K(\mathcal D)^G.
\end{equation}
Put $x=g^{-1}$. In the map of \cite[(8.14)]{Kirillov} the
equivariant structure map $gX\to X$ contributes to $T_x(\mathcal D)$.
Thus a charge representation acts on $T_x$ by its character evaluated
at $g$. Under \eqref{explicit:80658146e8:eq0119} the block over evaluation
at $[g]$ identifies with $T_x(\mathcal D)^Q$, by taking the $x$-component
of an invariant element supported on the conjugacy class of $x$.

The cyclic group $L=\langle x\rangle$ acts trivially on $T_x$. Indeed,
the twisted-trace relation, with the structural identifications of the
action understood, is
\begin{equation}\label{explicit:80658146e8:eq0120}
 [\,x(f)\varphi f^{-1}\,]=[\varphi]
 \quad\text{for }f:X\xrightarrow{\sim}Y,\quad
 \varphi:X\longrightarrow xX.
\end{equation}
On a fixed simple $X$, any nonzero $\varphi$ is an isomorphism. Taking
$f=\varphi$ gives $[x(\varphi)]=[\varphi]$, and these morphisms span
$T_x$. Thus the cyclic version of the same identification is
\begin{equation}\label{explicit:80658146e8:eq0121}
 T_{g^{-1}}(\mathcal D)\simeq K(\mathcal D^L)_g.
\end{equation}
In particular $T_x$ is a finite-dimensional commutative semisimple
algebra. The primitive idempotents of $T_x^Q$ are the sums over
$Q$-orbits of primitive idempotents of $T_x$.

Let $\sigma_x$ be the coefficient of the canonical unit morphism
$\one\to x\one$ in $T_x$. For a fusion category $\mathcal A$, let
$\tau_{\mathcal A}$ be the coefficient of the unit in $K(\mathcal A)$.
The trace identity and inverse map of
\cite[Theorem~8.17(2)--(3) and Corollary~8.18]{Kirillov} give
\begin{equation}\label{explicit:80658146e8:eq0122}
 \tau_{\mathcal C}|_{[g]}=\frac1{|Q|}\sigma_x,\qquad
 \tau_{\mathcal D^L}|_g=\frac1{|L|}\sigma_x .
\end{equation}
For the first equality the invariant element has $|g^G|$ conjugate
components with equal unit coefficients, and the inverse map divides
their total coefficient by $|G|$. For the second the conjugacy class
in $L$ has one element.

For a primitive idempotent $e_\nu$ of $T_x$ this yields
\begin{equation}\label{explicit:80658146e8:eq0123}
 \sigma_x(e_\nu)
 =\frac{|L|}{f_\nu}
 =\frac{\FPdim(Z_{\nu}^{\mathcal D^L})}{\FPdim(\mathcal D)}.
\end{equation}
The functional $\sigma_x$ is $Q$-invariant. If the orbit of $\nu$ has
length $t=[Q:Q_\nu]$, the sum of its primitive idempotents has
$\mathcal C$-trace
\begin{equation}\label{explicit:80658146e8:eq0124}
 \frac{t}{|Q|}
 \frac{\FPdim(Z_{\nu}^{\mathcal D^L})}{\FPdim(\mathcal D)}.
\end{equation}
Multiplying by $\FPdim(\mathcal C)=|G|\FPdim(\mathcal D)$ converts the
trace to the upstairs class dimension, giving the all-fiber formula~\eqref{explicit:80658146e8:eq0117}.
For $g=1$, it reduces to the identity-fiber formula~\eqref{explicit:80658146e8:eq0014} of Lemma~\ref{lem:class-equivariantization-20261026}.
\end{proof}

\begin{proposition}[Pointed bases and perfect acting groups]
\label{prop:pointed-perfect-equivariantization-connected-20261027}
Let $\mathcal D$ be a pointed braided fusion category over $\kk$,
and let a finite group $G$ act categorically on $\mathcal D$ by braided
tensor autoequivalences. Fix $g\in G$. Write $A=\operatorname{Inv}(\mathcal D)$ and
$A^g=\{a\in A:g(a)=a\}$. Let $\mu$ correspond to the
$C_G(g)$-orbit of a character $\nu$ in the $g$-fiber of
$K(\mathcal D^{\langle g\rangle})$ as in
Proposition~\ref{prop:equivariant-class-fiber-formula-20261027}. Then
\begin{equation}\label{explicit:80658146e8:eq0125}
 \boxed{
 \FPdim(Z_{\mu}^{\mathcal D^G})
 =|g^G|\,[C_G(g):C_G(g)_\nu]\,|A:A^g|.}
\end{equation}
If $G$ is nontrivial and perfect, then $\Dcl(\mathcal D^G)$ is connected.
More generally, the same conclusion holds whenever the ordinary class
graph of $G$ is nonempty, connected, and has vertex set $\pi(|G|)$.
\end{proposition}
\begin{proof}
Put $L=\langle g\rangle$, and take a character $\nu$ in the
generator-$g$ fiber of $K(\mathcal D^L)$.
A simple lying over an underlying element $a\in A$ not fixed by $g$
has proper inertia group $L_a<L$.  Choose a character $\chi$ of $L$
trivial on $L_a$ and with $\chi(g)\ne1$.  Tensoring by its canonical
charge fixes that simple, so its value under $\nu$ is zero.

Each $a\in A^g$ has exactly $|L|$ equivariant lifts since
$H^2(L,\kk^\times)=0$ by~\eqref{eq:cyclic-even-cohomology}.  These lifts are invertible; their character
values therefore have absolute value one.  The formal-codegree sum~\eqref{eq:traceexpansion} and the class-dimension formula~\eqref{explicit:80658146e8:eq0003} give
\begin{equation}\label{explicit:80658146e8:eq0126}
 f_\nu=|L||A^g|,\qquad
 \FPdim(Z_{\nu}^{\mathcal D^L})=\frac{|L||A|}{|L||A^g|}=|A:A^g|.
\end{equation}
Substitution in the all-fiber formula~\eqref{explicit:80658146e8:eq0117} proves the first assertion.

Suppose the group class graph is nonempty and connected with vertex
set $\pi(G)$.
Braided divisibility~\eqref{eq:braided-divisibility} of Theorem~\ref{thm:braided-divisibility} places this graph in one upstairs component
$\Sigma$.  Let $q$ be a prime in another component.  Then $q\nmid|G|$,
and the pointed-base formula~\eqref{explicit:80658146e8:eq0125} shows that $q$ could occur only in the factor
$|A:A^g|$.

For every noncentral $g$, the factor $|g^G|>1$ has a prime divisor in
$\pi(G)$.  Each corresponding fiber is nonzero, so~\eqref{explicit:80658146e8:eq0125} would
join $q$ to $\Sigma$ if $q\mid|A:A^g|$.  Thus every noncentral $g$ fixes
the full $q$-Sylow subgroup $A_q$ pointwise.  Noncentral elements
generate a nonabelian group: for $z\in Z(G)$ and $x\notin Z(G)$, both
$x$ and $xz$ are noncentral and $z=x^{-1}(xz)$.
Consequently all of $G$ fixes $A_q$.  Then $q$ divides none of the
indices $|A:A^g|$, and it divides neither of the group-theoretic factors.
It is absent from every class dimension, contrary to its choice.

A nontrivial perfect group is nonsolvable, so its ordinary class graph
is connected by Theorem~\ref{thm:group-class-graph-classical-20261027}.  If a prime divisor of $|G|$ were absent from that graph,
its Sylow subgroup would be central.  A central Sylow subgroup splits
off as an abelian direct factor, contradicting perfectness.
Thus the full-support hypothesis holds for nontrivial perfect $G$.
\end{proof}

\begin{remark}
Proposition~\ref{prop:pointed-perfect-equivariantization-connected-20261027} applies to an arbitrary braided action on the pointed
category, with its given associator. The orbit factor in
\eqref{explicit:80658146e8:eq0117} is computed from the centralizer
action on the fiber algebra. For example, the central extension of
$G=C_2\times C_2$ by $C_2$ with total group $Q_8$ realizes
$\Rep(Q_8)$ as an equivariantization of $\Rep(C_2)$. Conjugation fixes
both simples of $\Rep(C_2)$. For $g\ne1$ the preimage of
$\langle g\rangle$ is $C_4$. Its two lifts of $g$ give singleton
classes in $C_4$ and are interchanged by conjugation in $Q_8$.
Thus $|g^G|=1$, $[C_G(g):C_G(g)_\nu]=2$, and each nonidentity fiber
contains a single class of size two. The identity fiber contributes
the two classes of size one. For the strict trivial action all eight
classes have size one.

A prime $q\ne p$ occurring downstairs can disappear from every nonidentity
fiber of a $C_p$-equivariantization, as in
Example~\ref{ex:nonidentity-fiber-prime-persistence-fails-20261027}.
\end{remark}

\begin{corollary}[The component bound for equivariantizations of pointed categories]
\label{cor:pointed-equivariantization-components-20261102}
Let $\mathcal D$ be a pointed braided fusion category over $\kk$
with group of invertible objects $A=\operatorname{Inv}(\mathcal D)$,
and let a finite group $G$ act categorically on $\mathcal D$ by braided
tensor autoequivalences.
Then $\Dcl(\mathcal D^G)$ has at most two components. If $G=C_p$ has
prime order and acts nontrivially on $A$, the graph is the union of
the single vertex $p$ and the complete graph on
$\pi(|A:A^{C_p}|)$. In particular it has two components exactly when
$p\nmid|A:A^{C_p}|$. If $G=C_p$ acts trivially on $A$, the graph is empty.
\end{corollary}
\begin{proof}
For $G=C_p$, the identity-fiber dimensions are the orbit sizes of
the characters of $\kk[A]$ by the identity-fiber formula~\eqref{explicit:80658146e8:eq0014} of Lemma~\ref{lem:class-equivariantization-20261026}. All have size one
if the action on $A$ is trivial. Otherwise a character moves, since
characters separate the elements of $A$, and its orbit has size $p$.
Every $g\ne1$ generates $C_p$. By
\eqref{explicit:80658146e8:eq0126}, every class in its nonempty fiber
has dimension $|A:A^{C_p}|$. These are all the fibers, which gives the
asserted graph and the criterion for two components.

For solvable $G$, take a composition series. A prime cyclic step with
pointed base produces at most two components by the preceding paragraph.
At any step with non-pointed base
Proposition~\ref{prop:class-graph-prime-equivariantization-20261027}
bounds the number of components by that of the preceding graph. This
proves the bound by induction. For nonsolvable $G$, first equivariantize
by its nilpotent hypercenter $N$. The intermediate graph has at most
two components. If $\mathcal D^N$ is non-pointed, apply
Theorem~\ref{thm:class-graph-connected-equivariantization-20261027}
to the induced $G/N$-action. If $\mathcal D^N$ is pointed, the centerless
nonsolvable group $G/N$ has a connected class graph with vertex set
$\pi(|G/N|)$, and
Proposition~\ref{prop:pointed-perfect-equivariantization-connected-20261027}
gives connectedness.
\end{proof}

\subsection{The two-component bound}

\begin{theorem}[The general two-component bound]
\label{thm:general-two-components-20261102}\label{prob:general-two-components-20261027}
For every weakly integral braided fusion category $\C$ over $\kk$,
\begin{equation}\label{eq:general-class-component-bound-20261102}
 \#\pi_0\bigl(\Dcl(\C)\bigr)\leq 2.
\end{equation}
Equality is attained by $\Rep(D(S_3))$.
\end{theorem}
\begin{proof}
Write $Z_2(\C)\simeq\Rep(G,z)$, allowing $z=1$, and let
$H=G/\langle z\rangle$. The Tannakian subcategory $\Rep(H)$ is the
even part of $Z_2(\C)$. Its central de-equivariantization
$\mathcal B=\C_H$ is weakly integral and braided, and
$\C\simeq\mathcal B^H$.
Under the centralizer correspondence for this de-equivariantization,
$Z_2(\mathcal B)$ is the de-equivariantization of $\Rep(G,z)$ by
$\Rep(H)$, hence is $\Vec$ if $z=1$ and $\sVec$ otherwise
\cite[Corollary~4.27 and Proposition~4.30(iii)]{DGNO}.
Theorems~\ref{thm:class-graph-modular-two-20261027}
and~\ref{thm:class-graph-slightly-degenerate-20261027} therefore
bound $\#\pi_0(\Dcl(\mathcal B))$ by two.

If $\mathcal B$ is non-pointed, apply
Theorem~\ref{thm:class-graph-connected-equivariantization-20261027}
to its braided $H$-action. If $\mathcal B$ is pointed,
Corollary~\ref{cor:pointed-equivariantization-components-20261102}
gives the bound for $\mathcal B^H$ directly. Finally, the nontrivial class dimensions of
$\Rep(D(S_3))$ are $4$ and $9$, giving the two isolated vertices
$2$ and $3$.
\end{proof}

\section{Centralizer descent and solvability}\label{sec:unified-solvability}

\subsection{Symmetric class images and prime Tannakian centralizers}
\label{sec:prime-centralizer-splitting-20261027}

For a noncentral $\Rep(C_p)\subseteq\C$, the centralizer controls the
components of $\Dcl(\C)$. The proof uses the structure of connected
separable algebras in symmetric fusion categories.

\begin{lemma}[Prime-free connected algebras in a symmetric category]
\label{lem:symmetric-class-prime-free-algebra-pointed-20261027}
Let $\mathcal S$ be a symmetric fusion category over $\kk$ (hence integral),
and let $A$ be a connected separable algebra object in $\mathcal S$. Suppose that every prime
dividing $\FPdim(A)$ is absent from $\Dcl(\mathcal S)$. Then the object
$A$ is pointed. In particular, if $A$ generates $\mathcal S$ as a fusion
category, then $\mathcal S$ is pointed.
\end{lemma}
\begin{proof}
Write $\mathcal S=\Rep(H,z)$, with its possibly super braiding. The
underlying tensor category is $\Rep(H)$, so a connected separable algebra
has the form
\begin{equation}\label{explicit:80658146e8:eq0127}
 A=\operatorname{Ind}_{L}^{H}\End(V),\qquad
 n:=\FPdim(A)=[H:L](\dim V)^2,
\end{equation}
where $L\le H$ and $V$ is an irreducible projective representation of
$L$; see \cite[Section~3.4]{OstrikModules}. The prime divisors of a
projective irreducible degree belong to the prime divisors of the group
order. This follows, for example, by realizing the projective
representation on a finite central covering group whose order has the
same prime divisors as the original group.

For every prime dividing $n$, the corresponding Sylow subgroup of $H$
is central, since that prime divides no conjugacy-class size of $H$.
Their product $P$ is therefore an abelian central Hall subgroup. A
complement to $P$ gives
\begin{equation}\label{explicit:80658146e8:eq0128}
 H=P\times H_0,\qquad \gcd(n,|H_0|)=1.
\end{equation}
The index $[H:L]$ divides $n$. Its coprimality with $|H_0|$ implies
that the projection of $L$ onto $H_0$ is onto. Goursat's lemma and the
coprimality of $|P|$ and $|H_0|$ then give
\begin{equation}\label{explicit:80658146e8:eq0129}
 L=P_1\times H_0\qquad(P_1\le P).
\end{equation}
Projective cocycles on this coprime direct product split into cocycles
on its factors: the mixed term, a pairing between their abelianizations,
is trivial. Consequently
\begin{equation}\label{explicit:80658146e8:eq0130}
 V\simeq V_1\otimes V_0,
\end{equation}
with irreducible projective representations of $P_1$ and $H_0$,
respectively. Since $(\dim V_0)^2$ divides $n$, while the prime divisors
of $\dim V_0$ divide $|H_0|$, we have $\dim V_0=1$.

Thus $H_0$ acts trivially on $H/L$ and by scalars on $V_0$. It acts
trivially on the induced algebra $A$. The underlying representation of
$A$ therefore factors through the abelian group $P$, and all its simple
constituents are invertible. If these constituents generate
$\mathcal S$, the latter category is pointed as well.
\end{proof}

\begin{proposition}[The component exception for a prime centralizer]
\label{prop:prime-centralizer-class-components-20261027}
Let $\C$ be a weakly integral braided fusion category over $\kk$.
Let $p$ be a prime, and let $\mathcal E\simeq\Rep(C_p)$ be a Tannakian
fusion subcategory of $\C$ not contained in $Z_2(\C)$.
Let $\mathcal B=\mathcal E'_{\C}$ be its M\"uger centralizer in $\C$. Then $\C$ has a faithful $C_p$-grading
with neutral component $\mathcal B$. For each algebra character $\mu$
of $K(\C)$,
\begin{equation}\label{explicit:80658146e8:eq0131}
 \FPdim(Z_{\mu}^{\C})=s_{\mu}\,
 \FPdim(Z_{\mu|_{\mathcal B}}^{\mathcal B}),\qquad s_{\mu}\in\{1,p\}.
\end{equation}
Every vertex and edge of $\Dcl(\mathcal B)$ survives in
$\Dcl(\C)$. The only possible increase in the number of connected
components is a new isolated vertex $p$, absent from
$\Dcl(\mathcal B)$.

If this exceptional configuration occurs, every class dimension
divisible by $p$ equals $p$. Its class image has the form
\begin{equation}\label{explicit:80658146e8:eq0132}
 \FPdim(Z_{i})=p,\qquad
 Y_i=\bigoplus_{g\in K_i}g,\qquad
 K_i\le\operatorname{Inv}(Z_2(\mathcal B)),\quad |K_i|=p.
\end{equation}
For odd $p$, the pointed subcategory with simple objects $K_i$ is
Tannakian; for $p=2$ it may instead be $\sVec$.
\end{proposition}
\begin{proof}
Monodromy with the invertible objects of $\mathcal E$ defines a grading
by $\widehat{\operatorname{Inv}(\mathcal E)}\simeq C_p$. Noncentrality
makes this grading nontrivial, hence faithful because $p$ is prime; its
neutral component is exactly $\mathcal E'$. The grading formula~\eqref{explicit:80658146e8:eq0018} of Lemma~\ref{lem:class-dimensions-grading-20261027} gives~\eqref{explicit:80658146e8:eq0131}.

Every character downstairs extends by Proposition~\ref{prop:extent}.
For any such extension, \eqref{explicit:80658146e8:eq0131} gives the divisibility of class dimensions. Moreover, \eqref{explicit:80658146e8:eq0131} introduces
no new prime other than $p$. It can only add edges incident with $p$;
such edges cannot increase the number of components. An increase can
therefore occur only when $p$ is new and isolated.

In this case, a $p$-divisible class dimension has prime support
$\{p\}$. Its restricted class dimension is $p$-free, so~\eqref{explicit:80658146e8:eq0131}
forces $\FPdim(Z_{i})=p$ and $\FPdim(Z_{\mu_i|_{\mathcal B}}^{\mathcal B})=1$.
Since $Y_i\in\C_{\ad}\subseteq\mathcal B$,
Lemma~\ref{lem:class-restriction-collapse-20261027} gives
$Y_i\in Z_2(\mathcal B)$. Apply
Lemma~\ref{lem:symmetric-class-prime-free-algebra-pointed-20261027}
to its connected separable algebra structure: the only prime in
$\FPdim(Y_i)$ is absent from the class graph of
$\langle Y_i\rangle\subseteq\mathcal B$. Thus $Y_i$ is pointed.
The internal-end identity~\eqref{eq:offdiagonal} of Proposition~\ref{prop:diagonal} and Schur's lemma identify its constituents
with the subgroup
\begin{equation}\label{explicit:80658146e8:eq0133}
 K_i=\{g\in\operatorname{Inv}(\C):g^+Z_i\simeq Z_i\},
\end{equation}
each with multiplicity one. Its order is consequently $p$. The last
assertion follows from the classification of symmetric pointed
categories. Notice that $K_i$ need not be the original copy of $C_p$
defining $\mathcal E$.
\end{proof}

\begin{lemma}[Relative autoequivalences of a prime extension]
\label{lem:relative-prime-graded-autoequivalences-20261027}
Let $p$ be a prime, and let
\begin{equation}\label{explicit:80658146e8:eq0134}
 \mathcal R=\bigoplus_{g\in C_p}\mathcal R_g
\end{equation}
be a faithfully $C_p$-graded fusion category over $\kk$.
Suppose its neutral
component $\mathcal A=\mathcal R_e$ has $p$-power Frobenius--Perron
dimension. In particular $\mathcal R$ is weakly integral, since
$\FPdim(\mathcal R)=p\FPdim(\mathcal A)$. The group of isomorphism classes of pairs $(T,\eta)$, where
$T$ is a tensor autoequivalence of $\mathcal R$ preserving each $\mathcal R_g$, and
$\eta:T|_{\mathcal A}\xrightarrow{\sim}\operatorname{Id}_{\mathcal A}$
is a monoidal natural isomorphism, is a finite $p$-group. Isomorphisms
of pairs are required to respect the specified neutral trivializations.

Consequently, let $H$ be a finite group of order prime to $p$ acting
categorically on $\mathcal R$ by tensor autoequivalences preserving each
$\mathcal R_g$. Every coherent trivialization of its restricted action
on $\mathcal A$ extends to a coherent trivialization on $\mathcal R$.
\end{lemma}
\begin{proof}
Every $\mathcal R_g$ is an invertible $\mathcal A$-bimodule category,
and multiplication induces bimodule equivalences
\begin{equation}\label{explicit:80658146e8:eq0135}
 M_{g,h}:\mathcal R_g\boxtimes_{\mathcal A}\mathcal R_h
 \xrightarrow{\sim}\mathcal R_{gh};
\end{equation}
see \cite[Theorem~6.1]{ENO-homotopy}. The neutral trivialization
$\eta$ turns each $T|_{\mathcal R_g}$ into an
$\mathcal A$-bimodule autoequivalence. Restriction therefore defines a
group homomorphism
\begin{equation}\label{explicit:80658146e8:eq0136}
 \operatorname{Aut}_{\otimes}^{C_p}(\mathcal R;\mathcal A)
 \longrightarrow
 \prod_{g\in C_p}
 \operatorname{Aut}_{\mathcal A-\mathcal A}(\mathcal R_g).
\end{equation}
Here and below the automorphism groups consist of isomorphism classes.
The endomorphism category of each invertible bimodule is equivalent to
$Z(\mathcal A)$, so every group on the right is isomorphic to
$\operatorname{Inv}Z(\mathcal A)$
\cite[Corollary~6.2 and Section~8.2]{ENO-homotopy}. This is a finite
$p$-group, since $\FPdim Z(\mathcal A)=(\FPdim\mathcal A)^2$ is a
power of $p$.

We check injectivity. For a pair in the kernel, choose bimodule natural
isomorphisms
$\alpha_g:T|_{\mathcal R_g}\xrightarrow{\sim}
\operatorname{Id}_{\mathcal R_g}$, with $\alpha_e=\eta$.
Transporting the tensor structure of $T$ along these isomorphisms
leaves bimodule natural automorphisms of the equivalences $M_{g,h}$.
Such automorphisms are scalars: the endomorphism space of an invertible
bimodule functor is the endomorphism space of the unit in its fusion
category of bimodule endofunctors, hence is one-dimensional.
Write these scalars as $\lambda(g,h)\in\kk^\times$.
The unit constraints make them normalized, and the pentagon is exactly
\begin{equation}\label{explicit:80658146e8:eq0137}
 \lambda(g,h)\lambda(gh,k)
 =\lambda(h,k)\lambda(g,hk).
\end{equation}
Thus $\lambda$ is a normalized group $2$-cocycle. Since
$H^2(C_p,\kk^\times)=0$ by~\eqref{eq:cyclic-even-cohomology}, rescaling the $\alpha_g$ makes all
$\lambda(g,h)$ equal to one. The pair is then isomorphic to the
identity pair. This proves the claimed injection and the first assertion.

For the second assertion, choose a coherent trivialization of the
$H$-action on $\mathcal A$. It makes the action a homomorphism from
$H$ to the preceding group of relative isomorphism classes. The
homomorphism is trivial by coprimality. Choose relative monoidal
isomorphisms between all acting functors and the identity. Their
failure to respect multiplication in $H$ is a $2$-cocycle with values in
\begin{equation}\label{explicit:80658146e8:eq0138}
 \{\gamma\in\operatorname{Aut}_{\otimes}(\operatorname{Id}_{\mathcal R})
       :\gamma|_{\mathcal A}=1\}
 \simeq\widehat{C_p}.
\end{equation}
For this identification, a natural monoidal automorphism trivial on
$\mathcal A$ is scalar on each indecomposable homogeneous
$\mathcal A$-module category; its scalars multiply according to the
grading. Since $|\widehat{C_p}|=p$ is coprime to $|H|$, we have
$H^2(H,\widehat{C_p})=0$ by
\cite[Chap.~XII, \S2, Cor.~2.7, p.~237]{CE}. Altering the chosen isomorphisms by a $1$-cochain makes the
trivialization coherent, without changing it on $\mathcal A$.
\end{proof}

\begin{lemma}[Coprime splitting across a prime Tannakian centralizer]
\label{lem:prime-centralizer-coprime-splitting-20261027}
Let $\C$ be a weakly integral braided fusion category over $\kk$.
Let $p$ be a prime, and let $\mathcal E\simeq\Rep(C_p)$ be a Tannakian
fusion subcategory of $\C$ not contained in $Z_2(\C)$. Suppose
$\mathcal P,\mathcal B\subseteq\mathcal E'_{\C}$ are full fusion
subcategories and their tensor product induces the braided decomposition
\begin{equation}\label{explicit:80658146e8:eq0139}
 \mathcal E'_{\C}=\mathcal P\boxtimes\mathcal B,
 \qquad \mathcal P\text{ pointed},\qquad
 \FPdim\mathcal P=p^a,\qquad p\nmid\FPdim\mathcal B.
\end{equation}
Then there is a braided $p$-category $\mathcal R_p\subseteq\C$ such
that
\begin{equation}\label{explicit:80658146e8:eq0140}
 \C\simeq\mathcal R_p\boxtimes\mathcal B,
 \qquad \FPdim\mathcal R_p=p^{a+1}.
\end{equation}
The equivalence extends the specified inclusion of $\mathcal B$.
\end{lemma}
\begin{proof}
Put $\mathcal S=Z_2(\mathcal B)$. We first show that $\mathcal S$ is
transparent in $\C$. It centralizes both factors of $\mathcal E'$, so
the double-centralizer identity~\eqref{eq:muger-centralizer-identities} gives
\begin{equation}\label{explicit:80658146e8:eq0141}
 \mathcal S\subseteq(\mathcal E')'
 =\mathcal E\vee Z_2(\C).
\end{equation}
Noncentrality and primality give $\mathcal E\cap Z_2(\C)=\Vec$.
As these two subcategories centralize each other, their join is the
braided product
\begin{equation}\label{explicit:80658146e8:eq0142}
 \mathcal E\vee Z_2(\C)\simeq
 \mathcal E\boxtimes Z_2(\C).
\end{equation}
The first factor grades this product by $C_p$. If its restriction to
$\mathcal S$ were nontrivial, it would be a faithful $C_p$-grading,
forcing $p\mid\FPdim\mathcal S$. This is impossible because
$\FPdim\mathcal S\mid\FPdim\mathcal B$. Hence
\begin{equation}\label{explicit:80658146e8:eq0143}
 \mathcal S\subseteq Z_2(\C),\qquad
 \mathcal B\cap Z_2(\C)=\mathcal S.
\end{equation}

Let $\mathcal Q=\mathcal B'_{\C}$. The centralizer dimension formula~\eqref{eq:muger-centralizer-identities}
and the faithful monodromy grading give
\begin{equation}\label{explicit:80658146e8:eq0144}
 \FPdim\mathcal Q=p^{a+1}\FPdim\mathcal S,
 \qquad \mathcal Q\cap\mathcal B=\mathcal S,
 \qquad \mathcal Q\cap\mathcal E'=\mathcal P\boxtimes\mathcal S.
\end{equation}
In particular, the grading of $\C$ restricts to a faithful $C_p$-grading
of $\mathcal Q$ with neutral component
$\mathcal P\boxtimes\mathcal S$. The join dimension formula~\eqref{eq:join-fp-dimension} also gives
\begin{equation}\label{explicit:80658146e8:eq0145}
 \mathcal B\vee\mathcal Q=\C.
\end{equation}

Suppose first that $\mathcal S=\Rep(H)$ is Tannakian. Its central
de-equivariantization
\begin{equation}\label{explicit:80658146e8:eq0146}
 \mathcal R=\mathcal Q_H
\end{equation}
is braided and faithfully $C_p$-graded, with neutral component
$\mathcal P$ and dimension $p^{a+1}$. The canonical $H$-action
preserves the grading and is coherently trivial on $\mathcal P$:
on the neutral component it is the de-equivariantization of the
specified product $\mathcal P\boxtimes\Rep(H)$ over
$\Rep(H)$. Since $p\nmid|H|$, Lemma~\ref{lem:relative-prime-graded-autoequivalences-20261027}
trivializes the entire action relative to $\mathcal P$. Therefore,
as a braided category over the specified symmetric factor,
\begin{equation}\label{explicit:80658146e8:eq0147}
 \mathcal Q\simeq\mathcal R\boxtimes\mathcal S.
\end{equation}

Suppose next that $\mathcal S$ is super-Tannakian but not Tannakian.
Then $p$ is odd. Write
$\mathcal S=\Rep(H,z)$ with $z\ne1$, and let
$\mathcal S_0=\Rep(L)$ be its even part, where
$L=H/\langle z\rangle$. De-equivariantize by this central even part:
\begin{equation}\label{explicit:80658146e8:eq0148}
 \mathcal R=\mathcal Q_L,
 \qquad \mathcal R_e\simeq\mathcal P\boxtimes\sVec,
 \qquad \FPdim\mathcal R=2p^{a+1}.
\end{equation}
The neutral component is pointed, so
$\mathcal R_{\ad}\subseteq\mathcal R_e$ is pointed and
$\mathcal R$ is nilpotent. Its braided Sylow decomposition gives
\begin{equation}\label{explicit:80658146e8:eq0149}
 \mathcal R\simeq\mathcal R_p\boxtimes\sVec;
\end{equation}
the second factor is the specified central fermionic subcategory,
which already has the full $2$-dimension. We use the uniqueness of the
Sylow factors in \cite[Theorem~1.1]{DGNO-nilpotent}.

The $L$-action preserves these factors. Its restriction to
$\mathcal R_p$ preserves the $C_p$-grading, whose neutral component
is $\mathcal P$, and is coherently trivial on $\mathcal P$.
Lemma~\ref{lem:relative-prime-graded-autoequivalences-20261027}
therefore trivializes this restricted action. The action on the
$\sVec$ factor, including its coherence, is precisely the action
whose equivariantization reconstructs $\mathcal S$ from
$\mathcal S_0$. The braided product functor is equivariant because
the action preserves both factors. It follows that
\begin{equation}\label{explicit:80658146e8:eq0150}
 \mathcal Q\simeq\mathcal R^L
 \simeq\mathcal R_p\boxtimes(\sVec)^L
 \simeq\mathcal R_p\boxtimes\mathcal S.
\end{equation}
In both cases we have obtained a $p$-category $\mathcal R_p$ inside
$\mathcal Q$. It centralizes $\mathcal B$, and its intersection with
$\mathcal B$ is trivial by coprimality. Since
$\mathcal S\subseteq\mathcal B$ and
$\mathcal B\vee\mathcal Q=\C$, these two categories generate
$\C$. Their braided product is therefore the asserted equivalence.
\end{proof}

\begin{corollary}[Connectedness across a prime Tannakian centralizer]
\label{cor:prime-centralizer-connectedness-20261027}
Let $\C$ be a weakly integral braided fusion category over $\kk$.
Let $p$ be a prime, let $\mathcal E\simeq\Rep(C_p)$ be a Tannakian
fusion subcategory of $\C$ not contained in $Z_2(\C)$, and let
$\mathcal B=\mathcal E'_{\C}$ be its M\"uger centralizer. If $\mathcal B$ is non-pointed, then
\begin{equation}\label{explicit:80658146e8:eq0151}
 \#\pi_0\Dcl(\C)\le\#\pi_0\Dcl(\mathcal B).
\end{equation}
In particular a nonempty connected centralizer graph gives a connected
graph upstairs. If $\mathcal B$ is pointed, then $\C$ is nilpotent
and its class graph is empty or complete.
\end{corollary}
\begin{proof}
By Proposition~\ref{prop:prime-centralizer-class-components-20261027},
an increase could only be a new isolated $p$. In that case $p$ is absent
from the graph of $\mathcal B$, so Theorem~\ref{thm:main} gives
\begin{equation}\label{explicit:80658146e8:eq0152}
 \mathcal B=\mathcal P_p\boxtimes\mathcal B_0,
 \qquad p\nmid\FPdim\mathcal B_0,
\end{equation}
with $\mathcal P_p$ pointed. The coprime splitting of Lemma~\ref{lem:prime-centralizer-coprime-splitting-20261027} gives
$\C=\mathcal R_p\boxtimes\mathcal B_0$. Since $\mathcal B_0$
is non-pointed, Proposition~\ref{prop:class-graph-product-20261026} shows that a non-pointed $\mathcal R_p$ would make the product
graph connected, while a pointed $\mathcal R_p$ would make $p$
absent upstairs. Neither possibility gives a new isolated $p$.
If $\mathcal B$ is pointed, the monodromy grading gives
$\C_{\ad}\subseteq\mathcal B$, so $\C$ is nilpotent. Apply
Proposition~\ref{prop:class-graph-nilpotent-20261026}.
\end{proof}

\subsection{Disconnected class graphs force solvability}
\label{sec:general-class-graph-solvability-20261027}

In this section a minimal counterexample means a nonsolvable weakly
integral braided fusion category of least Frobenius--Perron dimension
among those with disconnected class graph.

\begin{theorem}[Solvable M\"uger center]
\label{thm:disconnected-class-graph-solvable-center-20261027}
Let $\C$ be a weakly integral braided fusion category over $\kk$.
Suppose that $Z_2(\C)$ is equivalent as a symmetric fusion category to
$\Rep(G,z)$, where $G$ is a finite solvable group and $z\in Z(G)$
satisfies $z^2=1$. If
$\Dcl(\C)$ is disconnected, then $\C$ is solvable.
\end{theorem}
\begin{proof}
Let $\Rep(H)$ be the maximal even Tannakian subcategory of
$\Rep(G,z)$; thus $H=G$ when $z=1$, and
$H=G/\langle z\rangle$ otherwise. It is central in $\C$, so
\begin{equation}\label{explicit:80658146e8:eq0153}
 \mathcal D=\C_H\quad\text{is braided},\qquad
 \C\simeq\mathcal D^H.
\end{equation}
Central de-equivariantization gives
\begin{equation}\label{explicit:80658146e8:eq0154}
 Z_2(\mathcal D)\simeq
 \begin{cases}
 \Vec,&z=1,\\
 \sVec,&z\ne1.
 \end{cases}
\end{equation}
If $\C$ were nonsolvable, then $\mathcal D$ would be nonsolvable,
because $H$ is solvable and solvability is preserved by such
equivariantizations \cite[Proposition~4.5]{ENO}. The modular and
slightly degenerate solvability results
(Theorem~\ref{thm:disconnected-class-graph-modular-solvable-20261027}
and Corollary~\ref{cor:disconnected-class-graph-minimal-center-solvable-20261027})
then imply that $\Dcl(\mathcal D)$ is connected and nonempty.
Theorem~\ref{thm:class-graph-connected-equivariantization-20261027}
makes $\Dcl(\C)$ connected as well, a contradiction.
\end{proof}

\begin{proposition}[A minimal counterexample is weakly group-theoretical]
\label{prop:class-graph-minimal-counterexample-wgt-20261027}
Let $\C$ be a nonsolvable weakly integral braided fusion category
over $\kk$ whose class graph has at least two nonempty connected
components. Suppose $\FPdim(\C)$ is minimal among categories with
these properties. Then $\C$ is weakly group-theoretical. More
precisely, $Z_2(\C)$ contains a Tannakian subcategory equivalent
to $\Rep(H)$ for a finite centerless nonsolvable group $H$. De-equivariantizing $\C$ by this subcategory gives a weakly integral
braided fusion category $\mathcal D=\C_H$ over $\kk$ with a
categorical $H$-action by braided tensor autoequivalences, and
\begin{equation}\label{explicit:80658146e8:eq0155}
 \C\simeq\mathcal D^H,
 \qquad \mathcal D=\C_H\text{ is solvable}.
\end{equation}
\end{proposition}
\begin{proof}
Theorem~\ref{thm:disconnected-class-graph-solvable-center-20261027} makes the underlying group of $Z_2(\C)$
nonsolvable. The group of its even Tannakian part is also nonsolvable,
since the omitted subgroup has order at most two. Quotient this group
by its hypercenter. The resulting group $H$ is centerless and
nonsolvable, and its representation category embeds centrally in
$\C$. The hypercenter is nilpotent, so quotienting by it cannot make
a nonsolvable group solvable.

Let $\mathcal D=\C_H$. It is weakly integral and braided, has
strictly smaller dimension, and $\C\simeq\mathcal D^H$.
If $\mathcal D$ were nonsolvable, minimality would make its class
graph connected. The graph is nonempty because an empty class graph
forces pointedness by Corollary~\ref{cor:im-class-support-20261024}. Theorem~\ref{thm:class-graph-connected-equivariantization-20261027}
would then make the class graph of $\C$ connected, a contradiction.
Thus $\mathcal D$ is solvable. Since weakly group-theoretical
categories are closed under equivariantizations by arbitrary finite
groups \cite[Proposition~4.1]{ENO}, $\C$ is weakly group-theoretical.
\end{proof}

\begin{lemma}[Tannakian subcategories of a minimal counterexample]
\label{lem:minimal-class-counterexample-central-tannakians-20261027}
Let $\C$ be a nonsolvable weakly integral braided fusion category
over $\kk$ whose class graph has at least two nonempty connected
components. Suppose $\FPdim(\C)$ is minimal among categories with
these properties. Then $\C$ contains no Tannakian subcategory equivalent
to $\Rep(C_p)$ for any prime $p$. Every Tannakian subcategory of
$\C$ is contained in $Z_2(\C)$. In particular the group of every nontrivial
Tannakian subcategory is perfect and nonsolvable.
\end{lemma}
\begin{proof}
Suppose first that $\mathcal E=\Rep(C_p)$ is central. The braided
category $\mathcal D=\C_{C_p}$ is smaller and must be nonsolvable:
otherwise its $C_p$-equivariantization $\C$ would be solvable.
Minimality gives a nonempty connected class graph for $\mathcal D$.
Theorem~\ref{thm:class-graph-connected-equivariantization-20261027}
then contradicts disconnectedness of $\Dcl(\C)$.

Suppose instead that $\mathcal E=\Rep(C_p)$ is noncentral. Its
centralizer $\mathcal B=\mathcal E'$ is proper and nonsolvable.
Indeed, if
$\mathcal B$ were solvable, its de-equivariantization by
$\Rep(C_p)$ would be solvable, and $\C_{C_p}$ is an extension of
this neutral category by a subgroup of $C_p$. Re-equivariantizing would
make $\C$ solvable. This is also a special case of
\cite[Proposition~5.2]{EscanuelaNatale}.
Minimality makes $\Dcl(\mathcal B)$ connected and nonempty.
Corollary~\ref{cor:prime-centralizer-connectedness-20261027}
again makes $\Dcl(\C)$ connected. We have excluded all
$\Rep(C_p)$.

Every nontrivial Tannakian subcategory $\mathcal E=\Rep(G)$
therefore has $G$ perfect: a nontrivial abelianization has a quotient
$C_p$, producing a forbidden subcategory $\Rep(C_p)$. A nontrivial
finite perfect group is nonsolvable. Suppose that such an $\mathcal E$
is noncentral. Its centralizer $\mathcal L=\mathcal E'$ is proper
and nonsolvable, because it contains the nonsolvable category
$\Rep(G)$. Thus $\Dcl(\mathcal L)$ is nonempty and connected,
and lies in a single component $K$ of $\Dcl(\C)$.

The ordinary class graph of $G$ has $2$ as a vertex. Indeed, a
nontrivial perfect group has even order by the odd-order theorem, and
it cannot have a central Sylow subgroup: such a subgroup would split
off as a nontrivial abelian direct factor. Since
$\mathcal E\subseteq\mathcal L$, class divisibility~\eqref{eq:braided-divisibility} of Theorem~\ref{thm:braided-divisibility} gives $2\in K$.

Choose a class object $Z_i$ whose dimension $\FPdim(Z_{i})>1$ belongs to a
component outside $K$. Individual class divisibility~\eqref{eq:braided-divisibility} gives
\begin{equation}\label{explicit:80658146e8:eq0156}
 \FPdim(Z_{\mu_i|_{\mathcal E}}^{\mathcal E})=1,
 \qquad \FPdim(Z_{\mu_i|_{\mathcal L}}^{\mathcal L})=1.
\end{equation}
Apply the modular quotient formula~\eqref{eq:modular-quotient} of Theorem~\ref{thm:modular-quotient} in
$Z(\C)$ to $\mathcal E^+$ and the simple object $Z_i$. The first
equality forces all of $Z_iZ_i^*$ to centralize $\mathcal E^+$.
For every $X\prec Y_i$, the internal-end identity~\eqref{eq:offdiagonal} places $X^+$ in
$Z_iZ_i^*$. Hence $Y_i$ centralizes $\mathcal E$, that is,
$Y_i\in\mathcal L$.
Now Lemma~\ref{lem:class-restriction-collapse-20261027}, using the
second equality, gives
\begin{equation}\label{explicit:80658146e8:eq0157}
 Y_i\in Z_2(\mathcal L).
\end{equation}

Let $\mathcal S=\langle Y_i\rangle$. It is symmetric, and its class
graph is a subgraph of $\Dcl(\mathcal L)$, whereas every prime
dividing $\FPdim(Z_{i})$ lies outside $K$.
By Proposition~\ref{prop:diagonal}, the internal end identified in~\eqref{eq:offdiagonal} carries a
connected separable algebra structure. Lemma~\ref{lem:symmetric-class-prime-free-algebra-pointed-20261027}
therefore makes $Y_i$ pointed. Its invertible constituents form the
stabilizer subgroup
\begin{equation}\label{explicit:80658146e8:eq0158}
 T_i=\{a\in\operatorname{Inv}(\C):a^+Z_i\simeq Z_i\},
\end{equation}
each with multiplicity one. Consequently
$\FPdim\mathcal S=|T_i|=\FPdim(Z_{i})>1$. Since $2\in K$, this order is
odd. A symmetric pointed category of odd dimension is Tannakian and
contains a subcategory $\Rep(C_p)$ for a prime divisor $p$ of its
dimension. This contradicts the first part of the proof. Every
Tannakian subcategory is therefore central.
\end{proof}

\begin{theorem}[Disconnected conjugacy-class prime graph implies solvability]
\label{thm:disconnected-class-graph-braided-solvable-20261027}
Let $\C$ be a weakly integral braided fusion category over $\kk$. If
$\Dcl(\C)$ has at least two nonempty connected components, then
$\C$ is solvable.
Equivalently, the class prime graph of every nonsolvable weakly
integral braided fusion category is nonempty and connected.
\end{theorem}
\begin{proof}
Suppose that the assertion fails, and choose a counterexample $\C$
of least Frobenius--Perron dimension. By
Proposition~\ref{prop:class-graph-minimal-counterexample-wgt-20261027},
$\C$ is weakly group-theoretical. By
Lemma~\ref{lem:minimal-class-counterexample-central-tannakians-20261027},
it contains no $\Rep(C_p)$, and every Tannakian subcategory is
central.

Choose a maximal Tannakian subcategory
$\mathcal E=\Rep(G)$. It is nontrivial by
Proposition~\ref{prop:class-graph-minimal-counterexample-wgt-20261027}.
The group $G$ is perfect and nonsolvable. Since $\mathcal E$ is
central, the whole de-equivariantization
\begin{equation}\label{explicit:80658146e8:eq0159}
 \mathcal D=\C_G
\end{equation}
is braided and coincides with the core of $\C$ relative to
$\mathcal E$. The core theorem for weakly group-theoretical braided
categories
\cite[Theorem~1.1]{Natale-core} gives
\begin{equation}\label{explicit:80658146e8:eq0160}
 \mathcal D\simeq\mathcal P
 \quad\text{or}\quad
 \mathcal D\simeq\mathcal I\boxtimes\mathcal P,
\end{equation}
where $\mathcal P$ is pointed and weakly anisotropic, and
$\mathcal I$ is an Ising braided category. In the second case the class graph of $\mathcal D$ is the single
vertex $\{2\}$: the Ising class dimensions are $1,1,2$, and a pointed
factor preserves the graph by Proposition~\ref{prop:class-graph-product-20261026}. Theorem~\ref{thm:class-graph-connected-equivariantization-20261027}
then makes the graph of $\C\simeq\mathcal D^G$ connected.
In the first case, Proposition~\ref{prop:pointed-perfect-equivariantization-connected-20261027}
gives the same conclusion because $G$ is nontrivial and perfect.
Both cases contradict the defining property of $\C$.

Finally, Corollary~\ref{cor:im-class-support-20261024} makes a
category with empty class graph pointed, hence solvable. This proves the equivalent formulation as well.
\end{proof}

\begin{remark}
The group of a maximal Tannakian subcategory in the proof is perfect.
Proposition~\ref{prop:pointed-perfect-equivariantization-connected-20261027}
therefore applies in the pointed-core case even when this group has a
nontrivial center.
\end{remark}

\section{Vanishing across the two components}\label{sec:component-bound-vanishing-20261102}

Theorems~\ref{thm:general-two-components-20261102}
and~\ref{thm:disconnected-class-graph-braided-solvable-20261027} show that a disconnected class graph has exactly two
nonempty components and its category is solvable. In particular,
if the class graph of $\C_{\mathrm{int}}$ has two components supported
on odd primes, Corollary~\ref{cor:class-graph-odd-primes-20261027}
and~\eqref{eq:general-class-component-bound-20261102} force any new
vertex $2$ to be adjacent to an odd prime.

\begin{problem}[Cross-component vanishing]\label{prob:cross-component-vanishing-20261027}
Let $\C$ be a weakly integral braided fusion category over $\kk$.
If $\FPdim(Z_i),\FPdim(Z_j)>1$ belong to different components of
$\Dcl(\C)$, must $\mu_i(Y_j)=0$?
\end{problem}

This is open for general weakly integral braided fusion categories.
By Lemma~\ref{lem:coprime-class-vanishing-dichotomy-20261027}, it
amounts to excluding centralization of the two corresponding class
objects in $Z(\C)$. For a category satisfying this vanishing for
all such pairs, the cross-component regular-product formula~\eqref{eq:cross-regular-product} of Proposition~\ref{prop:class-graph-vanishing-criterion-20261027} gives
\begin{equation}\label{explicit:80658146e8:eq0162}
 [Y_i][Y_j]=\frac{\FPdim(Z_i)\FPdim(Z_j)}{\FPdim(\C_{\ad})}
 R_{\C_{\ad}}.
\end{equation}
This imposes a numerical condition on every class in either component
$\Pi$. For a positive integer $n$, write
$n_{\Pi}=\prod_{p\in\Pi}p^{v_p(n)}$. Fix a class $Z_i$ in $\Pi$ and
take $Z_j$ in the other component. The unit coefficient in
\eqref{explicit:80658146e8:eq0162} gives
$\FPdim(\C_{\ad})\mid\FPdim(Z_i)\FPdim(Z_j)$, while no prime in
$\Pi$ divides $\FPdim(Z_j)$. The formal codegree
$\FPdim(\C)/\FPdim(Z_i)$ is an integer. Taking the $\Pi$-parts
therefore gives
\begin{equation}\label{eq:component-class-valuation-20261102}
 \bigl(\FPdim(\C_{\ad})\bigr)_{\Pi}
 \ \mid\ \FPdim(Z_i)\ \mid\
 \bigl(\FPdim(\C_{\ad})\bigr)_{\Pi}\,|U(\C)|_{\Pi}.
\end{equation}
By Theorem~\ref{thm:main}, every prime of $\Pi$ divides
$\FPdim(\C_{\ad})$. Thus vanishing across components forces every
nontrivial class in $\Pi$ to contain all primes of $\Pi$ and forces
each component to be complete. A class whose dimension fails the
first divisibility in \eqref{eq:component-class-valuation-20261102}
would give a negative answer to
Problem~\ref{prob:cross-component-vanishing-20261027}.

\subsection{Survival of components under equivariantization}

Suppose $\Dcl(\mathcal B)$ has two components. Under a braided action
of a finite group $G$, Theorem~\ref{thm:class-graph-connected-equivariantization-20261027}
shows that each old vertex and edge survives equivariantization and each upstairs
component contains an old vertex.
Thus the graph of $\mathcal B^G$ is disconnected exactly when the two
old components remain separate. To detect possible paths between
them, one can use the all-fiber equivariantization formula~\eqref{explicit:80658146e8:eq0117} of Proposition~\ref{prop:equivariant-class-fiber-formula-20261027}.
For $g\in G$, $L=\langle g\rangle$, $Q=C_G(g)$, and a character $\nu$
in the $g$-fiber of $K(\mathcal B^L)$, the character $\mu$ indexed by
its $Q$-orbit has class dimension
\begin{equation}\label{explicit:80658146e8:eq0161}
 \FPdim(Z_{\mu}^{\mathcal B^G})
 =|g^G|\,[Q:Q_\nu]\,\FPdim(Z_{\nu}^{\mathcal B^L}).
\end{equation}
The following example distinguishes persistence in the whole graph
from persistence in a specified fiber.

\begin{example}[A prime need not persist in a nonidentity fiber]
\label{ex:nonidentity-fiber-prime-persistence-fails-20261027}
Let $K=C_7\rtimes C_6$ be the faithful affine group, and let
$N=C_7\rtimes C_3$ be its normal subgroup of index two. The split
extension gives
\begin{equation}\label{explicit:80658146e8:eq0163}
 \Rep(K)\simeq\Rep(N)^{C_2}.
\end{equation}
The nonidentity elements of the kernel $C_7$ form two $N$-classes of
size $3$. The outer involution interchanges them, and they become one
$K$-class of size $6$. Every element of $K\setminus N$ has nontrivial
image in $C_6$ and lies in a class of size $7$: the complement acts
fixed-point freely on the nonzero elements of $C_7$, and the
centralizer of such an element has order $6$.

Thus $3$ occurs in the downstairs class graph, but every class in the
nonidentity $\Rep(C_2)$-fiber upstairs has dimension $7$. In particular,
the implication
\begin{equation}\label{explicit:80658146e8:eq0164}
 q\ne p,\quad q\in\operatorname{Vert}\Dcl(\mathcal B)
 \quad\Longrightarrow\quad
 q\text{ occurs in a nonidentity fiber of }\mathcal B^{C_p}
\end{equation}
is false without further hypotheses. Here $p=2$ and $q=3$.
The orbit factor in the identity fiber supplies the actual connection,
through the class dimension $6$.
\end{example}

\subsection{Finite groups, untwisted doubles and the two components}

\begin{proposition}[Vanishing for finite-group representation categories]
\label{prop:symmetric-cross-component-vanishing-20261102}
Let $\C$ be a fusion category over $\kk$ tensor equivalent to
$\Rep(G)$ for a finite group $G$, so $\C$ is integral. The graph and class images
are formed using the underlying tensor category. Suppose the conjugacy-class dimensions
$\FPdim(Z_i),\FPdim(Z_j)>1$ belong to different components of
$\Dcl(\C)$. Then $\mu_i(Y_j)=0$, where $Y_j=F(Z_j)$.
In particular, cross-component vanishing holds for every symmetric
fusion category over $\kk$.
\end{proposition}
\begin{proof}
Realize $Z(\Rep(G))$ as the category of conjugation-equivariant
vector bundles on $G$. In this model $I(\one)$ has a trivial line
over each group element. Its simple summands are the bundles of
trivial lines over the conjugacy classes. Thus a class supported on
$y^G$ has image $\kk[y^G]$ with its permutation action. Every
character of $K(\Rep(G))$ is evaluation at a conjugacy class $x^G$,
and its formal codegree is $|C_G(x)|$ by column orthogonality~\eqref{eq:traceexpansion} for $\Rep(G)$.
If $\mu_i$ is this evaluation, its class dimension is $|x^G|$ and
\begin{equation}\label{eq:symmetric-class-image-fixed-points-20261102}
 \mu_i(\kk[y^G])=|y^G\cap C_G(x)|.
\end{equation}
These identifications use the underlying tensor category.

Put $Z=Z(G)$. By \cite[Theorem~2.5(ii)]{NPR}, $G/Z$ is a Frobenius
group whose kernel $N/Z$ and complement $H/Z$ have abelian inverse
images $N$ and $H$ in $G$. For $x\in N\setminus Z$, the Frobenius
action gives $C_{G/Z}(xZ)=N/Z$. Since $N$ is abelian, $C_G(x)=N$.
Every $x\notin N$ lies in a conjugate $H^u$ of $H$. The centralizer
of $xZ$ in $G/Z$ is contained in $H^u/Z$ by the Frobenius property,
and $H^u$ is abelian, so $C_G(x)=H^u$. The two nontrivial class sizes
are consequently $|H/Z|$ on $N\setminus Z$ and $|N/Z|$ on its
complement in $G\setminus Z$. They are coprime. Conjugates of
elements giving opposite sizes cannot commute since
$N\cap H^u=Z$. Equation~\eqref{eq:symmetric-class-image-fixed-points-20261102}
now gives the required vanishing. Finally a symmetric fusion
category over $\kk$ is super-Tannakian, and $\Rep(G,z)$ has underlying
tensor category $\Rep(G)$.
\end{proof}

\begin{proposition}[Pointed equivariantizations with odd-order monodromy]
\label{prop:class-graph-odd-pointed-equivariantization-20261027}
Let $\mathcal B$ be a pointed braided fusion category over $\kk$,
and let a finite group $G$ act categorically on $\mathcal B$ by braided
tensor autoequivalences. Suppose the monodromy bicharacter on
$\operatorname{Inv}(\mathcal B)$ determined by its braiding has odd order. Then $\mathcal B^G$ admits a symmetric braiding on its
underlying fusion category. In particular, its class graph is the conjugacy-class-size graph of a finite group.
If the graph is disconnected, the category satisfies cross-component
vanishing.
\end{proposition}
\begin{proof}
Let $b$ be the monodromy bicharacter on the group of simple objects, and
let $n$ be its odd order. Choose an integer $r$ with $2r\equiv1\pmod n$.
For homogeneous objects replace $c_{X,Y}$ by
$b(X,Y)^{-r}c_{X,Y}$. Multiplying a braiding by a bicharacter preserves
the hexagon identities. The new monodromy is $b^{1-2r}=1$, so the
new braiding is symmetric. Every braided autoequivalence in the given
action preserves $b$, and therefore preserves this new braiding, with
the same tensor and action coherence maps. The equivariantization is
thus symmetric. By the super-Tannakian description, its underlying
fusion category is tensor equivalent to $\Rep(K)$ for a finite group
$K$, with a possible super braiding. By~\eqref{explicit:80658146e8:eq0003}, class dimensions are determined
by the Frobenius--Perron dimension and the formal codegrees of the
Grothendieck ring, so they equal the conjugacy-class sizes of $K$.
Proposition~\ref{prop:symmetric-cross-component-vanishing-20261102} gives the vanishing assertion.
\end{proof}

\begin{proposition}[Vanishing for untwisted Drinfeld doubles]
\label{prop:untwisted-double-cross-component-vanishing-20261102}
Let $G$ be a finite group and $\C=\Rep(D(G))$. If $X,Y\in\Irr(\C)$
label nontrivial classes in different components of $\Dcl(\C)$, then
$S_{X,Y}=0$ and $\mu_X(Y\otimes Y^*)=0$.
\end{proposition}
\begin{proof}
The graph on class dimensions of $\C$ is the graph on its simple
dimensions, since the former are the squares of the latter by
Lemma~\ref{lem:class-square-minimal-center-20261026}\textup{(1)}.
By \cite[Theorem~8.5(iii)]{NPR}, this disconnected graph is the
class-size prime graph of $G$. With $Z=Z(G)$, the quotient $G/Z$
is Frobenius, and the inverse images $N,H$ of its kernel and a
complement are abelian. Put $m=|N/Z|$ and $n=|H/Z|$.
The centralizer computation in
Proposition~\ref{prop:symmetric-cross-component-vanishing-20261102}
gives class sizes $n$ on $N\setminus Z$ and $m$ on $G\setminus N$.
In both cases the centralizer is abelian, so every simple with that
flux has the indicated class size as its simple dimension. A simple
with central flux has dimension $\chi(1)$ for some
$\chi\in\Irr(G)$, which divides $[G:N]=n$ by It\^o's
abelian-normal-subgroup theorem \cite{Ito1951}. Thus one member of an
opposite-component pair has outer flux $y\notin N$, and the other has
either noncentral kernel flux or central flux with a nonlinear charge.

In the first case the two flux classes have no commuting
representatives. Their $S$-entry vanishes by the commuting-pair
formula \cite[(2.12)--(2.13)]{CGR}. In the second case the same
formula reduces vanishing to $\chi(y)=0$. Conjugate $H$ so that
$y\in H$. The endomorphism $a\mapsto yay^{-1}a^{-1}$ of the abelian
group $N$ has kernel $Z$ and induces an automorphism of $N/Z$ by the
Frobenius property. Its image $K$ therefore has order $m$ and
$N=KZ$. The subgroup $K$ is normalized by $N$ and by $H$, since $H$
is abelian and centralizes $y$. Now $G=KH$ and $H$ is abelian, so
$[G,G]\subseteq K$. The reverse inclusion follows from the definition
of $K$. It follows that $|G/[G,G]|=|H|=|C_G(y)|$. In the column orthogonality relation~\eqref{eq:traceexpansion} at $y$, the linear characters already contribute
$|C_G(y)|$, hence every nonlinear $\chi$ vanishes there. Finally,
the character evaluation in~\eqref{eq:modular-cross-component-class-value-20261102} gives
$\mu_X(Y\otimes Y^*)=|S_{X,Y}|^2/\FPdim(X)^2=0$.
\end{proof}

For twisted doubles the simple objects can carry projective
centralizer representations. It remains to determine which twists admit the same vanishing argument.

For group-theoretical modular categories, the identification with the
prime graph of a twisted double in the proof of \cite[Theorem~1.2(ii)]{NPR}
gives two complete components whenever the graph is disconnected
\cite[Theorem~8.5(iii)]{NPR}. In the connected case its diameter is at
most three \cite[Theorem~1.2(ii)]{NPR}. For arbitrary weakly integral
braided categories, must each component of a disconnected class graph be
complete? Can one characterize its two prime sets categorically, as the
kernel and complement of the Frobenius quotient do for finite groups?
Does a uniform diameter bound hold for connected class graphs?

A related comparison concerns the Tannakian category in
Proposition~\ref{prop:common-tannakian-localization-20261027}.
For group elements $x,y$ with coprime class sizes, the intersection of
their normal closures is abelian \cite[Theorem~B]{Camina}. For the
pure-class objects attached to $x$ and $y$, it would be of interest to
describe how this intersection is reflected in their common adjoint
constituents.

\appendix

\section{Relative vanishing and the fibers of the forgetful functor}
\label{sec:general-relative}

For a fusion category $\mathcal A$ with commutative Grothendieck ring,
write $\widehat K(\mathcal A)$ for the set of algebra characters of
$K(\mathcal A)$. The following estimate locates a zero on a prescribed
set of characters and bounds the total weight of such zeros.
Its proof uses the first-moment argument for zeros of $S$-characters
\cite[Theorem~4.3]{SerreZeros} and the Galois-orbit argument
of~\cite[Theorem~4.17]{BP}.

\subsection{A relative Burnside theorem}

\begin{theorem}[Relative vanishing]\label{thm:general-relative}
Let $\mathcal A$ be a fusion category over $\kk$ with commutative
Grothendieck ring. Let $B\in\mathcal A$ have a simple summand $Z$,
and let $\widehat K(\mathcal A)$ be the set of algebra characters
of its complexified Grothendieck ring. For each
$\rho\in\widehat K(\mathcal A)$, let $n_\rho$ be its formal codegree.
Suppose that $k_\rho=\rho(B)\in\mathbb Z_{\ge0}$ for every $\rho$ and
that $\rho(Z)\in\mathbb R_{\ge0}$ whenever $k_\rho>0$.
The subcategory $\mathcal H=\langle B\rangle$ is integral. Moreover,
\begin{equation}\label{eq:general-multiplicity-bound}
 \sum_{\substack{\rho(Z)=0\\k_\rho>0}}\frac{k_\rho}{n_\rho}
 \ \ge\ [B:\one]-[B:Z]+
 \frac{\FPdim(B)}{\FPdim(\mathcal H)}\bigl(\FPdim(Z)-1\bigr).
\end{equation}
Equality in \eqref{eq:general-multiplicity-bound} holds if and only if
every character $\psi\ne\FPdim$ of $K(\mathcal H)$ with $\psi(B)>0$
satisfies $\psi(Z)\in\{0,1\}$. In particular, if $[B:\one]=[B:Z]=1$, then
\begin{equation}\label{eq:general-zero-bound}
 \sum_{\substack{\rho(Z)=0\\k_\rho>0}}\frac{k_\rho}{n_\rho}
 \ \ge\ \frac{\FPdim(B)}{\FPdim(\mathcal H)}
                 \bigl(\FPdim(Z)-1\bigr).
\end{equation}
A character with $\rho(B)>0$ and $\rho(Z)=0$ therefore exists whenever
$\FPdim(Z)>1$ and $[B:\one]\ge[B:Z]$.
\end{theorem}

\begin{proof}
Since $\rho(B^*)=\overline{\rho(B)}=\rho(B)$ for every character,
$B\simeq B^*$. In particular $[B:Z^*]=[B:Z]$, and every simple
$X\in\mathcal H$ occurs in some $B^{\otimes m}$. The dimension character
gives $\FPdim(B)\in\mathbb Z$. For any Galois automorphism $\gamma$
of a splitting field of the fusion matrices,
$|\gamma(\FPdim(X))|\le\FPdim(X)$ by the spectral-radius bound.
Applying $\gamma$ to the simple decomposition of $\FPdim(B)^m$ and
taking real parts shows that each constituent attains equality in
$\Re\gamma(\FPdim(X))\le\FPdim(X)$. Each $\FPdim(X)$ is consequently
a rational algebraic integer.

Every character $\psi$ of $K(\mathcal H)$ extends to $K(\mathcal A)$
since the inclusion is unital and both algebras are commutative
semisimple. Thus the two hypotheses on character values hold also in
$\mathcal H$. Write $f_\psi$ for its formal codegree. The corresponding
class object of $\mathcal H$ has dimension
$\FPdim(\mathcal H)/f_\psi$ by
\cite[Theorem~2.13]{OstrikRankThree}, which applies to the canonical
spherical structure of the integral category $\mathcal H$.
This dimension is an integer since its forgetful image belongs to
$\mathcal H$. In particular $f_\psi$ is rational. Set
$w_\psi=\FPdim(\mathcal H)\psi(B)/f_\psi$.
The unit-coefficient trace~\eqref{eq:app-class-idempotents-20261102} gives
\begin{equation}\label{eq:general-moments}
 \sum_\psi w_\psi=\FPdim(\mathcal H)[B:\one],
 \qquad
 \sum_\psi w_\psi\psi(Z)=\FPdim(\mathcal H)[B:Z].
\end{equation}
The nonnegative integer weights $w_\psi$ are Galois-invariant since
$\psi(B)$ is integral and $f_{\gamma\psi}=\gamma(f_\psi)$.
The dimension character $\psi_0$ is fixed by Galois, with
$w_{\psi_0}=\FPdim(B)$. On every other Galois orbit where
$w_\psi>0$ and $\psi(Z)\ne0$, the values $\psi(Z)$ are positive
algebraic integers. Their product is a positive rational algebraic
integer, so their arithmetic mean is at least one.
In the difference of the two moments in~\eqref{eq:general-moments},
these orbits thus contribute a nonnegative amount to
$\sum_\psi w_\psi(\psi(Z)-1)$. The dimension character contributes
$\FPdim(B)(\FPdim(Z)-1)$, and a zero contributes $-w_\psi$.
Rearranging and dividing by $\FPdim(\mathcal H)$ gives the bound
in $\mathcal H$. Equality requires each positive-weight orbit outside
the dimension character on which $\psi(Z)\ne0$ to have arithmetic
mean one. Its geometric mean is at least one, so equality in the
arithmetic--geometric mean inequality forces every value on the orbit
to be one. Conversely, this condition makes every such orbit contribute
zero to the difference.

For transfer to $\mathcal A$, the primitive idempotent of $\psi$
is the sum of the idempotents of its extensions. Taking the
unit-coefficient trace~\eqref{eq:app-class-idempotents-20261102} gives
\begin{equation}\label{eq:general-restriction-trace}
 \sum_{\rho|_{K(\mathcal H)}=\psi}\frac1{n_\rho}=\frac1{f_\psi}.
\end{equation}
Both $B$ and $Z$ belong to $\mathcal H$, so their values are constant
on these fibers. Summation preserves both the bound and its equality
condition.
\end{proof}

\begin{example}[Equality and the multiplicity hypothesis]\label{rem:relative-S3}
Let $Z$ be the two-dimensional irreducible of $S_3$. On the identity,
a transposition and a three-cycle its values are $2,0,-1$. For
$B=\one\oplus Z$ the corresponding values are $3,1,0$. The only
supported zero is the transposition class, whose formal codegree is
two. Its contribution $1/2$ attains~\eqref{eq:general-zero-bound}.
For the regular representation $B$, the values are $6,0,0$ and
$[B:Z]=2$. The supported zero set is empty and both sides
of~\eqref{eq:general-multiplicity-bound} are zero.
\end{example}

\subsection{The canonical algebra and three character fibers}
\label{sec:relative}

Return to a weakly integral braided category $\C$. For its $i$-th
class object, define
\begin{equation}\label{eq:B}
 \E_i=\langle Z_i\rangle\subseteq\Z(\C),\qquad
 \D_i=\langle Y_i\rangle\subseteq\C,\qquad
 B_i=I(\one)\cap\E_i=\bigoplus_{Z_j\in\E_i}Z_j .
\end{equation}
The restriction $F_i\colon\E_i\to\D_i$ of forgetting is a dominant
central functor. Indeed, all simples of $\D_i$ occur in tensor words
in $Y_i$ and $Y_i^*$, and the half-braidings of $\E_i$ restrict to
$\D_i$. By adjunction $B_i$ is the value at $\one$ of its right adjoint.
It is a connected \'etale algebra, and
\cite[Lemmas~3.5 and~3.11]{DMNO} give
\begin{equation}\label{eq:image-dimension}
 \FPdim(\D_i)=\frac{\FPdim(\E_i)}{\FPdim(B_i)}.
\end{equation}
The canonical algebra is self-dual because the multiplicities of $X$
and $X^*$ in $I(\one)$ both equal $[F(X):\one]$ by adjunction and
duality. Thus $B_i$ is self-dual, and multiplicity-freeness gives
$[B_i:\one]=[B_i:Z_i]=1$. Since $Z_i$ is a summand of $B_i$,
$\langle B_i\rangle=\E_i$.

For $\rho\in\widehat K(\E_i)$ and $\sigma\in\widehat K(\D_i)$ write
$n_\rho$ and $m_\sigma$ for their respective formal codegrees, and
put $k_\rho=\rho(B_i)$. Let $\chi_W$ be the Verlinde character of
$W\in\Irr(\Z(\C))$. Define the fibers and their weights by
\begin{align}
 \Omega_\rho&=\{W\in\Irr(\Z(\C)): \chi_W|_{K(\E_i)}=\rho\},
 & r_\rho&=\sum_{j:\,Z_j\in\Omega_\rho}\FPdim(Z_j),
 \label{eq:relative-full-fiber-data}\\
 \Pi_\rho&=\{\sigma\in\widehat K(\D_i): \sigma\circ F_{i*}=\rho\},
 & s_\rho&=\sum_{\sigma\in\Pi_\rho}\FPdim(Z_\sigma^{\D_i}).
 \label{eq:relative-image-fiber-data}
\end{align}
Here $F_{i*}$ denotes the induced algebra map.

\begin{lemma}[Full and canonical-algebra fibers]\label{lem:fibres}
Let $\C$ be a weakly integral braided fusion category over $\kk$. For the forgetful functor
$F\colon\Z(\C)\to\C$ with right adjoint $I$, write
$I(\one)=\bigoplus_j Z_j$ and $Y_j=F(Z_j)$. Fix $i$ and set
$\E_i=\langle Z_i\rangle\subseteq\Z(\C)$,
$\D_i=\langle Y_i\rangle\subseteq\C$, and $B_i=I(\one)\cap\E_i$.
Write $F_i\colon\E_i\to\D_i$ for the restriction of $F$.
For an algebra character $\rho$ of $K(\E_i)$ let $n_\rho$
be its formal codegree and set $k_\rho=\rho(B_i)$.
Let $\Omega_\rho$ consist of those $W\in\Irr(\Z(\C))$ whose
Verlinde character for the canonical modular structure restricts to
$\rho$ on $K(\E_i)$.
Put $r_\rho=\sum_{j:\,Z_j\in\Omega_\rho}\FPdim(Z_j)$.
Then
\begin{align}
 \sum_{W\in\Omega_\rho}\FPdim(W)^2
    &=\frac{\FPdim(\C)^2}{n_\rho},\label{eq:fullfibre}\\
 r_\rho&=\frac{\FPdim(\C)k_\rho}{n_\rho},
 \qquad k_\rho\in\mathbb Z_{\ge0}.\label{eq:Afibre}
\end{align}
Let $\mu_j$ be the algebra character of $K(\C)$ corresponding to $Z_j$.
The supported characters are characterized by
\begin{equation}\label{eq:survival}
 k_\rho>0\quad\Longleftrightarrow\quad
 \rho=(\mu_j|_{K(\D_i)})\circ F_{i*}\quad\text{for some }j .
\end{equation}
For such $j$, $\rho(Z_i)=\mu_j(Y_i)\ge0$.
\end{lemma}
\begin{proof}
The primitive idempotent $e_\rho$ of $K(\E_i)$ maps in
$K(\Z(\C))$ to the sum of the idempotents corresponding to
$W\in\Omega_\rho$. Their unit traces are $\FPdim(W)^2/\FPdim(\C)^2$,
which gives~\eqref{eq:fullfibre}. On the other hand, for
$v\in K(\E_i)$ adjunction and self-duality of $B_i$ give
\begin{equation}\label{eq:relative-forgetful-trace}
 \tau_{\C}(F_{i*}(v))=\tau_{\D_i}(F_{i*}(v))
                      =\tau_{\E_i}(v[B_i]).
\end{equation}
Indeed on a simple $X$ both sides equal $[F(X):\one]=[B_i:X]$.
By~\eqref{eq:app-master-20261102}, $\mu_j\circ F_*=\chi_{Z_j}$.
Hence in $K(\C)$ the image of $e_\rho$ is the sum of the
idempotents corresponding to those $j$ for which $Z_j\in\Omega_\rho$.
Its trace is $r_\rho/\FPdim(\C)$, while the right side
of~\eqref{eq:relative-forgetful-trace} gives $k_\rho/n_\rho$.
Since $n_\rho$ is rational and $k_\rho$ is an algebraic integer,
\eqref{eq:Afibre} gives $k_\rho\in\mathbb Z_{\ge0}$.
The last assertion follows from Proposition~\ref{prop:positivity}.
\end{proof}

\begin{proposition}[Three compatible fiber identities]\label{prop:image-fibres}
Let $\C$ be a weakly integral braided fusion category over $\kk$. For the forgetful functor
$F\colon\Z(\C)\to\C$ with right adjoint $I$, write
$I(\one)=\bigoplus_j Z_j$ and $Y_j=F(Z_j)$. Fix $i$ and set
$\E_i=\langle Z_i\rangle\subseteq\Z(\C)$,
$\D_i=\langle Y_i\rangle\subseteq\C$, and $B_i=I(\one)\cap\E_i$.
Write $F_i\colon\E_i\to\D_i$ for the restriction of $F$.
For an algebra character $\rho$ of $K(\E_i)$ let $n_\rho$
be its formal codegree and set $k_\rho=\rho(B_i)$.
Let $\Pi_\rho$ consist of the algebra characters $\sigma$ of
$K(\D_i)$ whose pullback along $F_{i*}$ equals $\rho$.
For these $\sigma$, let $m_\sigma$ denote the formal codegree in
$\D_i$ and put $s_\rho=\sum_{\sigma\in\Pi_\rho}\FPdim(Z_\sigma^{\D_i})$.
Write $r_\rho=\sum_j\FPdim(Z_j)$, summing over precisely those $j$
for which the central Verlinde character $\chi_{Z_j}$ restricts to
$\rho$ on $K(\E_i)$. Then
\begin{align}
 \sum_{\sigma\in\Pi_\rho}\frac1{m_\sigma}
      &=\frac{k_\rho}{n_\rho},\label{eq:image-codegrees}\\
 \FPdim(B_i)\,s_\rho&=\FPdim(Z_\rho^{\E_i})k_\rho,
      \label{eq:image-classes}\\
 r_\rho&=\frac{\FPdim(\C)}{\FPdim(\D_i)}s_\rho.
      \label{eq:image-to-C}
\end{align}
In particular $\Pi_\rho\ne\varnothing$ precisely when $k_\rho>0$, and
\begin{equation}\label{eq:image-divisibilities}
 \FPdim(B_i)\mid\FPdim(Z_\rho^{\E_i})k_\rho,\qquad
 n_\rho\mid\FPdim(\D_i)k_\rho,\qquad
 \frac{\FPdim(\C)}{\FPdim(\D_i)}\mid r_\rho.
\end{equation}
\end{proposition}
\begin{proof}
Evaluation on the characters of $K(\D_i)$ identifies $F_{i*}(e_\rho)$
with the sum of the primitive idempotents of $\sigma\in\Pi_\rho$.
Taking its trace in~\eqref{eq:relative-forgetful-trace}
proves~\eqref{eq:image-codegrees}.
Multiplying by $\FPdim(\D_i)$ gives
$s_\rho=\FPdim(\D_i)k_\rho/n_\rho$.
Now use~\eqref{eq:image-dimension} and
$\FPdim(Z_\rho^{\E_i})=\FPdim(\E_i)/n_\rho$ to
obtain~\eqref{eq:image-classes}. Finally, the fibers of restriction
from $K(\C)$ to $K(\D_i)$ over the members of $\Pi_\rho$ partition
$\{\mu_j:Z_j\in\Omega_\rho\}$ by~\eqref{eq:app-master-20261102}.
Applying the extent formula~\eqref{explicit:80658146e8:eq0008} of Proposition~\ref{prop:extent} to
each of them gives~\eqref{eq:image-to-C}. Both $s_\rho$ and the index
$\FPdim(\C)/\FPdim(\D_i)$ are integers, which yields the divisibilities.
\end{proof}

\begin{corollary}[Divisibility in the two categories]\label{prop:Ei-divisibility}
Let $\C$ be a weakly integral braided fusion category over $\kk$. For the forgetful functor
$F\colon\Z(\C)\to\C$ with right adjoint $I$, write
$I(\one)=\bigoplus_j Z_j$ and $Y_j=F(Z_j)$. Fix $i$ and set
$\E_i=\langle Z_i\rangle\subseteq\Z(\C)$,
$\D_i=\langle Y_i\rangle\subseteq\C$, and $B_i=I(\one)\cap\E_i$.
Write $F_i\colon\E_i\to\D_i$ for the restriction of $F$.
For an algebra character $\rho$ of $K(\E_i)$ let $n_\rho$
be its formal codegree and set $k_\rho=\rho(B_i)$.
Let $\Omega_\rho$ consist of those $W\in\Irr(\Z(\C))$ whose
Verlinde character for the canonical modular structure restricts to
$\rho$ on $K(\E_i)$.
Let $\Pi_\rho$ consist of the algebra characters $\sigma$ of
$K(\D_i)$ whose pullback along $F_{i*}$ equals $\rho$.
For $W\in\Omega_\rho$ one has
\begin{equation}\label{eq:Ei-quotient}
 \frac{\FPdim(W)^2}{\FPdim(Z_\rho^{\E_i})}
 =\FPdim\bigl((WW^*)\cap(\E_i)'_{\Z(\C)}\bigr)\in\mathbb Z_{>0}.
\end{equation}
Let $\mu_j$ denote the algebra character of $K(\C)$ corresponding to
$Z_j$. If $\mu_j|_{K(\D_i)}=\sigma\in\Pi_\rho$, then
\begin{equation}\label{eq:two-divisors}
 \FPdim(Z_\rho^{\E_i})\mid\FPdim(Z_j)^2,\qquad
 \FPdim(Z_\sigma^{\D_i})\mid\FPdim(Z_j).
\end{equation}
Thus a prime absent from all class dimensions in $\C$ is absent
from all class dimensions in $\D_i$, and also from the dimensions
$\FPdim(Z_\rho^{\E_i})$ for $k_\rho>0$.
\end{corollary}
\begin{proof}
Apply the modular quotient formula~\eqref{eq:modular-quotient} of Theorem~\ref{thm:modular-quotient} to $\E_i\subseteq\Z(\C)$.
Under the stated restriction, $Z_j\in\Omega_\rho$,
so~\eqref{eq:Ei-quotient} applies with $W=Z_j$.
Apply the individual braided divisibility~\eqref{eq:braided-divisibility} of Theorem~\ref{thm:braided-divisibility} to $\D_i\subseteq\C$ for the
second divisibility. Every $\sigma$ extends to $K(\C)$ by Proposition~\ref{prop:extent},
and every supported $\rho$ has such a lift by Lemma~\ref{lem:fibres}.
\end{proof}

\subsection{Quantitative vanishing for a class image}

\begin{corollary}[Zero bounds in the center and its forgetful image]
\label{cor:relative}\label{cor:classzero}\label{cor:image-zero}
Let $\C$ be a weakly integral braided fusion category over $\kk$. For the forgetful functor
$F\colon\Z(\C)\to\C$ with right adjoint $I$, write
$I(\one)=\bigoplus_j Z_j$ and $Y_j=F(Z_j)$. Fix $i$ and set
$\E_i=\langle Z_i\rangle\subseteq\Z(\C)$,
$\D_i=\langle Y_i\rangle\subseteq\C$, and $B_i=I(\one)\cap\E_i$.
For an algebra character $\rho$ of $K(\E_i)$ let $n_\rho$
be its formal codegree and set $k_\rho=\rho(B_i)$.
Let $\sigma$ range over the algebra characters of $K(\D_i)$, and
let $\mu_j$ be the algebra character of $K(\C)$ corresponding to $Z_j$.
Then
\begin{align}
 \sum_{\substack{\rho(Z_i)=0\\k_\rho>0}}\frac{k_\rho}{n_\rho}
   &\ge\frac{\FPdim(Z_i)-1}{\FPdim(\D_i)},\label{eq:zerobound}\\
 \sum_{\sigma(Y_i)=0}\FPdim(Z_\sigma^{\D_i})
   &\ge \FPdim(Z_i)-1,\label{eq:image-zero}\\
 \sum_{\mu_j(Y_i)=0}\FPdim(Z_j)
   &=\frac{\FPdim(\C)}{\FPdim(\D_i)}
       \sum_{\sigma(Y_i)=0}\FPdim(Z_\sigma^{\D_i})
     \ge \frac{\FPdim(\C)}{\FPdim(\D_i)}(\FPdim(Z_i)-1).
     \label{eq:global-class-zero-bound}
\end{align}
In particular every $Z_i$ of dimension greater than one admits a
character $\mu_j$ of $K(\C)$ with $\mu_j(Y_i)=0$.
\end{corollary}
\begin{proof}
Apply Theorem~\ref{thm:general-relative} in $\E_i$ with $B=B_i$ and
$Z=Z_i$. The multiplicities are both one and $\langle B_i\rangle=\E_i$.
Lemma~\ref{lem:fibres} supplies the nonnegative integral weights and
positivity on their support. Equation~\eqref{eq:image-dimension}
therefore gives~\eqref{eq:zerobound}.
For $\sigma\in\Pi_\rho$, one has $\sigma(Y_i)=\rho(Z_i)$.
Summing~\eqref{eq:image-codegrees} of Proposition~\ref{prop:image-fibres} over the zero fibers and
multiplying by $\FPdim(\D_i)$ gives~\eqref{eq:image-zero}.
Summing~\eqref{eq:image-to-C} over the same fibers gives
\eqref{eq:global-class-zero-bound}.
\end{proof}

\section*{Acknowledgments}
The author used the assistance of ChatGPT Pro (OpenAI) in developing proof strategies, examining examples, and drafting and revising portions of this paper. The author takes full responsibility for the statements, proofs, references, and final text.


\end{document}